\documentclass[12pt,fleqn, a4]{article}

\usepackage[french]{babel}
\usepackage{inputenc} 
\usepackage[T1]{fontenc}

\usepackage{color}
\usepackage{amssymb}
\usepackage{amsmath}
\usepackage{dsfont}
\usepackage{fancyheadings}

\newcommand{\dcb}{\begin{array}{lll}}
\newcommand{\dce}{\end{array}}
\newcommand{\ebe}{\begin{enumerate}\setlength{\baselineskip}{17pt}\setlength{\parskip}{0pt}}
\newcommand{\dbe}{\end{enumerate}}

\newcommand{\ibegin}{\begin{itemize}\setlength{\baselineskip}{19pt}\setlength{\parskip}{7pt}}
\newcommand{\iend}{\end{itemize}}
\newcommand{\ok}{\rule{4pt}{6pt}}%

\newtheorem{Theorem}{Theorem}[section]
\newtheorem {Cor}[Theorem]{Corollary}
\newtheorem {definition}[Theorem]{Definition}
\newtheorem {pro}[Theorem]{Proposition}
\newtheorem {Lemma}[Theorem]{Lemma}
\newtheorem {rem}[Theorem]{Remark}
\newtheorem {assumption}[Theorem]{Assumption}

\newcommand {\bd}{\begin{definition}}
\newcommand {\ed}{\end{definition}}
\newcommand {\bl}{\begin{Lemma}}
\newcommand {\el}{\end{Lemma}}

\newcommand {\bpro}{\begin{pro}}
\newcommand {\epro}{\end{pro}}
\newcommand {\bcor}{\begin{Cor}}
\newcommand {\ecor}{\end{Cor}}
\newcommand {\brem }{\begin{rem} \rm }
\newcommand {\erem }{\end{rem}}
\newcommand{\bethe}{\begin{Theorem}}
\newcommand{\ethe}{\end{Theorem}}
\newcommand {\bassumption}{\begin{assumption}}
\newcommand {\eassumption}{\end{assumption}}

\def\proof{\noindent {\bf Proof. $\, $}}

\def \ind{1\!\!1\!}

\def\cro#1{\langle #1\rangle}

\begin{document}

\begin{center}
\textbf{\large 
\bfseries
Multi-dimensional BSDEs whose terminal values are  bounded and have bounded Malliavin derivatives}\footnote{A third arXiv version with some improvements and new results}

Shiqi Song\footnote{This research has benefited from the support of the ``Chair Markets in Transition'',
F\'ed\'eration Bancaire Fran\c caise, and of the ANR project 11-LABX-0019.}\\
Laboratoire Analyse et Probabilit\'{e}s\\ 
Universit\'{e} d'Evry Val D'Essonne, France
\end{center}

\

\begin{quote}
\textbf{Abstract.} \itshape
We consider a class of multi-dimensional BSDEs on a finite time horizon (containing in particular Lipschitzian-quadratic BSDEs), whose terminal values are bounded as well as their corresponding Malliavin derivatives. We prove two results. The first one is an exponential integrability condition which determines when a BSDE in this class has a solution up to a given time horizon. In the second result, via an ordinary differential equation, we compute a minimum horizon up to which any BSDE of this class has a solution. The combination of these two results leads to a new scheme to solve quadratic BSDEs.
\end{quote}

\

\

It is not easy for a novice to follow the recent development in the computation of multi-dimensional non Lipschitzian BSDEs. The classical scheme based on the \textit{\`a-priori} estimations is no more valid. Instead, a big variety of techniques take place, which are valid only for specific divisions of BSDEs. The ideal is to have a global and centralized vision of the development, with a \texttt{"}map\texttt{"} of all the techniques, organized according to the way the techniques are involved, according to the reasons for which they are introduced, according to the consequences they yield. The present paper is a consequence of the effort to constitute such a map.

Many aspects can be used to classify the computation techniques. We can mention the Markov property (cf. \cite{HR, XZ}), the Girsanov transformation technique (cf. \cite{CN, HT}), the smallness techniques (cf. \cite{JKL, T}), the sliceability (cf. \cite{HR}) and the various others such as exponential transformation, Malliavin derivative, Lyapunov function, etc..  
Trying to establish logical connections of the various techniques, we notice somethings missing from the literature. A first point is the exponential integrability condition, which has been largely discussed in the case of one-dimensional quadratic BSDEs (cf. \cite{BEK, BH, MW}), but is absent in the multi-dimensional case. A second missing point is about local solutions. In fact, like ordinary differential equations, multi-dimensional non Lipschitzian BSDEs can have local solutions (with a possible explosion time (cf. \cite{FR})). The existence of local solutions can be very helpful in the computation of multi-dimensional BSDEs. From the very beginning of the BSDE literature, people already has the idea to solve a BSDE in computing it recursively on each of the small sub-intervals $[T-(k+1)h, T-h]$, whenever possible. See \cite[proof of Theorem 5.1]{EKPQ}. However, despite its popularity, the local existence property seems not to have been systematically studied (in computing the local solution and the explosion time).

In fact, we expect to organize the computations of multi-dimensional non Lipschitzian BSDEs around some general principles, just like the computations of Lipschitzian BSDEs are based on the \textit{\`a-priori} estimations. We believe that the exponential integrability and the local solvability constitute two such fundamental principles. 
In the present paper, we try to prove these two principles for the class of multi-dimensional BSDEs whose terminal values are bounded as well as their corresponding Malliavin derivatives (cf. Section \ref{theorem-section} for the exact definition of this class). In Theorem \ref{exp-integrability} we show that, effectively, the resolution of a BSDE in that class is characterized by an exponential integrability condition. In Theorem \ref{solution-extend}, we show that the BSDEs in that class have always non trivial local solutions, and we compute a minimum horizon up to which bounded solutions exist for all the BSDEs of that class. The proofs of Theorem \ref{exp-integrability} and of Theorem \ref{solution-extend} necessitate a property of linear BSDEs with unbounded coefficients. It seems that, despite its elementary natural, no complete result exists in the literature for the resolution of such BSDEs. To fill the gap, we will prove Lemma \ref{linear-existence} together with Corollary \ref{varrho-estimation}.

Section \ref{definition-lyapunov} is dedicated to illustrating how the two general principles, i.e. the exponential integrability condition and the local solvability, can be applied to solve BSDEs in concrete situations. Already, the exponential integrability condition gives a nice explanation why the sliceability can be so useful in BSDE computations, while the smallness technique is a direct consequence of the formula in Theorem \ref{solution-extend} of the explosion times. More specifically, we will illustrate the two principles in the well-known situation of sub-quadratic BSDEs. After that, we will prove a new result on quadratic BSDEs whose growth rates do not exceed the power $1+\alpha$, for an $0\leq \alpha<1$. By the way, we take this opportunity to examine another technique: the Lyapunov function. This technique (combined with the Gaussian estimates) has proved its efficiency in case of Markovian quadratic BSDEs (cf. \cite{XZ}). We would like to do the same under the Malliavin differentiability.

The paper is not organized in a chronological order. We have prefered to present the main results in Section \ref{theorem-section} and the applications in Section \ref{definition-lyapunov}, and write their proofs with elements which will only be considered in latter sections. In fact, we believe profitable to separate, in these proofs, the general reasoning, valid in various other situations, from those specifical to the proofs. Section \ref{LBSDE} is a typical example of general reasoning, where the results on multi-dimensional linear BSDEs with unbounded coefficients has clearly an independent interest. Likewise, from Section \ref{epsilon-neighbor} to Section \ref{L-S}, one can find various estimates which are valid in general.

The present paper is issued from a review of the existing literature on quadratic BSDEs. Some computations are classical. However, for a better readability, the computations of the paper will be presented in a self-contained manner, and organized according to our own vision of the BSDEs. The literature on the theory of BSDEs is large. We can not give here a complete presentation of it. One can instead find a good r\' esum\' e in the recent papers \cite{HR, JKL} on that part of the literature about the multi-dimensional quadratic BSDEs.

\section{Preliminary}

Let $B$ be a $n$-dimensional ($n\in\mathbb{N}^*$) Brownian motion living on a probability space $(\Omega, \mathcal{A}, \mathbb{P})$, where $(\Omega,\mathcal{A})$ is a measurable space, $\mathbb{P}$ is a probability measure on the $\sigma$-algebra $\mathcal{A}$. Let $\mathbb{F}=(\mathcal{F}_{t})_{t\in\mathbb{R}_{+}}$ be the natural filtration of the Brownian motion $B$, augmented by the usual conditions.

\subsection{Backward stochastic differential equation}

A $d$-dimensional ($d\in\mathbb{N}^*$) backward stochastic differential equation (BSDE in short) with parameters $[T, \xi, f]$ is an expression of following type:
\begin{equation}\label{theBSDE}
\left\{
\dcb
dY_{t} = -f(t, Y_{t}, Z_{t}) dt + Z_{t} dB_{t},\\
\\
Y_{T}=\xi.
\dce
\right.
\end{equation}
Here 

\noindent\textbf{-}
$T>0$ denotes a positive real number, that will be called the terminal time of the BSDE, 

\noindent\textbf{-}
$\xi=(\xi_{i})_{1\leq i\leq d}$ denotes a $d$-dimensional $\mathcal{F}_{T}$ measurable vector of $d$ real random variables, that will be called the terminal value of the BSDE,  

\noindent\textbf{-}
$f(t,y,z)$ denotes a $d$-dimensional parametered predictable process, i.e., a function from the space $\Omega\times \mathbb{R}_{+}\times \mathbb{R}^d\times \mathbb{R}^{d\times n}$ into $\mathbb{R}^d$, which is measurable with respect to the $\sigma$-algebra $\mathcal{P}(\mathbb{F})\otimes\mathcal{B}(\mathbb{R}^d)\otimes\mathcal{B}(\mathbb{R}^{d\times n})$ (where $\mathcal{P}(\mathbb{F})$ denotes the predictable $\sigma$-algebra of $\mathbb{F}$), that will be called the driver of the BSDE, while $f(s,0,0)$ will be called the rout of the driver,

\noindent\textbf{-}
$Y$ will be called the value process of the BSDE,

\noindent\textbf{-}
$Z$ will be called the martingale coefficient process.

A solution of the BSDE$[T,\xi, f]$(\ref{theBSDE}) is a pair $(Y,Z)$, where $Y$ is a $d$-dimensional $(\mathbb{F}, \mathbb{P})$ semimartingale and $Z=(Z_{i,j})_{1\leq i\leq d, 1\leq j\leq n}$ is a $d\times n$-dimensional matrix valued predictable process, such that 
\begin{equation}
\int_{0}^T |f(i, s, Y_{s}, Z_{s})| ds<\infty, \
\int_{0}^T Z_{i,j, s}^2\ ds < \infty, \ \forall  1\leq i\leq d, 1\leq j\leq n,
\end{equation}
and
\begin{equation}\label{BSDE-integral}
Y_{i,t}
=
\xi_{i} + \int_{t}^T f(i, s, Y_{s}, Z_{s})ds - \int_{t}^T Z_{i,s} dB_{s}, \ \forall 0\leq t\leq T, \ 1\leq i\leq d,
\end{equation}
where $Y_{i}$ denotes the $i^e$ component of the process $Y$, $f(i)$ denotes the $i^e$ component of the parametered process $f$, $Z_{i}$ denotes the $i^e$ row of $Z$. Obviously, the notion of solution of the BSDE$[T,\xi, f]$(\ref{theBSDE}) can be defined on any time interval $[b,T]$ for $0\leq b\leq T$.

\subsection{Index convention}

A family of indexed objects is considered as a function defined on the index sets (for example, $\phi_{i,i',a,t}$, $h_{i,i',j,t}$, $Z_{i,j,t}$, $D_{j,\theta}Y_{i,t}$). As usual, the notations $i,j$ indicate integer indices, while the others like $t, \theta$ indicate positive real numbers. When all indices are present in an expression, it indicates the value of the corresponding object at indices. When some indices are omitted, for example $D_{\theta}Y_{t}$, it indicates, as in the example, the sub-family $(D_{j,\theta}Y_{i,t})_{i,j}$. There will no confusion possible, when the expressions $Z_{t}=(Z_{i,j,t})_{i,j}$ and $Z_{i}=(Z_{i,j,t})_{j,t}$ are used, because they are distinguished by the different indices $t$ and $i$ which vary in different definition sets. When two indices run in the same index set, as $(i,i')$ in the example $h_{i,i',j,t}$ or $(a,t)$ in $\phi_{i,i',a,t}$, the expression $h_{i,j,t}$ or $\phi_{a}$ is only used with respect to the first of the two indices: $h_{i,j,t}=(h_{i,i',j,t})_{i'}$ and $\phi_{a}=(\phi_{i,i',a,t})_{i,i',t}$.

The notation $N$ will represent a positive integer and will be reserved to denote the index of an approximation sequence. The consideration of the definition sets for the indices allows us to exclude the confusions between, for example, the expression $Y_{t}$ denoting the value at $t\in\mathbb{R}_{+}$ of the limit process $Y$, and the expression $Y_{N}$ denoting an element in an approximation sequence ($Y_{N}=(Y_{N,t})_{t\geq 0} \rightarrow Y$).  

The notations such as $i,i'$ represent always indices in the set $\{1,\ldots,d\}$. The notations such as $j,j'$ indicate always indices in the set $\{1,\ldots,n\}$. The notations such as $s,t,\theta, a,b$ are used to present positive real numbers. In particular, the notation $\theta$ is especially used in relation with a Malliavin derivative.

\subsection{Matrix notations}\label{matrix}

For any real vector $u=(u_{1}, \ldots, u_{k})$ ($k\in\mathbb{N}^*$), we define $
|u| = \sqrt{\sum_{i=1}^ku_{i}^2},
$
and for a matrix $z\in\mathbb{R}^{k\times k'}$, we define $
\|z\| = \sqrt{\sum_{a=1}^k \sum_{b=1}^{k'} z_{a,b}^2}=\sqrt{\sum_{a=1}^k |z_{a}|^2}.
$
When two matrix $g,z$ of suitable formats are put side by side: $gz$, it represents the matrix product. We may need scalar product between two matrix $z,z'$ of same format, denoted by $\cro{z \boldsymbol{|} z'}=\sum_{a=1}^k \sum_{b=1}^{k'} z_{a,b}z_{a,b}'$. If necessary, a number or a vector can be considered as a matrix so that the notations $|\cdot|$ and $\|\cdot\|$ are not exclusionary.

Exception is made in the equations (\ref{theBSDE-malliavin}) and (\ref{theBSDE-malliavin-*}), where $\partial_{z}f(t, Y_{t}, Z_{t})D_{j,\theta}Z_{t}$ denotes actually $\cro{\partial_{z}f(t, Y_{t}, Z_{t})\boldsymbol{|}D_{j,\theta}Z_{t}}=\sum_{i'=1}^d \sum_{j'=1}^{n}\partial_{z_{i',j'}}f(t, Y_{t}, Z_{t})D_{j,\theta}Z_{i',j',t}$ (cf. the notation system of \cite[Proposition 2.4]{EKPQ}, or the notation system of \cite{HR}). This exception is made for us to write exactly the same formulas as in \cite{EKPQ}.

A family of objects can be indexed by three integer indices as in the example $h_{i,i',j,t}$. We can have different matrix issued from this family: $h_{j,t}=(h_{i,i',j,t})_{1\leq i,i'\leq d}$ or $h_{i,,t}=(h_{i,i',j,t})_{1\leq i'\leq d, 1\leq j\leq n}$. It is to notice that, in this example, we have $
\sum_{j=1}^n\|h_{j,t}\|^2 = \sum_{i=1}^d\|h_{i,t}\|^2.
$

\subsection{Growth rate and increment condition}\label{growth-increment-conditions}

In this paper we always consider the BSDE(\ref{theBSDE}) with a driver $f$ which is Lipschitzian in $y$. This means that, for every $0\leq t\leq T$, there is a positive process $B\geq 0$ (independent of $y$ and of $z$) such that, for $1\leq i\leq d$, for $y, y'$ in $\mathbb{R}^{d}$, for $z\in\mathbb{R}^{d\times n}$,
\begin{equation}\label{lipsch-driver-in-y}
|f(i,t,y',z) - f(i,t,y,z)| \leq B_{t} \sum_{i'=1}^d|y'_{i'}-y_{i'}|.
\end{equation} 
The processes $B$ will be called the Lipschitzian index in $y$. If we can take $B$ a positive constant: $B=\beta>0$, we say that $f$ is uniformly Lipschitzian in $y$ with index $\beta$.

In the same way, we can define the Lipschitzian condition on the driver $f$ in the variable $z$. But, the present paper is especially a study of the BSDE(\ref{theBSDE}), when the increment condition of the driver $f$ in the variable $z$ takes various different forms.

We introduce therefore the notion of $\varkappa$-rate condition in $z$. Let $\varkappa$ be a (deterministic) function defined on $\mathbb{R}_{+}$ which is continuous, nonnegative, increasing and locally Lipschitzian. We say that $f$ has an increment rate in $z$ uniformly less than $\varkappa$, if, for $1\leq i\leq d$, for $y\in\mathbb{R}^{d}$, for $z, z'\in\mathbb{R}^{d\times n}$,
\begin{equation}\label{kappa-increment}
|f(i,t,y,z) - f(i,t,y,z')|\leq \varkappa(\|z\|+\|z'\|)\|z-z'\|.
\end{equation}
We will call the function $\varkappa$ a rate function in $z$.

When $\varkappa$ can be taken to be a positive constant $\eta>0$, we have the uniform Lipschitzian condition in $z$. When $\varkappa$ is taken to be an affine function with positive coefficients, we have a quadratic increment condition in $z$.

Notice that, in the BSDE literature, there is a distinction between the increment condition such as the conditions (\ref{lipsch-driver-in-y}) or (\ref{kappa-increment}), from the growth condition (cf. \cite{bahlali} where the growth condition is used distinctly). For example, we say that the driver has an affine growth rate, if, for $1\leq i\leq d$, for $y$ in $\mathbb{R}^{d}$, for $z\in\mathbb{R}^{d\times n}$,
\begin{equation}\label{growth-affine-in-y}
|f(i,t,y,z)| \leq A_{t}+B_{t} \sum_{i'=1}^d|y_{i'}| + C_{t}\|z\|,
\end{equation} 
where $A,B,C$ are positive processes independent of $y$ and $z$, that we call growth indices. We also can say that the driver $f$ has a linear growth rate in $y$ and a $\varkappa$ growth rate in $z$, if
\begin{equation}\label{growth-rate}
|f(i,t,y,z)| \leq A_{t}+\beta \sum_{i'=1}^d|y_{i'}| + \varkappa(\|z\|)\|z\|.
\end{equation} 
The increment conditions are fundamentally linked with the (local) existence problem of the BSDE(\ref{theBSDE}), while the growth rate conditions are used to control the integrability of the solutions of the BSDE(\ref{theBSDE}). In this paper, however, we will not study the specificity of the growth rate condition (except in Section \ref{alpha-growth}). We will content ourselves with the growth rate condition induced by the increment condition.

\subsection{Spaces of processes}\label{the-spaces}

For $0\leq b\leq T<\infty$, we consider the following spaces of processes.
$$
\dcb
\mathcal{S}_{\infty}[b,T]
=
\{X=(X_{t})_{b\leq t\leq T}: \mbox{$X$ is continuous}, X_{t}\in\mathbb{R}^d, X_{t}\in\mathcal{F}_{t}, \|\sup_{b\leq t\leq T}\sum_{i=1}^d|X_{i,s}|\|_{\infty}<\infty\},\\
\\
\mathcal{Z}_{\mbox{\tiny BMO}}[b,T]
=
\{X=(X_{t})_{b\leq t\leq T}: X_{t}\in\mathbb{R}^{d\times n} \mbox{ predictable}, \sup_{b\leq \sigma\leq T}\|\mathbb{E}[\int_{\sigma}^T \|X_{s}\|^2 ds\ |\mathcal{F}_{\sigma}] \|_{\infty}<\infty\},
\dce
$$
where $\sigma$ runs over the class of all stopping times. 
The elements in $\mathcal{Z}_{\mbox{\tiny BMO}}[0,T]$ are naturally associated with the BMO martingales. See \cite[Chapter X]{HWY}  and \cite{K}. In fact, we will use the notation (the BMO-norm):$$
\dcb
\|X\|_{\star}:=\sqrt{\sup_{\sigma\leq T}\|\mathbb{E}[\int_{\sigma}^T \|X_{s}\|^2 ds\ |\mathcal{F}_{\sigma}]\|_{\infty}}
\dce
$$
for any predictable matrix valued process $X$.

\bl
If $\int_{0}^T \|X_{s}\|^2 ds$ is integrable, we have$$
\dcb
\|X\|_{\star}
=
\sqrt{\sup_{0\leq t\leq T}\|\mathbb{E}[\int_{t}^T \|X_{s}\|^2 ds\ |\mathcal{F}_{t}]\|_{\infty}}.
\dce
$$
\el

\proof
Let $$
\dcb
c=\sqrt{\sup_{0\leq t\leq T}\|\mathbb{E}[\int_{t}^T \|X_{s}\|^2 ds\ |\mathcal{F}_{t}]\|_{\infty}}.
\dce
$$
Clearly $c \leq \|X\|_{\star}$. For any discrete stopping time $\sigma$, there exists a sequence of positive real numbers $(a_{k})_{k\geq 1}$ such that $$
\sigma = \sum_{k=1}^\infty a_{k}\ind_{\{\sigma = a_{k}\}}+\ind_{\{\sigma = \infty\}}.
$$
We have 
$$
\dcb
&&\mathbb{E}[\int_{\sigma}^T \|X_{s}\|^2 ds\ |\mathcal{F}_{\sigma}]
=
\sum_{k=1}^\infty\mathbb{E}[\int_{\sigma}^T \|X_{s}\|^2 ds\ |\mathcal{F}_{\sigma}]\ind_{\{\sigma = a_{k}\}}\\

&=&
\sum_{k=1}^\infty\mathbb{E}[\int_{a_{k}}^T \|X_{s}\|^2 ds\ |\mathcal{F}_{a_{k}}]\ind_{\{\sigma = a_{k}\}}
\leq
\sum_{k=1}^\infty c^2\ \ind_{\{\sigma = a_{k}\}} \leq c^2. 
\dce
$$
This inequality combined with \cite[Theorem 3.7]{HWY} and \cite[Chapter II, Corollary 2.4]{RY} implies $\|X\|_{\star}^2\leq c^2$. \ok

\subsection{Malliavin derivative}

We will employ the notion of Malliavin calculus. We refer to \cite{N} for a detailed account of the theory of Malliavin calculus. The Malliavin calculus here is effected with respect to the $n$-dimensional Brownian motion $B$. We make use of the usual notions of Malliavin calculus. In particular, the Malliavin derivative is denoted by $D$ and 
$$
\dcb
\mathbb{D}^{1,2}
&=\{\zeta: &\mbox{$\zeta$ is a real square integrable random variable,} \\
&&\mbox{whose Malliavin derivative exists and is square integrable}\}.
\dce
$$
In \cite[Proposition 5.3]{EKPQ}, the Malliavin derivative of the BSDE$[T,\xi,f]$(\ref{theBSDE}) is studied. We recall here the conditions of \cite[Proposition 5.3]{EKPQ} on the driver $f$ for the Malliavin differentiability of the BSDE$[T,\xi,f]$(\ref{theBSDE}). We divide the conditions into two parts.
\begin{equation}\label{CR-C1}
\dcb
\mbox{\textbf{-}} & \mbox{the driver $f(t,y,z)$ has first derivatives in $y,z$ and the derivatives}\\
& \mbox{are continuous in $y,z$,}\\

\mbox{\textbf{-}} & \mbox{for every $i, t, y,z$, $f(i,t,y,z)$ is Malliavin differentiable in $\mathbb{D}^{1,2}$,}\\

\mbox{\textbf{-}} 
& \mbox{the function $D_{j,\theta}f(\omega,i,t, y,z)$ is measurable with respect to all parameters,}\\ 
 
\mbox{\textbf{-}} 
& \mbox{for every $j,\theta, i, y,z$, the process $D_{j,\theta}f(i,\cdot, y,z)$ is predictable,}\\ 
 
\mbox{\textbf{-}}
& \mbox{for every $y,z$, $\mathbb{E}[\int_{0}^T\int_{0}^T\|D_{\theta}f(s, y,z)\|^2ds\ d\theta]<\infty$,}\\

\mbox{\textbf{-}}
& \mbox{for every $j,\theta, i,t$, the random function $D_{j,\theta}f(i,t, y,z)$ is Lipschitzian in $y,z$}\\
& \mbox{with common index $K(\theta,t)$ such that $\int_{0}^T\mathbb{E}[(\int_{0}^T K(\theta,s)^2 ds)^2]d\theta <\infty$}.
\dce
\end{equation}
To have the Malliavin differentiability, \cite[Proposition 5.3]{EKPQ} supposes, in addition of (\ref{CR-C1}), 
\begin{equation}\label{CR-C1-complement}
\dcb
\mbox{\textbf{-}} &
\xi\in\mathbf{L}^4,\\

\mbox{\textbf{-}} & \mbox{the derivatives of $f(i, t,y,z)$ in $y,z$ are bounded,}\\

\mbox{\textbf{-}} &
\sum_{i=1}^d \mathbb{E}[(\int_{0}^{T} f(i, s, 0, 0)^2ds)^2]<\infty,\\

\mbox{\textbf{-}}
& \mbox{$\mathbb{E}[\int_{0}^T\int_{0}^T\|D_{\theta}f(s, Y_{s},Z_{s})\|^2ds\ d\theta]<\infty$,}
\dce
\end{equation}
where $(Y,Z)$ is the solution of the BSDE$[T,\xi,f]$(\ref{theBSDE}), which exists, because the boundedness of the derivatives of $f$ implies that the driver $f$ is Lipschitzian.

We note that the conditions (\ref{CR-C1}) are purely differentiability conditions of the driver $f$, while the conditions (\ref{CR-C1-complement}) are mixed. In the present paper, the conditions (\ref{CR-C1}) will be applied as independent conditions. As for the conditions (\ref{CR-C1-complement}), they will be provided via other packages of conditions. Below, we will call the combination of conditions (\ref{CR-C1}) and (\ref{CR-C1-complement}) the strong version of (\ref{CR-C1}).

\subsection{ Gronwall's type inequalities}\label{gronwall}

In the case of a Lipschitzian BSDE(\ref{theBSDE}), a basic method to control the level of the integrability of the solution is to apply the Gronwall inequality. The following lemma formulates the Gronwall inequality in a particular form which makes its applications in the present paper easier.

\bl\label{gronwall-ineq}
For real number $c\geq 0$, for non negative functions $v(t), u(t), h(t)$ defined on $0\leq t\leq T$, if $u(t), h(t), h(t)v(t)$ are integrable on $0\leq t\leq T$, if the inequality $$
v(t) \leq c + \int_{t}^T u(s) ds + \int_{t}^T h(s) v(s) ds,\ 0\leq t\leq T,
$$
holds, we have $
v(t) \leq e^{\int_{t}^T h(s) ds}(c + \int_{t}^T u(s) ds ).
$
For non negative integrable random variable $\zeta$, for non negative $\mathbb{F}$-adapted measurable process $V$, for non negative measurable process $U$, if $\mathbb{E}[\int_{0}^T U_{s} + h(s)V(s)ds]$$<\infty$ and 
$$
V_{t} \leq \mathbb{E}[\zeta+\int_{t}^T U_{s} ds + \int_{t}^T h(s) V_{s} ds\ |\mathcal{F}_{t}], \ 0\leq t\leq T,
$$
we have $
V_{t} \leq e^{\int_{t}^T h(s) ds}\mathbb{E}[\zeta + \int_{t}^T U_{s} ds\ |\mathcal{F}_{t}].
$
\el

\proof
The first part of the lemma follows from \cite[Corollary 6.62]{PR}. To prove the second part of the lemma, we can write, for $0\leq a\leq t\leq T$,
$$
\dcb
\mathbb{E}[V_{t}|\mathcal{F}_{a}]
&\leq&
\mathbb{E}[\zeta + \int_{t}^T U_{s}ds|\mathcal{F}_{a}] + 
\int_{t}^T h_{s} \mathbb{E}[V_{s} |\mathcal{F}_{a}] ds.
\dce
$$
The point $a$ being fixed, we can write this inequality with a regular conditional probability kernel of $\mathbb{E}[\ \cdot\ |\mathcal{F}_{a}]$. The first part of the lemma is applicable on the time interval $[a,T]$. $$
\dcb
\mathbb{E}[V_{t}|\mathcal{F}_{a}]
\leq
e^{\int_{t}^Th_{s}ds}\mathbb{E}[\zeta + \int_{t}^T U_{s} ds|\mathcal{F}_{a}],\ a\leq t\leq T.
\dce
$$
Take $a=t$ to finish the proof of the lemma. \ok

\subsection{Two fundamental identities of the BSDE(\ref{theBSDE})}\label{fundamental-identities}

Suppose a solution $(Y,Z)$ of the BSDE(\ref{theBSDE}) exists. Then, for $1\leq i\leq d$, $0\leq t\leq T$,
$$
\dcb
d|Y_{i,t}|
=
-\mbox{sig}(Y_{i,t})f(i,t, Y_{t}, Z_{t}) dt + \mbox{sig}(Y_{i,t})Z_{i,t}dB_{t} + dL_{i,t},\\

d Y_{i,t}^2
=
- 2Y_{i,t}f(i,t,Y_{t},Z_{t})dt +  2 Y_{i,t}Z_{i,t}dB_{t} +  |Z_{i,t}|^2 dt,
\dce
$$
where $L_{i}$ denotes the local time of $Y_{i}$ at zero (cf. \cite[Chapter VI]{RY}). If $\xi$ is integrable, if\\ $\mathbb{E}[\int_{0}^T \sum_{i=1}^d |f(i,s,Y_{s},Z_{t})|ds]<\infty$, if $Z\!\centerdot\! B$ is a true martingale, we can write
\begin{equation}\label{ine-y}
\dcb
&&\mathbb{E}[|Y_{i,t}| +\int_{t}^T dL_{i,s}|\mathcal{F}_{a}]
= 
\mathbb{E}[|\xi_{i}|+\int_{t}^T \mbox{sig}(Y_{i,t})f(i,s, Y_{s},Z_{s}) ds|\mathcal{F}_{a}],\ 0\leq a\leq t\leq T.
\dce
\end{equation}
If moreover the value process $Y$ is bounded (so that $\xi$ is bounded), 
\begin{equation}\label{ine-z}
\dcb
&&\mathbb{E}[Y_{i,t}^2 + \int_{t}^T  |Z_{i,s}|^2 ds\ |\mathcal{F}_{a}]
=
\mathbb{E}[\xi_{i}^2+\int_{t}^T 2Y_{i,s}f(i,s,Y_{s},Z_{s})ds\ |\mathcal{F}_{a}].
\dce
\end{equation}

\subsection{The stopping functions $r_{v}(x)$ and driver approximations}\label{approximation-drivers}

For any driver $f$, we consider the following approximation sequence. We introduce first the stopping functions $r_{v}(x)$. Here $v$ denotes a positive real number $v\in\mathbb{R}^*$ and $r_{v}(x), x\in\mathbb{R},$ denotes a non decreasing odd $C_{b}^2$-function such that $|r_{v}'(x)|\leq 1$, $r_{v}(x)=x$ for $-v\leq x\leq v$ and $r_{v}(x) = -v-1$ for $x\leq -v-4$ and $r_{v}(x) = v+1$ for $x\geq v+4$. Note that, for this function,$$
|r_{v}(x)| = r_{v}(|x|) \leq |x|,\ \ |r_{v}(x') - r_{v}(x)|\leq |x'-x|.
$$
For a matrix $z=(z_{i,j})\in\mathbb{R}^{d\times n}$, we define $
\mathbf{r}_{v}(z)
$
to be the matrix of $r_{v}(z_{i,j})$. We can verify the inequalities:$$
\|\mathbf{r}_{v}(z)\|\leq \|z\|,\ \ \|\mathbf{r}_{v}(z')-\mathbf{r}_{v}(z)\|\leq \|z'-z\|.
$$
We introduce the modified drivers$$
\bar{f}_{N}(t,y,z) = f(t, y,\mathbf{r}_{N}(z)),\ \ N\in\mathbb{N}^*,  y\in\mathbb{R}^d,  z\in\mathbb{R}^{d\times n}.
$$
It is to note that, if $f$ is a driver with Lipschitzian condition in $y$ and $\varkappa$-rate condition in $z$, for all $N>0$, $\bar{f}_{N}$ is a Lipschitzian driver in the two variables $y,z$.

Clearly, $f$ is the pointwise limit of $\bar{f}_{N}$. If $f(t,y,z)$ is Malliavin differentiable, we have $D\bar{f}_{N}(t,y,z)=Df(t,y,\mathbf{r}_{N}(z))$. Notice especially the following point.

\brem\label{strong-stationary}
If the BSDE$[T,\xi,\bar{f}_{N}]$(\ref{theBSDE}) with the driver $\bar{f}_{N}$ has a bounded solution $(Y_{N},Z_{N})$ such that $\sup_{0\leq s\leq T}\|Z_{N,s}\|_{\infty}<N$, the process $(Y,Z):=(Y_{N},Z_{N})$ is actually a bounded solution of the BSDE$[T,\xi,f]$(\ref{theBSDE}) with the original driver $f$, and vice versa.
\erem

\

\section{The main results}\label{theorem-section}

We introduce the following package of conditions:
\begin{equation}\label{package1}
\dcb
\mbox{\textbf{-}} & \mbox{$\xi\in\mathcal{F}_{T}$, $\|\xi\|_{\infty}<\infty$, $\xi\in\mathbb{D}^{1,2}$ and $\sup_{1\leq j\leq n}\sup_{1\leq i\leq d}\sup_{0\leq \theta<\infty}\|D_{j,\theta}\xi_{i}\|_{\infty}<\infty$,}\\

\mbox{\textbf{-}} 
& \mbox{the driver $f(i,t,y,z)$ is uniformly Lipschitzian in $y$ with index $\beta$ and }\\
&\mbox{has an increment rate in $z$ less than $\varkappa$, while its root is bounded by a constant $\upsilon$,}\\ 

\mbox{\textbf{-}} & \mbox{the conditions (\ref{CR-C1}),}\\
 
\mbox{\textbf{-}}
& \mbox{$\|D_{j,\theta}f(i,t,y,z)\|_{\infty}$ is uniformly bounded on any bounded set of $(j,\theta,i,t,y,z)$.}
\dce
\end{equation}
Notice that the above conditions are almost the minimum conditions that we have to assume to have the results of this paper. For example, the conditions (\ref{CR-C1}) is required to have the Malliavin differentiability of the BSDE. The boundedness of $f(t,0,0)$ is required, because we hope that the solution of the BSDE preserves the boundedness of $\xi$. The last condition in (\ref{package1}) holds, if $\|D_{j,\theta}f(i,t,y,z)\|_{\infty}$ is continuous in its parameters. Another sufficient condition for the last term in (\ref{package1}) is the following:
\begin{equation}\label{explicit-Df}
\dcb
& \mbox{$D_{j,\theta}f(i,t,y,z)$, $1\leq i\leq d, 1\leq j\leq n, 0\leq \theta\leq T,$ are uniformly Lipschitzian in $y$ }\\
&\mbox{with a common index $\hat{\beta}$ and have increment rates in $z$ less than a common rate }\\
&\mbox{function $\hat{\varkappa}$, while their roots $D_{j,\theta}f(i,t,0,0)$ are bounded by a common constant $\hat{\upsilon}$.}
\dce
\end{equation}
There is some redundancy among the above conditions. But we keep all conditions entire to make clear the reason why the conditions are introduced. There will be some advantages from this strategy. For example, the conditions (\ref{CR-C1}) can be replaced by any other sufficient conditions for Malliavin differentiability such as that in \cite{IMPR}.

\subsection{An exponential integrability condition}

In this section, we are interested in the solution-existence issue under $\varkappa$-rate condition. Ideally, we expect a \texttt{"}size level\texttt{"} to exist such that a BSDE under $\varkappa$-rate condition will have a solution, if and only if the size of $\varkappa$ does not exceed the level. We notice happily that it is exactly what happens under the conditions (\ref{package1}).

For any $N\in\mathbb{N}^*$ let $(Y_{N}, Z_{N})$ be the solution of the BSDE$[T, \xi, \bar{f}_{N}]$(\ref{theBSDE}). We introduce the conditions:
\begin{equation}\label{exp-boundedness}
\dcb
\sup_{N\in\mathbb{N}^*}\sup_{0\leq t\leq T}\|\ \mathbb{E}[e^{28 d^2\times n\int_{t}^T \varkappa(2\|Z_{N,s}\|)^2ds}|\mathcal{F}_{t}]\ \|_{\infty}<\infty,\\
\\
\sup_{N\in\mathbb{N}^*}\sup_{1\leq j\leq n}\sup_{0\leq \theta<\infty}\sup_{1\leq i\leq d} \sup_{0\leq t\leq T}\| \mathbb{E}[\ \int_{t}^T\  |D_{j,\theta}\bar{f}_{N}(i,s,Y_{N,s},Z_{N,s})|ds\ |\mathcal{F}_{t}] \|_{\infty}< \infty.
\dce
\end{equation}

\

\bethe\label{exp-integrability}
Suppose that the parameters of the BSDE$[T,\xi, f]$(\ref{theBSDE}) satisfy the conditions (\ref{package1}). 
Then, the BSDE$[T, \xi, f]$(\ref{theBSDE}) has a unique solution $(Y,Z)$ in $\mathcal{S}_{\infty}[0,T]\times \mathcal{Z}_{\mbox{\tiny BMO}}[0,T]$ with the property that 
\ebe
\item[$\centerdot$]
the Malliavin derivatives $D_{j,\theta}Y_{t}, 1\leq j\leq n, 0\leq \theta \leq T, 0\leq t\leq T$, exist and 
\item[$\centerdot$]
$D_{\theta}Y_{\theta}= Z_{\theta}$ for $0\leq \theta\leq T$, and
\item[$\centerdot$]
the processes $(D_{j,\theta}Y)_{1\leq j\leq n, 0\leq \theta \leq T}$ form a bounded family in the space $\mathcal{S}_{\infty}[0,T]$,
\dbe
if and only if the conditions (\ref{exp-boundedness}) hold. In this case, the Malliavin derivatives $D_{j,\theta}Y$ satisfy the linear 
BSDE$[T, \xi, f, j,\theta]$(\ref{theBSDE-malliavin-*}) on $[\theta,T]$.
\ethe

\proof
Notice at first that, under the conditions (\ref{package1}), $\bar{f}_{N}$ is Lipschitzian in the two variables with
$$
\dcb
|\partial_{y_{i'}}\bar{f}_{N}(i,t, Y_{N,t}, Z_{N,t})|\leq \beta,\ 1\leq i\leq d, 1\leq i'\leq d,\\

|\partial_{z_{i',j'}}\bar{f}_{N}(i,t, Y_{N,t}, Z_{N,t})| \leq  \varkappa(2\|\mathbf{r}_{N}(Z_{N,t})\|),\ 1\leq i\leq d, 1\leq i'\leq d,\ 1\leq j'\leq n.
\dce
$$ 
Consequently, for any $N\in \mathbb{N}^*$, $(Y_{N}, Z_{N})$ exists (cf. Proposition \ref{Lipschitzian-slice-collage}). Moreover, $Y_{N}$ is bounded on the time interval $[0,T]$ (cf. Proposition \ref{value-bound}) so that $D_{j,\theta}\bar{f}_{N}(i,s,Y_{N,s},Z_{N,s})$$=D_{j,\theta}f(i,s,Y_{N,s},\mathbf{r}_{N}(Z_{N,s}))$ is bounded, which implies that $Z_{N}$ is bounded (cf. Corollary \ref{MD-Z-bound}).

\textit{Necessity.} If the BSDE$[T,\xi, f]$(\ref{theBSDE}) has a solution $(Y,Z)$ as described in the theorem, the martingale coefficient process $Z$ will be bounded. So, as pointed out in Remark \ref{strong-stationary} of Section \ref{approximation-drivers}, for $N>\sup_{0\leq s\leq T}\|Z_{s}\|_\infty$, $(Y,Z)$ is a solution of the BSDE$[T, \xi, \bar{f}_{N}]$(\ref{theBSDE}). By uniqueness property, $(Y_{N}, Z_{N})$ coincides with $(Y,Z)$ for all $N$ big enough. The supremum in (\ref{exp-boundedness}) will be finite, because they concern only a finite number of $Z_{N}$, which are bounded by Corollary \ref{MD-Z-bound}.

\textit{Sufficiency.} Conversely, suppose the condition (\ref{exp-boundedness}) to hold. 
We apply Corollary \ref{varrho-estimation} to the Malliavin derivatives $D_{j,\theta}Y_{N}$ (cf. Section \ref{malliavin-derivative}) via the linear 
BSDE$[T, \xi, \bar{f}_{N}, j,\theta]$(\ref{theBSDE-malliavin-*}). We see that the processes $D_{j,\theta}Y_{N}$ are bounded on $[0,T]$ uniformly with respect to $j,\theta,N$. The identities $D_{\theta}Y_{N,\theta} = Z_{N,\theta}, 0\leq \theta\leq T$ (because \cite[Proposition 5.3]{EKPQ} is applicable to the BSDE$[T, \xi, \bar{f}_{N}]$(\ref{theBSDE})) implies that the processes $Z_{N}$ are also bounded uniformly with respect to $N$. Hence, according to Remark \ref{strong-stationary}, for any positive integer $N_{0}>\sup_{N\in\mathbb{N}^*}\sup_{0\leq s\leq T}\|Z_{N,s}\|_{\infty}$, the process $(Y_{N_{0}}, Z_{N_{0}})$ is in fact a solution of the BSDE$[T, \xi, f]$(\ref{theBSDE}) with the property described in the theorem. The solution is unique, because of the Lipschitzian property of the $\bar{f}_{N_{0}}$. \ok

\brem
Theorem \ref{exp-integrability} gives also an extension of \cite[Proposition 5.3]{EKPQ} to a family of non Lipschitzian BSDEs. 
\erem

Henceforth, for any $0\leq b\leq T$, we say that a solution $(Y,Z)$ in $\mathcal{S}_{\infty}[b,T]\times \mathcal{Z}_{\mbox{\tiny BMO}}[b,T]$ of BSDE$[T,\xi, f]$(\ref{theBSDE}) on $[b,T]$ satisfy the $\mathfrak{M}d$-property, if 
\ebe
\item[$\centerdot$]
the Malliavin derivatives $D_{j,\theta}Y_{t}, 1\leq j\leq n, 0\leq \theta \leq T, b\leq t\leq T$, exist, 
\item[$\centerdot$]
$D_{\theta}Y_{\theta}= Z_{\theta}$ for $b\leq \theta\leq T$, and
\item[$\centerdot$]
the processes $(D_{j,\theta}Y)_{1\leq j\leq n, 0\leq \theta \leq T}$ form a bounded family in the space $\mathcal{S}_{\infty}[b,T]$.
\dbe
Notice that, when $(Y,Z)$ satisfies the $\mathfrak{M}d$-property, by Remark \ref{strong-stationary}, $(Y,Z)$ is also the solution of a BSDE$[T, \xi, \bar{f}_{N}]$(\ref{theBSDE}) on $[b,T]$ for $N$ big enough, so that $D_\theta Y$ satisfies the BSDE$[T, \xi, f, \theta, j]$(\ref{theBSDE-malliavin}) on $[b\vee \theta,T]$.

\subsection{Estimation of the horizon up to which the BSDE has a solution}

We present now the second result. Consider the BSDE$[T,\xi, f]$(\ref{theBSDE}). If a solution $(Y,Z)$ of the BSDE$[T,\xi, f]$(\ref{theBSDE}) exists on $[t,T]$ as well as its corresponding Malliavin derivatives, we define, for $0\leq t\leq T$,
\begin{equation}\label{the-lambdas}
\dcb
\Lambda_{t}
=
\sup_{t\leq s\leq T}\sum_{i=1}^d\|Y_{i,s}\|_{\infty},\\

\hat{\Lambda}_{t}
=
\sup_{t\leq s\leq T}\sup_{0\leq \theta<\infty} \|D_{\theta}Y_{s}\|_{\infty}
=
\sup_{t\leq s\leq T}\sup_{0\leq \theta<\infty}\|\sqrt{\sum_{j=1}^n\sum_{i=1}^d (D_{j,\theta}Y_{i,s})^2}\|_{\infty}.
\dce
\end{equation}
Suppose the conditions (\ref{package1}) and (\ref{explicit-Df}). Let $\varkappa$ and $\hat{\varkappa}$ be the two rate functions given in (\ref{package1}) and (\ref{explicit-Df}). Define $\mathtt{g}^\circ(t,u), 0\leq t\leq T, 0\leq u<\infty,$ to be the real function: 
$$
\dcb
&&
u\times
\big(
d\times 2\beta 
+
d^2\times n\times \varkappa(2\sqrt{u})^2
\big)\\
&&
+\sqrt{u}2\sqrt{d\times n}\big( 
\hat{\upsilon}+\hat{\beta}e^{d \times \beta\times(T-t)}
\big(
\sum_{i=1}^d\|\xi_{i}\|_{\infty}
+
d\times (\upsilon+ \varkappa(\sqrt{u})
\sqrt{u})(T-t) \big) 
+ \hat{\varkappa}(\sqrt{u})\sqrt{u}
\big).
\dce
$$
Let $\mathtt{g}(t,u)$ be any continuously differentiable function, defined on $(t,u)\in[0,T]\times \mathbb{R}_{+}$, such that $\mathtt{g}(t,u)\geq \mathtt{g}^\circ(t,u)$.
Consider the ordinary differential equation, defined on $[0,T]$: 
\begin{equation}\label{u-equation}
\left\{
\dcb
u_{0}&=&\sup_{0\leq \theta<\infty}\|D_{\theta}\xi\|_{\infty}^2,\\
\\
\frac{du_{t}}{dt}&=&\mathtt{g}(T-t,u_{t}).
\dce
\right.
\end{equation}
Let $\alpha_{\mathtt{g}}$ be the supremum of the $0\leq a\leq T$ such that the differential equation (\ref{u-equation}) has a finite classical solution on $[0, a]$. Let $u_{\mathtt{g}}$ be the solution of the differential equation (\ref{u-equation}) on the time interval $[0,\alpha_{\mathtt{g}})$.

\bethe\label{solution-extend}
Suppose that the parameters of the BSDE$[T,\xi, f]$(\ref{theBSDE}) satisfy the conditions (\ref{package1}) and (\ref{explicit-Df}). Then, for any $T-\alpha_{\mathtt{g}}<a\leq T$, a unique solution $(Y,Z)$ of the BSDE$[T,\xi, f]$(\ref{theBSDE}) exists in $\mathcal{S}_{\infty}[a,T]\times \mathcal{Z}_{\mbox{\tiny BMO}}[a,T]$ on the time interval $[a, T]$ with $\mathfrak{M}d$-property. On the time interval $(T-\alpha_{\mathtt{g}}, T]$, the functions $\Lambda$ and $\hat{\Lambda}$ are well-defined and continuous and satisfy $$
\dcb
\hat{\Lambda}_{t}^2\leq u_{\mathtt{g}, T-t}\ \mbox{ and }\\

\Lambda_{t}
\leq
e^{d \times \beta\times(T-t)}
\big(
\sum_{i=1}^d\|\xi_{i}\|_{\infty}
+
d\times (\upsilon +
\varkappa(\hat{\Lambda}_{t})
\hat{\Lambda}_{t})(T-t)\big).

\dce
$$
\ethe

\brem
Notice that the continuity of a bounded process $(X_t)_{t\geq 0}$ does not imply automatically the continuity of $\|X_t\|_\infty, t\geq 0$.
\erem

\proof
For any $N\in\mathbb{N}^*$ we consider the BSDE$[T,\xi, \bar{f}_{N}]$(\ref{theBSDE}) and its solution $(Y_{N}, Z_{N})$ (which exists and unique because $\bar{f}_{N}$ is Lipschitzian in the two variables). Denote by $\Lambda_{N}$ and $\hat{\Lambda}_{N}$ the lambda functions defined with respect to $Y_{N}$ and its corresponding Malliavin derivatives (cf. \cite[Proposition 5.3]{EKPQ}). The continuity of $\Lambda_{N}$ and $\hat{\Lambda}_{N}$ is proved in Lemma \ref{lambda-hat-lambda}. By the same lemma, we have
$$
\dcb
\Lambda_{N,t}-\Lambda_{N,T}
\leq
\int_{t}^T d \times \beta\times \Lambda_{N,s} + d\times \upsilon+ d\times\varkappa(\hat{\Lambda}_{N,s})
\hat{\Lambda}_{N,s}\ ds
\\
\leq
\int_{t}^T d \times \beta\times \Lambda_{N,s}\ ds

+d\times (\upsilon+ \varkappa(\hat{\Lambda}_{N,t})
\hat{\Lambda}_{N,t})(T-t).
\dce
$$
By Gronwall's inequality (cf. Lemma \ref{gronwall-ineq}),$$
\dcb
\Lambda_{N,t}

\leq
e^{d \times \beta\times(T-t)}
\big(
\sum_{i=1}^d\|\xi_{i}\|_{\infty}
+
d\times (\upsilon+ \varkappa(\hat{\Lambda}_{N,t})
\hat{\Lambda}_{N,t})(T-t) \big)
\dce
$$
(cf. Lemma \ref{Y-phi}). By Lemma \ref{lambda-hat-lambda} again, for $0<t\leq T$,
$$
\dcb
&&\limsup_{t'\uparrow t}\frac{1}{t-t'}\big(\hat{\Lambda}_{N,t'}^2-\hat{\Lambda}_{N,t}^2\big)
\\

&\leq&
\hat{\Lambda}_{N,t}^2\times
\big(
d\times 2\beta 
+
d^2\times n\times \varkappa(2\hat{\Lambda}_{N,t})^2
\big)

+\hat{\Lambda}_{N,t}2\sqrt{d\times n}\big( 
\hat{\upsilon}+\hat{\beta}\Lambda_{N,t} 
+ \hat{\varkappa}(\hat{\Lambda}_{N,t})\hat{\Lambda}_{N,t}
\big)\\

&\leq&
\hat{\Lambda}_{N,t}^2\times
\big(
d\times 2\beta 
+
d^2\times n\times \varkappa(2\hat{\Lambda}_{N,t})^2
\big)\\
&&
+\hat{\Lambda}_{N,t}2\sqrt{d\times n}\big( 
\hat{\upsilon}+\hat{\beta}e^{d \times \beta\times(T-t)}
\big(
\sum_{i=1}^d\|\xi_{i}\|_{\infty}
+
d\times (\upsilon+ \varkappa(\hat{\Lambda}_{N,t})
\hat{\Lambda}_{N,t})(T-t) \big) 
+ \hat{\varkappa}(\hat{\Lambda}_{N,t})\hat{\Lambda}_{N,t}
\big)\\

&\leq&
\mathtt{g}(t, \hat{\Lambda}_{N,t}^2).
\dce
$$
Let $g(s)=\hat{\Lambda}_{N,T-s}^2$ for $0\leq s\leq T$. The above inequality becomes$$
\dcb
\limsup_{s'\downarrow s}\frac{1}{s'-s}(g(s')-g(s))=\limsup_{s'\downarrow s}\frac{1}{(T-s)-(T-s')}\big(\hat{\Lambda}_{N,T-s'}^2-\hat{\Lambda}_{N,T-s}^2\big)\\
\leq
\mathtt{g}(T-s, \hat{\Lambda}_{N,T-s}^2)=\mathtt{g}(T-s,g(s)).
\dce
$$
The theorem of differential inequality implies $\hat{\Lambda}_{N,t}^2=g(T-t)\leq u_{\mathtt{g}, T-t}$ for $T-\alpha_{\mathtt{g}}< t \leq T$.

For any $0\leq a< \alpha_{\mathtt{g}}$, the function $u_{\mathtt{g}}$ is uniformly bounded on $[0,a]$ by, say, a constant $C_{a}>0$. By the relation $D_\theta Y_{N,\theta} = Z_{N,\theta}$ of \cite[Proposition 5.3]{EKPQ}, $\|Z_{N}\|^2$ is bounded by this constant on $[T-a, T]$. By Remark \ref{strong-stationary}, for any $N>C_{a}$, $(Y_{N}, Z_{N})$ is actually a solution of the BSDE$[T,\xi, f]$(\ref{theBSDE}) on $[T-a,T]$. This proves the existence of the solution of the BSDE$[T,\xi, f]$(\ref{theBSDE}) on $[T-a,T]$ with $\mathfrak{M}d$-property.

On the other hand, if a solution $(Y',Z')$ of the BSDE$[T,\xi, f]$(\ref{theBSDE}) exists on the interval $[T-a, T]$ with $\mathfrak{M}d$-property, then, for any $N>C_{a}\vee \hat{\Lambda}'_{T-a}$, $(Y',Z')$ is a solution of the BSDE$[T,\xi, \bar{f}_{N}]$(\ref{theBSDE}) on $[T-a, T]$. The driver $\bar{f}_{N}$ being Lipschitzian, $(Y',Z')$ must coincide with $(Y_{N},Z_{N})=(Y,Z)$. This proves the uniqueness of the solution. \ok

\section{Applications: existence results for some classes of BSDEs delimited by Lyapunov functions}\label{definition-lyapunov}

In the last section we have established the two principles: the exponential integrability condition and the local solvability. As mentioned in Introduction, these two principles are expected to be helpful in the resolution of non Lipschitzian BSDEs. It is the subject of this section.

Precisely, we will study the following resolution scheme: For a given BSDE$[T,\xi, f]$(\ref{theBSDE}), let $\mathfrak{B}$ be the set of $0\leq b\leq T$ such that the BSDE$[T, \xi, f]$(\ref{theBSDE}) has a unique solution $(Y,Z)$ in the space $\mathcal{S}_{\infty}[b,T]\times \mathcal{Z}_{\mbox{\tiny BMO}}[b,T]$ on $[b,T]$ with $\mathfrak{M}d$-property. Let $b_{0}=\inf\mathfrak{B}$. In the case that the principle of local solvability holds (for example, under the conditions (\ref{package1}) and (\ref{explicit-Df})), to solve BSDE$[T,\xi, f]$(\ref{theBSDE}), it is enough to prove $b_{0}$ itself is a point in $\mathfrak{B}$ (so that necessarily $b_0=0$).

To establish the relation $b_0\in\mathfrak{B}$, we note that, in the light of Remark \ref{strong-stationary}, it is equivalent to prove $\hat{\Lambda}_{b_0+}=\lim_{b\downarrow b_0}\hat{\Lambda}_{b} <\infty$, which, according to Corollary \ref{varrho-estimation}, amounts to verify the principle of exponential integrability. John-Nirenberg inequality constitutes the standard for verifying exponential integrability. We are lead to estimate the BMO norm of the process $\partial_zf(t, Y_t, Z_t)$, or more precisely its sliceability (cf. the next section for definition). That is the key point of the above resolution scheme.

\subsection{The two notions}

We say that a predictable matrix valued process $X$ is uniformly sliceable on the time interval $[0,T]$, if, for any $c>0$, there exists $\delta >0$, for any $0\leq a < b \leq a+\delta$, the BMO norm
$$
\dcb
\|\ind_{(a,b)}X\|_\star^2=\sup_{0\leq t \leq T}\|\mathbb{E}[\int_{t}^{T} \ind_{\{a< s< b\}}\|X_{s}\|^2 ds\ |\mathcal{F}_{t}]\|_{\infty} \leq c.
\dce
$$
The uniform sliceability can be defined on any sub-interval of $[0,T]$, whatever closed or open, in a similar way.

The property of sliceability goes back to the study of \cite{schachermayer} on the BMO martingales. (The uniform sliceability in the present paper is stronger than the sliceability of \cite{schachermayer}.) But the present paper does not depends on the result of \cite{schachermayer}. 

As explained previously, the sliceability is useful in the resolution of non Lipschitzian BSDEs, because of its combination with John-Nirenberg inequality (cf. \cite[Chapitre VI, $n^\circ$105, Theorem]{DM}, \cite[Theorem 2.2]{K}) leading to exponential integrability conditions like (\ref{exp-boundedness}) in Theorem \ref{exp-integrability} or (\ref{exponential-functional}) in Corollary \ref{varrho-estimation}.
(See also \cite{HR} for other considerations about the sliceability.) That being said, it seems very challenging to get the sliceability from a general non Lipschitzian BSDE with simply Ito's calculus. The most significant development on that issue should be the work \cite{XZ} (cf. the explication in \cite{xing}), where, in the Markovian setting, it is shown that Lyapunov function can be a very efficient tool in establishing the sliceability. We decide therefore to test the idea of Lyapunov functions in the non Markovian context of our paper. 

\

We say that a function $\mathtt{h}$ is a (global) Lyapunov function for a driver $f$, if the function $\mathtt{h}$ is defined on $\mathbb{R}^d$ taking values in $\mathbb{R}_{+}$, two-times continuously differentiable, such that
\begin{equation}\label{liapunov-function}
\dcb
\frac{1}{2}\sum_{i=1}^d\sum_{i'=1}^d \partial_{y_{i}}\partial_{y_{i'}}\mathtt{h}(y)\sum_{j=1}^nz_{i,j}z_{i',j}
-
\sum_{i=1}^d\partial_{y_{i}}\mathtt{h}(y)f(i,t, y, z)
&\geq& 
\|z\|^2 - \mathtt{k}(t,y), 
\dce
\end{equation}
for some predictable function $\mathtt{k}(t,y)\geq 0$, for $\forall 0\leq t\leq T, \ \forall z\in\mathbb{R}^{d\times n}, \ \forall y \in\mathbb{R}^d$. The function $\mathtt{k}$ will be called a lower bound function of the Lyapunov condition (\ref{liapunov-function}). In this paper, we always suppose 
\begin{equation}\label{k-beta}
\dcb
|\mathtt{k}(t,y)-\mathtt{k}(t,0)|\leq \bar{\beta}\sum_{i=1}^d|y_{i}|\ \mbox{ for some $\bar{\beta}>0$ (the growth condition)},\\

\|\mathtt{k}\|_{\centerdot}=\sup_{1\leq i\leq d}\sup_{0\leq t\leq T}\|\mathbb{E}[\int_{t}^T  |\mathtt{k}(s,0)| ds\ |\mathcal{F}_{t}]\|_{\infty}<\infty \ (\mbox{the potential bound}).
\dce
\end{equation}

\

It is important to notice that the definition of Lyapunov function can be different from one paper to another. Actually, the Lyapunov functions in \cite{HR, XZ} have different lower bound functions from ours, and, in addition, the Lyapunov function in \cite{XZ} expresses a local condition.
That is why \cite{XZ} needs an extra \textit{\`a-priori} boundedness result \cite[Theorem 2.14]{XZ} to make Lyapunov functions work well in concrete situations. The notion of Lyapunov functions in the present paper is a global notion, and the boundedness results will be the direct consequence of the Lyapunov functions themselves 

The general properties of BSDEs with Lyapunov functions are exposed in Section \ref{about-LF}. The present section shows only how Lyapunov functions can be applied in concrete examples.

\

\subsection{Sub-quadratic BSDEs under Lyapunov condition}

Here is our first example. We say that the BSDE$[T,\xi, f]$(\ref{theBSDE}) is sub-quadratic (in $z$), if the driver $f$ has a $z$-increment rate function $\varkappa(r)$ bounded by constant times $1+r^\alpha$ for a $0< \alpha <1$ (the sub-quadratic index).

To illustrate our idea of the two fundamental principles, of the resolution scheme, of the sliceability and of the Lyapunov functions, we find perfect to consider the sub-quadratic BSDE. Actually, the sub-quadratic BSDE is a well-known situation and have been carefully studied in the literature. See \cite{bahlali, bahlali-al, CN}. With it, people can easily contrast our approach in the face of the literature.

\subsubsection{Immediate consequences of the Lyapunov condition}

The computation on sub-quadratic BSDEs is not as problematic as the computation on quadratic BSDEs. This is because of the Holder inequality$$
\dcb
\int_{b}^{b'}\|Z_{s}\|^{2\alpha} ds
\leq
\left(\int_{b}^{b'} ds\right)^{1-\alpha}\ \left(\int_{b}^{b'}\|Z_{s}\|^{2} ds\right)^{\alpha}
\dce
$$ 
(which is effective only when $\alpha<1$). We have the following lemma.

\bl\label{Z-sliceability}
Consider the BSDE$[T, \xi, f]$(\ref{theBSDE}) with a driver $f$ sub-quadratic in $z$ with sub-quadratic index $0<\alpha < 1$. Suppose the existence of a Lyapunov function for $f$. Suppose also the conditions on the driver $f$ in Corollary \ref{bound-Lyapunov}. Then, for any $0\leq a\leq T$, for any solution $(Y,Z)$ on $[a,T]$ of the BSDE$[T, \xi, f]$(\ref{theBSDE}) such that $Y\in\mathcal{S}_{\infty}[a,T]$, the process $\|Z\|^{\alpha}$ is uniformly sliceable on the time interval $[a,T]$. More precisely, for some constant $c_\alpha>0$ independent of $a$ and of $(Y,Z)$, for all $a\leq b<b'\leq T$,
$$
\dcb
\mathbb{E}[\int_{b}^{b'}\|Z_{s}\|^{2\alpha} ds\ |\mathcal{F}_{b}]
\leq
(b'-b)^{1-\alpha} c_{\alpha}^{2\alpha}.
\dce
$$
\el

\proof
It is to notice that the sub-quadratic condition in $z$ is compatible with the conditions in Corollary \ref{bound-Lyapunov}. As a consequence of Corollary \ref{bound-Lyapunov}, there exists a constant $c_\alpha>0$ such that
$$
\|\ind_{[a,T]}.Z\|_{\star}\leq c_{\alpha},
$$ 
for all $0\leq a\leq T$ and for any solution $(Y,Z)$ on $[a,T]$. Fix $0\leq a\leq T$ and let $a\leq b<b'\leq T$. We have  $$
\dcb
&&\mathbb{E}[\int_{b}^{b'}\|Z_{s}\|^{2\alpha} ds\ |\mathcal{F}_{b}]
\leq
\mathbb{E}[\left(\int_{b}^{b'} ds\right)^{1-\alpha}\ \left(\int_{b}^{b'}\|Z_{s}\|^{2} ds\right)^{\alpha}\ |\mathcal{F}_{b}]\\
&\leq&
(b'-b)^{1-\alpha}\left(\mathbb{E}[\int_{b}^{b'}\|Z_{s}\|^{2} ds\ |\mathcal{F}_{b}]\right)^{\alpha}
\leq
(b'-b)^{1-\alpha}\|\ind_{[a,T]}.Z\|_{\star}^{2\alpha}
\leq
(b'-b)^{1-\alpha} c_{\alpha}^{2\alpha}. \ \ok
\dce
$$

\

\bcor\label{Np-bound}
Suppose the same conditions in Lemma \ref{Z-sliceability}. Let $p>1$, $\delta = \big(\frac{1}{2p c_{\alpha}^{2\alpha}}\big)^{\frac{1}{1-\alpha}}$ and $N_{p}$ be the first integer number bigger than $\frac{T}{\delta}$. Then, for any $0\leq a\leq  T$, for any solution $(Y,Z)$ on $[a,T]$ of the BSDE$[T, \xi, f]$(\ref{theBSDE}) such that $Y\in\mathcal{S}_{\infty}[a,T]$, $$
\dcb
\mathbb{E}[e^{p\int_{a}^T\|Z_{s}\|^{2\alpha}ds}|\mathcal{F}_{a}] \leq 2^{N_{p}}.
\dce
$$
\ecor

\proof
It is the consequence of Lemma \ref{slice-bound} applied to $X=\sqrt{p}\|Z\|^\alpha$, because, by Lemma \ref{Z-sliceability}, for the $\delta>0$ of the corollary, $\|\ind_{(b,b')}.X\|_\star^2\leq \frac{1}{2}$ for any $a< b<b'\leq T, b'-b\leq \delta$.  \ok

\

\subsubsection{An existence result for sub-quadratic BSDE}

The above results imply a uniform estimate on the linear BSDE$[T, \xi, f, j, \theta]$(\ref{theBSDE-malliavin-*}).

\bl\label{finite-lambda}
Suppose the conditions (\ref{package1}) and (\ref{explicit-Df}) on BSDE$[T,\xi, f]$(\ref{theBSDE}). Suppose that the driver $f$ is sub-quadratic in $z$ and the Malliavin derivatives $D_{j,\theta}f$ of the driver are quadratic in $z$. Suppose that the driver $f$ has a global Lyapunov function. Then, for any $0\leq a\leq  T$, for any solution $(Y,Z)$ on $[a,T]$ of the BSDE$[T, \xi, f]$(\ref{theBSDE}) such that $Y\in\mathcal{S}_{\infty}[a,T]$, for any solution $(\hat{Y}, \hat{Z})$ on $[a,T]$ of the linear BSDE$[T, \xi, f, j, \theta]$(\ref{theBSDE-malliavin-*}), there exists a constant $c$, independent of $a$, such that $$
\sup_{1\leq j\leq n}\sup_{0\leq \theta<\infty}\sup_{1\leq i\leq d}\sup_{a\leq s\leq T} \|\hat{Y}_{i,s}\|_{\infty} \leq c.
$$
\el

Notice that, in the proof below, \cite{EKPQ} is not invoked because it is not applicable.

\

\proof
We apply Corollary \ref{varrho-estimation} (applied on $[a, T]$, cf. Remark \ref{valid-for-all-b}) to prove the lemma. It is enough to prove that the various quantities considered in Corollary \ref{varrho-estimation} relative to BSDE$[T, \xi, f, j, \theta]$(\ref{theBSDE-malliavin-*}) on $[a, T]$ are universally bounded with respect to $i, j, \theta, a$. 

It is obviously the case for $\|D_{j,\theta}\xi_{i}\|_{\infty}$ which are universally bounded according to (\ref{package1}). 

Consider next the driver root $\ind_{[a,T]}D_{j,\theta}f(\cdot, Y_{\cdot}, Z_{\cdot})$ of the BSDE$[T, \xi, f, j, \theta]$(\ref{theBSDE-malliavin-*}). Under the conditions of the lemma, the upper bounds on $(Y,Z)$ of Corollary \ref{bound-Lyapunov} on $[a,T]$ are valid. With the conditions (\ref{explicit-Df}) ($\hat{\varkappa}$ being affine), we see that the potential bounds $\|\ind_{[a,T]}D_{j,\theta}f(i, \cdot, Y_{\cdot}, Z_{\cdot})\|_{\centerdot}$ are bounded universally with respect to $i, j,\theta, a$. 

Consider thirdly the exponential functional $\varrho_{p}$ (cf. (\ref{exponential-functional}) for definition) relative to BSDE$[T, \xi, f, j, \theta]$(\ref{theBSDE-malliavin-*}). Recall
$$
|\partial_{z_{i',j'}}f(i,t, Y_{t}, Z_{t})| \leq \varkappa(2\|Z_{t}\|),\ 1\leq i\leq d, 1\leq i'\leq d,\ 1\leq j'\leq n.
$$
Hence, for any $p>1$ the exponential functional of the BSDE$[T, \xi, f, j, \theta]$(\ref{theBSDE-malliavin-*}) on the time interval $[a,T]$ is bounded by
$$
\varrho_{p}\leq\sup_{a\leq t\leq T}\|\ \mathbb{E}[e^{p\times d^2\times n \int_{t}^T \varkappa(2\|Z_{s}\|)^2ds}|\mathcal{F}_{t}]\ \|_{\infty}
\leq\sup_{a\leq t\leq T}\|\ \mathbb{E}[e^{p\times d^2\times n \int_{t}^T C(1+\|Z_{s}\|^{2\alpha})ds}|\mathcal{F}_{t}]\ \|_{\infty},
$$
for a constant $C>0$ (depending only on $\varkappa$). Under the conditions of the lemma, Corollary \ref{Np-bound} is valid, according to which, $$
\varrho_{p}\leq e^{p\times d^2\times n \times C\times T} 2^{N_{p\times d^2\times n \times C}},
$$
which is also independent of $i, j,\theta, a$. The lemma is proved. \ok

\bethe\label{sub-quadratic-driver}
Suppose 
\ebe
\item
the conditions (\ref{package1}) and (\ref{explicit-Df}) on BSDE$[T,\xi, f]$(\ref{theBSDE}),
\item
that the rate function $\hat{\varkappa}$ is an affine function,
\item
that the driver $f$ has a global Lyapunov function,
\item
that the driver $f$ is sub-quadratic in $z$. 
\dbe
Then, BSDE$[T,\xi, f]$(\ref{theBSDE}) has a unique solution on the whole time horizon $[0,T]$ with $\mathfrak{M}d$-property. 
\ethe

\proof
Recall $b_{0}=\inf\mathfrak{B}$. Notice that there exists a unique pair of processes $(Y,Z)$ on $(b_0,T]$ such that, for any $b_0<b\leq T$, restricted on $[b,T]$, $(Y,Z)$ is the solution of BSDE$[T,\xi, f]$(\ref{theBSDE}) in the space $\mathcal{S}_{\infty}[b,T]\times \mathcal{Z}_{\mbox{\tiny BMO}}[b,T]$ with $\mathfrak{M}d$-property. The quantity $\hat{\Lambda}_{b_0+}$ is well-defined.

Theorem \ref{solution-extend} shows $b_{0}<T$. Lemma \ref{finite-lambda} shows that $\hat{\Lambda}_{b_{0}+}<\infty$. This means, according to Remark \ref{strong-stationary}, that the solution of the BSDE$[T, \xi, f]$(\ref{theBSDE}) on $(b_{0},T]$ coincides with that of BSDE$[T,\xi, \bar{f}_{N}]$(\ref{theBSDE}) for some $N\in\mathbb{N}$, which in fact is a solution on $[b_0, T]$ with $\mathfrak{M}d$-property. We conclude $b_{0}\in\mathfrak{B}$. 

We say that $b_{0}=0$, which proves the theorem. If it was not the case, because the BSDE$[b_{0},Y_{b_{0}}, f]$(\ref{theBSDE}) satisfy the conditions (\ref{package1}) and (\ref{explicit-Df}), by Theorem \ref{solution-extend} (the existence of local solutions), the set $\mathfrak{B}$ would contain a $0\leq b<b_{0}$ which is in contradiction with the definition of $b_{0}$. \ok

\

Notice that Theorem \ref{sub-quadratic-driver} has been proved for a pedagogical reason. The reader are recommended to look at \cite{CN} for a study of sub-quadratic BSDEs without Malliavin derivatives. Note also that an extension of Theorem \ref{sub-quadratic-driver} will be presented in the following sections.

\

\subsubsection{A complementary remark with the sub-quadratic BSDEs}

Besides the two fundamental principles and the resolution scheme, sub-quadratic BSDEs serves also to explain another point. In fact, Theorem \ref{solution-extend} is proved for whole a class of BSDEs. It is therefore not optimal in particular situations. 
For example, in the case of a sub-quadratic BSDE, we can take the function $C(1+u)^{1+\alpha}$ ($C>0$ a constant) as the function $\mathtt{g}(t,u)$ of Theorem \ref{solution-extend}. Resolve the differential equation (\ref{u-equation}) associated with this function $\mathtt{g}$: $u'_{t} = \mathtt{g}(u_{t})=C(1+u)^{1+\alpha}$. We have $$
\dcb
-\frac{1}{\alpha}(1+u_{t})^{-\alpha} = \mbox{constant} + C\times t \ 
\mbox{ and } u_{0}= \sup_{0\leq \theta<\infty}\|D_{\theta}\xi\|_{\infty}^2,
\dce
$$ 
$$
\dcb
\left(\frac{1}{(1+\sup_{0\leq \theta<\infty}\|D_{\theta}\xi\|_{\infty}^2)^{-\alpha} - \alpha \times C\times t}\right)^{1/\alpha}
=
u_{t}.
\dce
$$
Therefore, the resolution horizon $\alpha_{\mathtt{g}}$, predicted by Theorem \ref{solution-extend} with the differential equation (\ref{u-equation}) is finite: $$
\alpha_{\mathtt{g}} = \frac{1}{(1+\sup_{0\leq \theta<\infty}\|D_{\theta}\xi\|_{\infty}^2)^{\alpha} \times \alpha \times C} < \infty.
$$
However, we have just proved that a sub-quadratic BSDE may very well have a solution on any time horizon $[0,T]$.

\

\subsection{Beyond the case of sub-quadratic BSDEs}\label{alpha-growth}

We have spoken about the example of sub-quadratic BSDEs, for pedagogical reason. In fact, our resolution scheme allows us to prove more general results than Theorem \ref{sub-quadratic-driver}.
Let us introduce the following growth condition:
\begin{equation}\label{sub-growth}
|f(t, y, z)| \leq \rho(1+|y|+ |z|^{1+\alpha})
\end{equation}
for some constant $\rho>0$ and $0\leq \alpha < 1$. 

Notice that, in the present paper, as explained in Section \ref{growth-increment-conditions}, we do not expressly distinguish the growth condition from the increment condition. However, such a distinction can be very useful. The discussion of this section serves as a good example.
We will prove that the growth condition (\ref{sub-growth}) combined with Lyapunov functions leads to the exponential integrability. Notice that a quadratic BSDE can satisfy the growth condition (\ref{sub-growth}). The following driver gives an example:$$
\dcb
f(i,t,y,z) = |\sum_{i'\neq i}y_{i'}|\cos(|z_i|^2)+\left|1+\sum_{i',j'}\sum_{i'',j''}z_{i',j'}z_{i'',j''}\right|^{\frac{1+\alpha}{2}}u(\left|1+\sum_{i',j'}\sum_{i'',j''}z_{i',j'}z_{i'',j''}\right|^{\frac{1-\alpha}{2}}),
\dce
$$
where
$$
u(r) = \int_0^r\sum_{m=0}^\infty(-1)^m\ind_{\{m<s\leq m+1\}} ds, \ r\geq 0.
$$

\bl\label{root-sub-quadratic} 
Consider the BSDE$[T,\xi, f]$(\ref{theBSDE}). Suppose the existence of a Lyapunov function and the conditions in Corollary \ref{bound-Lyapunov}. Suppose the growth condition (\ref{sub-growth}) on the driver $f$. Then, there exists a constant $C_1>0$ such that, for any $0\leq a\leq T$, for any $(Y,Z)$ solution of the BSDE$[T,\xi, f]$(\ref{theBSDE}) in $\mathcal{S}_{\infty}[a,T]\times \mathcal{Z}_{\mbox{\tiny BMO}}[a,T]$ on the time interval $[a, T]$, we have 
$$
\mathbb{E}[\int_{t}^{b}  |Z_{i,s} - \mu^b_{i, s}|^2 ds\ |\mathcal{F}_{t}]
\leq
4C_1 \rho ((1+C_1)(b-t)+ C_1^{1+\alpha}(b-t)^{\frac{1-\alpha}{2}}),\
a\leq t<b \leq T,
$$
where the process $\mu^b$ is defined by the relation
$$
\mathbb{E}[Y_b|\mathcal{F}_{t}]=\mathbb{E}[Y_b|\mathcal{F}_{0}]+\int_{0}^t \mu^b_{s}dB_{s}, 0\leq t\leq T.
$$
\el

\proof
By Corollary \ref{bound-Lyapunov}, there exists a constant $C_1>0$, independent of $a$ and of $(Y,Z)$, such that
\begin{equation}\label{YZC1}
\sup_{a\leq t\leq T}\|Y_t\|_\infty \leq C_1,\
\|\ind_{[a,T]}.Z\|_\star\leq C_1.
\end{equation}
Let $
M^b_{t}=\mathbb{E}[Y_b|\mathcal{F}_{t}]$. 
Clearly, $\|M^b_t\|_\infty\leq C_1$ uniformly for $a\leq t\leq b\leq T$.
We consider the process $Y_t - M^b_t$ on $[a,b]$:
\begin{equation}\label{Y-Mb}
\left\{
\dcb
d(Y-M^b)_{t} = -f(t, Y_{t}, Z_{t}) dt + (Z-\mu^b)_t dB_{t},\\
\\
(Y-M^b)_{b}=0.
\dce
\right.
\end{equation}
The conditions (\ref{sub-growth}) and (\ref{YZC1}) imply the validity of formula (\ref{ine-z}) for the above BSDE. Recall $0\leq \alpha < 1$. We get, for $a\leq t\leq b$, 
$$
\dcb
&&\mathbb{E}[\int_{t}^b  |Z_{i,s} - \mu^b_{i, s}|^2 ds\ |\mathcal{F}_{t}]
\leq
\mathbb{E}[\int_{t}^b 2|(Y_{i,s}-M^b_{i,s})f(i,s,Y_{s},Z_{s})|ds\ |\mathcal{F}_{t}]\\
&&
\leq
4C_1\mathbb{E}[\int_t^b\rho(1+|Y_s|+ |Z_s|^{1+\alpha}) ds |\mathcal{F}_t]\\
&&
\leq 
4C_1 (\rho(1+C_1)(b-t)+\rho C_1^{1+\alpha}(b-t)^{\frac{1-\alpha}{2}}).
\dce
$$
The lemma is proved. \ok

\

\bl\label{lambda-b}
Suppose 
\ebe
\item
the conditions (\ref{package1}) and (\ref{explicit-Df}) on BSDE$[T,\xi, f]$(\ref{theBSDE}),
\item
that the rate functions $\varkappa, \hat{\varkappa}$ are affine functions,
\item
the existence of a global Lyapunov function for the driver $f$,
\item
the growth condition (\ref{sub-growth}) on the driver $f$.
\dbe
Then, we have $\hat{\Lambda}_{b_0+}<\infty$.
\el

\proof
Recall $b_{0}=\inf\mathfrak{B}$. 
Notice that the assumptions of the lemma imply the conditions in Corollary \ref{bound-Lyapunov} so that Lemma \ref{root-sub-quadratic} is applicable.

Theorem \ref{solution-extend} shows $b_{0}<T$. Fix a $b_0<b\leq T$. By \cite[Proposition 1.2.8]{N}, by \cite[Lemma 5.1]{EKPQ}, the processes $M^b$ and $\mu^b$ are Malliavin differentiable. Recall that the $\mathfrak{M}d$-property (combined with Remark \ref{strong-stationary}) implies formula (\ref{theBSDE-malliavin}) for $D_\theta Y$ on $(b_0,T]\cap[\theta,T]$, from which we obtain,
for $1\leq j\leq n$, $0\leq \theta \leq b$, 
\begin{equation}\label{Y-Mb-equation}
\left\{
\dcb
dD_{j,\theta}(Y_{t}-M^b_t) 
&=& -(\partial_{y}f(t, Y_{t}, Z_{t})D_{j,\theta}Y_{t}+\partial_{z}f(t, Y_{t}, Z_{t})D_{j,\theta}Z_{t} + D_{j,\theta}f(t, Y_{t}, Z_{t}) )dt\\
&& + D_{j,\theta}(Z_{t}-\mu^b_t) dB_{t}\\
&=& -(\partial_{y}f(t, Y_{t}, Z_{t})D_{j,\theta}(Y_{t}-M^b_t)+\partial_{z}f(t, Y_{t}, Z_{t})D_{j,\theta}(Z_{t}-\mu^b_t))dt\\
&&-(\partial_{y}f(t, Y_{t}, Z_{t})D_{j,\theta}M^b_{t}+\partial_{z}f(t, Y_{t}, Z_{t})D_{j,\theta}\mu^b_{t} + D_{j,\theta}f(t, Y_{t}, Z_{t}) )dt\\
&& + D_{j,\theta}(Z_{t}-\mu^b_t) dB_{t},\ b_0< t\leq b, \theta\leq  t\leq b\\
\\
D_{j,\theta}(Y_{b}-M^b_b)&=&0.
\dce
\right.
\end{equation}
This is a linear BSDE in the sense of (\ref{linear-BSDE}) of the next section. Its $g$-coefficient is bounded and its $h$-coefficient is given by $
\partial_{z_{i',j'}}f(i,t, Y_{t}, Z_{t})
$
which is bounded by $\varkappa(2\|Z_t\|)$. We have$$
\sum_{j=1}^n\|h_{j,s}\|^2 \leq d^2  n \varkappa(2\|Z_t\|)^2
\leq
d^2 n \tilde{\gamma}(1+\|Z_t\|^2),
$$
for a constant $\tilde{\gamma}$. Consider the exponential integrability $\varrho_{28}(h)$ on the time interval $(b_0,b]\cap [\theta,b]$ associated with the linear BSDE(\ref{Y-Mb-equation}). For any $b_0<t\leq b$, we have $$
\mathbb{E}[e^{28\int_{t}^b \sum_{j=1}^n\|h_{j,s}\|^2ds}|\mathcal{F}_t]
\leq 
\mathbb{E}[e^{28\int_{t}^bd^2 n \tilde{\gamma}(1+\|Z_s\|^2)ds} |\mathcal{F}_t]
\leq
\mathbb{E}[e^{28\int_{t}^bd^2 n \tilde{\gamma}(1+2\|\mu^b_s\|^2+2\|Z_s-\mu^b_s\|^2)ds} |\mathcal{F}_t].
$$
As $b_0<b$, $\mathfrak{M}d$-property applies to $Y_b$ so that $DY_b$ is bounded by a constant $C_2>0$. By Clark-Ocone formula,$$
\|\mu^b_t\|_\infty
=\|\mathbb{E}[D_tY_b |\mathcal{F}_t]\|_\infty
\leq
C_2,
$$
so that
$$
\mathbb{E}[e^{28\int_{t}^b \sum_{j=1}^n\|h_{j,s}\|^2ds}|\mathcal{F}_t]
\leq 
e^{28 d^2 n \tilde{\gamma}(1+2C_2^2)(b-t)}
\mathbb{E}[e^{56d^2 n \tilde{\gamma}\int_{t}^b\|Z_s-\mu^b_s\|^2ds} |\mathcal{F}_t].
$$
Now, suppose that the number $b_0<b\leq T$ has been chosen in such a way that $$
56d^2 n \tilde{\gamma} d4C_1 \rho ((1+C_1)(b-b_0)+ C_1^{1+\alpha}(b-b_0)^{\frac{1-\alpha}{2}}) \leq \frac{1}{2}.
$$
By Lemma \ref{root-sub-quadratic} (the conditions in Corollary \ref{bound-Lyapunov} being satisfied), $$
\|\ind_{[t,b]}\sqrt{56d^2 n \tilde{\gamma}}(Z-\mu^b)\|_\star^2\leq \frac{1}{2}
$$
so that, for all $b_0<t\leq b$, by John-Nirenberg's inequality,
$$
\mathbb{E}[e^{28\int_{t}^b \sum_{j=1}^n\|h_{j,s}\|^2ds}|\mathcal{F}_t]
\leq 
e^{28 d^2 n \tilde{\gamma}(1+2C_2^2)(b-t)}
\mathbb{E}[e^{56d^2 n \tilde{\gamma}\int_{t}^b\|Z_s-\mu^b_s\|^2ds} |\mathcal{F}_t] \leq 2 e^{28 d^2 n \tilde{\gamma}(1+2C_2^2)(b-b_0)}<\infty.
$$
We conclude $\varrho_{28}(h)<\infty$ and its validity is independent of $j, \theta$.

Consider the root process of the linear BSDE(\ref{Y-Mb-equation}):$$
\partial_{y}f(t, Y_{t}, Z_{t})D_{j,\theta}M^b_{t}+\partial_{z}f(t, Y_{t}, Z_{t})D_{j,\theta}\mu^b_{t} + D_{j,\theta}f(t, Y_{t}, Z_{t}) 
$$
and its potential bound on $(b_0,b]\cap [\theta,b]$ (cf. Lemma \ref{linear-existence}). There are three terms. The first them is bounded uniformly with respect to $t\in (b_0,b]\cap [\theta,b]$ and to $j, \theta$ so as its potential bound. By \cite[Proposition 1.2.8]{N}, $D_{j,\theta}M^b$ is a martingale on $(b_0,b]\cap [\theta,b]$ bounded by $C_2$ (so that a BMO martingale). By \cite[Lemma 5.1]{EKPQ}, $D_{j,\theta}\mu^b$ is the coefficient process in the martingale predictable representation of $D_{j,\theta}M^b$ on $(b_0,b]\cap [\theta,b]$. By \cite[Theorem 10.11]{HWY}, the BMO norm $\|\ind_{(b_0,b]\cap [\theta,b]}.D_{j,\theta}\mu^b\|_\star$ is bounded by $2C_2$. By (\ref{YZC1}) and the affine $\varkappa$-rate condition, the BMO norm of $\partial_{z}f(t, Y_{t}, Z_{t})$ on $(b_0,b]\cap [\theta,b]$ is finite, independent of $j, \theta$. As a consequence, the potential bound of the second term is finite, independent of $j, \theta$. The conditions (\ref{explicit-Df}) and (\ref{YZC1}), together with the affine $\hat{\varkappa}$-rate condition, show that the last term has equally a finite potential bound on $(b_0,b]\cap [\theta,b]$, independent of $j, \theta$.

Lemma \ref{linear-existence} and Corollary \ref{varrho-estimation} are now applicable. By Lemma \ref{linear-existence}, $(D_{j,\theta}(Y-M^b), D_{j,\theta}(Z-\mu^b))$ is the unique solution of the linear BSDE(\ref{Y-Mb-equation}) on $(b_0,b]\cap [\theta,b]$. By Corollary \ref{varrho-estimation}, $D_{j,\theta}(Y-M^b)$ is bounded on $(b_0,b]\cap [\theta,b]$ uniformly with respect to $j, \theta$, which leads to the uniform boundedness of $D_{j,\theta}Y$ on $(b_0,b]$ (taking into account of $D_{j,\theta}Y_t=0$ if $\theta>t>b_0$). By $\mathfrak{M}d$-property, $D_{j,\theta}Y$ is uniformly bounded on $[b,T]$. Putting together the boundedness on the two intervals, we prove thus $\widehat{\Lambda}_{b_0+}<\infty$. \ok

\

We now conclude.

\bethe
Suppose the same assumptions as in Lemma \ref{lambda-b}. 
Then, BSDE$[T,\xi, f]$(\ref{theBSDE}) has a unique solution on the whole time horizon $[0,T]$ in the space $\mathcal{S}_{\infty}[0,T]\times \mathcal{Z}_{\mbox{\tiny BMO}}[0,T]$ with $\mathfrak{M}d$-property. 
\ethe

\proof
This theorem can be proved in the same way as Theorem \ref{sub-quadratic-driver} has been proved. \ok

\

\section{Linear BSDE with unbounded coefficients}\label{LBSDE}

Various results have been used in the proofs of Theorem \ref{exp-integrability}, Theorem \ref{solution-extend} and Theorem \ref{sub-quadratic-driver}. The most important ones are Corollary \ref{varrho-estimation}, Lemma \ref{lambda-hat-lambda}.
The remainder of the present paper will be devoted to the proofs of these results. In these proofs, we make use of many classical properties. Because of the pedagogic consideration of the present paper, we will, nevertheless, present all properties in great detail.

We begin with a discussion on the linear BSDEs, which will lead to, in particular, Corollary \ref{varrho-estimation}. 
\begin{equation}\label{linear-BSDE}
\left\{
\dcb
dY_{t} = -(f(t, 0, 0)+g_{t}Y_{t} + h_{t}Z_{t}) dt + Z_{t} dB_{t},\\
\\
Y_{T}=\xi,
\dce
\right.
\end{equation}
where $f(t, 0, 0)$ is a predictable process pathwisely Lebesgue integrable on any finite interval, and $g_{t}=(g_{i,i',t})_{1\leq i,i'\leq d}$ and $h_{t}=(h_{i,i',j,t})_{1\leq i,i'\leq d, 1\leq j\leq n}$ are matrix of locally bounded predictable processes and$$
g_{i,t}Y_{t} = \sum_{i'=1}^dg_{i,i',t}Y_{i',t},\
h_{i,t}Z_{t} = \sum_{j'=1}^n\sum_{i'=1}^dh_{i,i',j',t}Z_{i',j',t}.
$$
(Later in our computations, the coefficient $g$ will be linked with the Lipschitzian index $\beta$ while the coefficient $h$ will be linked with the $\varkappa$-rate condition.)

\subsection{Linear transformation $\phi$}

For any $0\leq a\leq T$ we introduce the matrix valued process $(\phi_{a,t})_{a\leq t\leq T}$ defined by the SDE:
\begin{equation}\label{phi-SDE}
\left\{
\dcb
d\phi_{i,i'',a,t} &=& \sum_{i'=1}^d\phi_{i,i',a,t}g_{i',i'',t}dt + \sum_{j'=1}^n\sum_{i'=1}^d\phi_{i,i',a,t}h_{i',i'',j',t}dB_{j',t},\\
&&\  \  \ \mbox{for } a\leq t\leq T,\\

\phi_{i,i'',a,a}&=&\ind_{\{i=i''\}}.
\dce
\right.
\end{equation}
In other words,
$
d\phi_{a,t} = \phi_{a,t}g_{t}dt + \sum_{j=1}^n\phi_{a,t}h_{j,t}dB_{j,t}.
$
Note that, as the coefficients $g$ and $h$ are locally bounded, the solution of the SDE(\ref{phi-SDE}) exists and is unique.

\subsubsection{Transformation to a local martingale}

\bl\label{phi-Y}
For any value process $Y$ of the linear BSDE(\ref{linear-BSDE}) (if it exists), the processes $$
\phi_{i,a,t}Y_{t}+\int_{a}^t\sum_{i''=1}^d\phi_{i,i'',a,s}f(i'',s,0,0)ds,\ 1\leq i\leq d,
$$ 
are local martingales on $[a,T]$. 
\el

\proof
By the integration by parts formula applied to the process $\phi_{a,t} Y_{t}$ on $[a,T]$, we can write 
$$
\dcb
&&
d(\phi_{i,a,t}Y_{t})=d(\sum_{i''=1}^d\phi_{i,i'',a,t}Y_{i'',t})
=\sum_{i'=1}^d(\phi_{i,i'',a,t}dY_{i'',t}+Y_{i'',t}d\phi_{i,i'',a,t} + d\cro{\phi_{i,i'',a}, Y_{i''}}_{t})\\

&=&
-\sum_{i''=1}^d\phi_{i,i'',a,t}f(i'',t,0,0)dt 
 
+
\sum_{i''=1}^d\phi_{i,i'',a,t}\sum_{j'=1}^nZ_{i'',j',t}dB_{j',t}
\\
&&
+ \sum_{i''=1}^dY_{i'',t}\sum_{j'=1}^n\sum_{i'=1}^d\phi_{i,i',a,t}h_{i',i'',j',t}dB_{j',t},\ \ 1\leq i\leq d. \ \ok
\dce
$$

\subsubsection{The inverse of $\phi$}

For $0\leq a\leq T$, we introduce the matrix valued process $(\psi_{a,t})_{a\leq t\leq T}$ defined by the SDE:
\begin{equation}\label{phi-SDE-inverse}
\left\{
\dcb
d\psi_{i'',i',a,t} &=& -\sum_{i'''=1}^dg_{i'',i''',t}\psi_{i''',i',a,t}dt 
+\sum_{j'=1}^n\sum_{i'''=1}^d\sum_{i^\circ=1}^d h_{i'',i''',j',t}h_{i''',i^\circ,j',t}\psi_{i^\circ,i',a,t}dt\\
&& - \sum_{j'=1}^n\sum_{i'''=1}^dh_{i'',i''',j',t}\psi_{i''',i',a,t}dB_{j',t},\
  \  \ \mbox{for } a\leq t\leq T,\\

\psi_{i'',i',a,a}&=&\ind_{\{i''=i'\}}.
\dce
\right.
\end{equation}
In other words,
$
d\psi_{a,t} = (-g_{t}\psi_{a,t}+\sum_{j=1}^nh_{j,t}h_{j,t}\psi_{a,t})dt + \sum_{j=1}^nh_{j,t}\psi_{a,t}dB_{j,t}.
$
By the integration by parts formula, we can prove that $d(\sum_{i''=1}^d\phi_{i,i'',a,t}\psi_{i'',i',a,t})=0$. We obtain

\bl
The matrix valued process $\psi_{a}$ is the inverse of $\phi_{a}$ on $[a,T]$. Consequently, for $0\leq a\leq b\leq T$, $\psi_{a,b}\phi_{a, t}=\phi_{b,t}$ on $t\in[b,T]$.
\el

\subsubsection{The norm $\|\phi_{a,t}\|$}

Consider now the matrix norm $
\|\phi_{a,t}\| = \sqrt{\sum_{i,i'=1}^d \phi_{i,i',a,t}^2}
$ on $[a,T]$.
We have
$$
\dcb
&&
d\|\phi_{a,t}\|^2 = d\sum_{i,i'=1}^d\phi_{i,i',a,t}^2\\ 
&=& \sum_{i,i'=1}^d\sum_{i''=1}^d 2\phi_{i,i',a,t}\phi_{i,i'',a,t}g_{i'',i',t}dt + \sum_{j=1}^n \sum_{i,i'=1}^d\sum_{i''=1}^d 2\phi_{i,i',a,t}\phi_{i,i'',a,t}h_{i'',i',j,t}dB_{j,t}\\
&&
+\sum_{j=1}^n\sum_{i,i'=1}^d(\sum_{i''=1}^d\phi_{i,i'',a,t}h_{i'',i',j,t})^2dt,
\dce
$$
which implies
$$
\dcb
&&
d\|\phi_{a,t}\| = d\sqrt{\|\phi_{a,t}\|^2}
=
\frac{1}{2}\frac{1}{\sqrt{\|\phi_{a,t}\|^2}}d\|\phi_{a,t}\|^2 + \frac{1}{2}\frac{1}{2}(-\frac{1}{2})\frac{1}{\sqrt{\|\phi_{a,t}\|^2}^3}d\cro{\|\phi_{a}\|^2,\|\phi_{a}\|^2}_{t} \\ 
&=& 
\frac{1}{2}\frac{1}{\|\phi_{a,t}\|}\sum_{i,i'=1}^d\sum_{i''=1}^d 2\phi_{i,i',a,t}\phi_{i,i'',a,t}g_{i'',i',t}dt 
+\frac{1}{2}\frac{1}{\|\phi_{a,t}\|} \sum_{j=1}^n \sum_{i,i'=1}^d\sum_{i''=1}^d 2\phi_{i,i',a,t}\phi_{i,i'',a,t}h_{i'',i',j,t}dB_{j,t}\\
&&
+\frac{1}{2}\frac{1}{\|\phi_{a,t}\|} \sum_{j=1}^n\sum_{i,i'=1}^d(\sum_{i''=1}^d\phi_{i,i'',a,t}h_{i'',i',j,t})^2dt

-\frac{1}{8}\frac{1}{\|\phi_{a,t}\|^3}\sum_{j=1}^n (\sum_{i,i'=1}^d\sum_{i''=1}^d 2\phi_{i,i',a,t}\phi_{i,i'',a,t}h_{i'',i',j,t})^2dt\\

&=& 
\|\phi_{t}\|
\big(
\cro{\bar{\phi}_{a,t} \boldsymbol{|} \bar{\phi}_{a,t}g_{t}}dt
+
\sum_{j=1}^n \cro{\bar{\phi}_{a,t} \boldsymbol{|} o_{j,a,t}} dB_{j,t}
+
\frac{1}{2}\sum_{j=1}^n(\|o_{j,a,s}\|^2
-\cro{\bar{\phi}_{a,s} \boldsymbol{|} o_{j,a,s}}^2)ds
\big),
\dce
$$
where the notation $\cro{\cdot \boldsymbol{|} \cdot}$ is defined in Section \ref{matrix} and $$
\dcb
\bar{\phi}_{i,i',a}=\frac{\phi_{i,i',a}}{\|\phi_{i,i',a}\|},\ \
o_{i,i',j,a,t}=\sum_{i''=1}^d\bar{\phi}_{i,i'',a,t}h_{i'',i',j,t}.
\dce
$$
The last equation is a linear SDE, which shows that $\|\phi_{a,t}\|$ is a Dolean-Dade exponential on $[a,T]$. Let $$
\dcb
\bar{o}_{i,i',j,a,t}=\frac{o_{i,i',j,a,t}}{\|o_{j,a,t}\|},\\
\mathfrak{n}_{a,s}
=
\cro{\bar{\phi}_{a,s} \boldsymbol{|} \bar{\phi}_{a,s}g_{s}}
+\frac{1}{2} \sum_{j=1}^n\|o_{j,a,s}\|^2
(1-\cro{\bar{\phi}_{a,s} \boldsymbol{|} \bar{o}_{j,a,s}}^2),\\
\mathcal{E}_{a,t}:=\mathcal{E}\big(
\int_{a}^t \sum_{j=1}^n \cro{\bar{\phi}_{a,s} \boldsymbol{|} o_{j,a,s}} dB_{j,s}\ \big).
\dce
$$
Notice that $
\cro{\bar{\phi}_{a,t} \boldsymbol{|} o_{j,a,t}}=
\|o_{j,a,t}\|\cro{\bar{\phi}_{a,t} \boldsymbol{|} \bar{o}_{j,a,t}}.
$

\

\bl\label{phi-norm}
We have$$
\|\phi_{a,t}\|
=
\sqrt{d}\exp\{
\int_{0}^t \cro{\bar{\phi}_{a,s} \boldsymbol{|} \bar{\phi}_{a,s}g_{s}} ds \
+\int_{0}^t \sum_{j=1}^n \cro{\bar{\phi}_{a,s} \boldsymbol{|} o_{j,a,s}} dB_{j,s}\
+\int_{0}^t\frac{1}{2} \sum_{j=1}^n\|o_{j,a,s}\|^2
(1-2\cro{\bar{\phi}_{a,s} \boldsymbol{|} \bar{o}_{j,a,s}}^2)ds
\}
$$
for $a\leq t\leq T$. In other words, 
$
\|\phi_{a,t}\|=\sqrt{d}\exp\{\int_{a}^t \mathfrak{n}_{a,s'}ds'\}\mathcal{E}_{a,t},\ a\leq t\leq T.
$
\el

\subsection{An existence result}

For constant $p>0$ let us define the exponential functional 
\begin{equation}\label{exponential-functional}
\varrho_{p}(h):=\sup_{0\leq t\leq T}\|\ \mathbb{E}[e^{p\int_{t}^T \sum_{j=1}^n \|h_{j,s}\|^2ds}|\mathcal{F}_{t}]\ \|_{\infty}.
\end{equation}

\bl\label{linear-existence}
Suppose that $\varrho_{28}(h)<\infty$. Suppose that the terminal value $\xi$ is bounded. Suppose that the components $g_{i,j}$ are bounded by $\beta>0$ and $$
\dcb
\|f\|_{\centerdot}=\sup_{1\leq i\leq d}\sup_{0\leq t\leq T}\|\ \mathbb{E}[\ \int_{t}^T  |f(i,s,0,0)|ds\ |\mathcal{F}_{t}]\ \|_{\infty}<\infty.
\dce
$$
Then, for $1\leq i\leq d$, the semimartingale
\begin{equation}\label{linear-solution}
\dcb
Y_{i,a}&=&\mathbb{E}[\phi_{i,a,T}\xi+\int_{a}^T\sum_{i''=1}^d\phi_{i,i'',a,s}f(i'',s,0,0)ds|\mathcal{F}_{a}], \ 0\leq a\leq T,
\dce
\end{equation}
is a well-defined process and the process $Y$, together with its martingale coefficient $Z$, forms a solution of the linear BSDE(\ref{linear-BSDE}). The process $(Y,Z)$ is the unique solution of the linear BSDE(\ref{linear-BSDE}) which satisfies $\sup_{\sigma}\|Y_{\sigma}\|_{\mathbf{L}^2}<\infty$, where $\sigma$ runs over all stopping times $\sigma\leq T$. 
\el

\brem
For an one dimensional Lipschitzian linear BSDE, the formula (\ref{linear-solution}) is known (cf. \cite[Proposition 2.2]{EKPQ}). In that case, the BSDE possesses a unique solution $(Y,Z)$, while the formula (\ref{linear-solution}) is deduced with a simple Ito's calculus. In Lemma \ref{linear-existence}, however, no \textit{\`a-priori} solution exists. We have to compute a $Z$ for the semimartingale $Y$ of formula (\ref{linear-solution}), and prove that $(Y,Z)$ is a solution of the BSDE.
\erem

\proof
We begin with the $\mathbf{L}^p$-estimation of $\|\phi_{a}\|$. Note that $\|o_{j,a,s}\|^2\leq \|h_{j,s}\|^2$ and $|\cro{\bar{\phi}_{a} \boldsymbol{|} \bar{\phi}_{a}g}|\leq \|g\|\leq d\times \beta$. For any $p>1$, for stopping time $\sigma\leq T$, consider
$$
\dcb
&&
\mathbb{E}[\mathcal{E}\big(
\sum_{j=1}^n \cro{\bar{\phi}_{a} \boldsymbol{|} o_{j,a}} \ind_{[a,T]}\centerdot B_{j}\ \big)_{\sigma}^p |\mathcal{F}_{a}]\\

&=&
\mathbb{E}[ 
\exp\{p\int_{a}^\sigma \sum_{j=1}^n \cro{\bar{\phi}_{a,s} \boldsymbol{|} o_{j,a,s}} dB_{j,s}\
-p\frac{1}{2}\int_{a}^\sigma \sum_{j=1}^n\|o_{j,a,s}\|^2
\cro{\bar{\phi}_{a,s} \boldsymbol{|} \bar{o}_{j,a,s}}^2ds
\}|\mathcal{F}_{a}]\\

&\leq&
\left(\mathbb{E}[ 
\exp\{
-p\int_{a}^\sigma \sum_{j=1}^n\|o_{j,a,s}\|^2
\cro{\bar{\phi}_{N,a,s} \boldsymbol{|} \bar{o}_{j,a,s}}^2ds
+
2p^2\int_{a}^\sigma \sum_{j=1}^n\|o_{j,a,s}\|^2
\cro{\bar{\phi}_{N,a,s} \boldsymbol{|} \bar{o}_{j,a,s}}^2ds\}|\mathcal{F}_{a}]\right)^{1/2} \\
&&
\left(\mathbb{E}[ 
\exp\{2p\int_{a}^\sigma \sum_{j=1}^n \cro{\bar{\phi}_{N,a,s} \boldsymbol{|} o_{j,a,s}} dB_{j,s}\
-\frac{4p^2}{2}\int_{a}^\sigma \sum_{j=1}^n\|o_{j,a,s}\|^2
\cro{\bar{\phi}_{N,a,s} \boldsymbol{|} \bar{o}_{j,a,s}}^2ds
\}|\mathcal{F}_{a}]\right)^{1/2}\\

&\leq&
\left(\mathbb{E}[ 
\exp\{
(2p^2-p)\int_{a}^T \sum_{j=1}^n \|h_{j,s}\|^2
 ds\}|\mathcal{F}_{a}]\right)^{1/2} = \sqrt{\varrho_{(2p^2-p)}(h)}.
\dce
$$
By Burkholder-Davis-Gundy's inequality (cf. \cite[Theorem 10.36]{HWY}), by Doob's inequality (cf. \cite[Theorem 2.15]{HWY}), by localization, we see that the martingale $\mathcal{E}\big(
\sum_{j=1}^n \cro{\bar{\phi}_{a} \boldsymbol{|} o_{j,a}} \ind_{[a,T]}\centerdot B_{j}\ \big)$ is in $\mathcal{H}^p$-martingale space, if $(\varrho_{(2p^2-p)}(h))^{\frac{1}{2p}}$ is finite. The relationship
$
\|\phi_{a,t}\|=\sqrt{d}\exp\{\int_{a}^t \mathfrak{n}_{a,s'}ds'\}\mathcal{E}_{a,t}
$
implies then  $$
\dcb
&&
\mathbb{E}[\sup_{a\leq t\leq T}\|\phi_{a,t}\|^p |\mathcal{F}_{a}]
=
\mathbb{E}[\sup_{a\leq t\leq T}(\sqrt{d}\exp\{\int_{a}^t \mathfrak{n}_{a,s'}ds'\}\mathcal{E}_{a,t})^p |\mathcal{F}_{a}]\\

&\leq&
\sqrt{d}^p e^{p\times d\times\beta(T-a)}\mathbb{E}[\exp\{\int_{a}^T \frac{p}{2} \sum_{j=1}^n\|o_{j,a,s}\|^2 ds'\}\sup_{a\leq t\leq T}\mathcal{E}_{a,t}^p |\mathcal{F}_{a}]\\

&\leq&
\sqrt{d}^p e^{p\times d\times\beta(T-a)}(\mathbb{E}[\exp\{p\int_{a}^T \sum_{j=1}^n\|o_{j,a,s}\|^2 ds'\} |\mathcal{F}_{a}])^{1/2}
(\mathbb{E}[\sup_{a\leq t\leq T}\mathcal{E}_{a,t}^{2p} |\mathcal{F}_{a}])^{1/2}\\

&\leq&
\sqrt{d}^p e^{p\times d\times\beta(T-a)}(\varrho_{p}(h))^{1/2}
(\frac{2p}{2p-1}(\mathbb{E}[\mathcal{E}_{a,T}^{2p} |\mathcal{F}_{a}])^{\frac{1}{2p}})^{p}\
\leq
\sqrt{d}^p e^{p\times d\times\beta(T-a)}(\varrho_{p}(h))^{1/2}
(\frac{2p}{2p-1})^{p}\varrho_{(8p^2-2p)}(h)^{\frac{1}{4}}.
\dce
$$
We next consider the definition of the process $Y$, especially the random variables $
\int_{0}^T\phi_{i,0,s}f(s,0,0)ds.
$
$$
\dcb
&&
\mathbb{E}[\int_{0}^T\sum_{i=1}^d\sum_{i''=1}^d\ |\phi_{i,i'',0,s}f(i'',s,0,0)|ds\ ]\\

&\leq&
\mathbb{E}[\int_{0}^T\sum_{i=1}^d|\phi_{i,0,s}|\ |f(s,0,0)|\ ds\ ]\\

&\leq&
\sqrt{d}\mathbb{E}[\int_{0}^T\|\phi_{i,0,s}\|\ |f(s,0,0)|\ ds\ ]
\leq
\sqrt{d}\mathbb{E}[\sup_{0\leq t\leq T}\|\phi_{i,0,t}\|\ \int_{0}^T|f(s,0,0)|\ ds|\mathcal{F}_{a}]\\

&\leq&
\sqrt{d}
\left(
\mathbb{E}[\sup_{0\leq t\leq T}\|\phi_{i,0,t}\|^2\ ]\right)^{1/2}
\left(
\mathbb{E}[\ (\int_{0}^T  |f(s,0,0)|ds)^2\ ]\right)^{1/2}.
\dce
$$
On the one hand, $$
\dcb
\left(
\mathbb{E}[\sup_{0\leq t\leq T}\|\phi_{i,0,t}\|^2\ ]\right)^{1/2}
\leq
\sqrt{d} e^{d\times\beta \times T}\varrho_{2}(h)^{1/4}
\frac{4}{3}\varrho_{28}(h)^{\frac{1}{8}}.
\dce
$$
On the other hand, the energy inequality implies $$
\dcb
\left(
\mathbb{E}[\ (\int_{0}^T |f(s,0,0)|ds)^2\ ]\right)^{1/2}
\leq \sqrt{2}d\|f\|_{\centerdot}.
\dce
$$
This means that the following martingale$$
\dcb
M_{i,a}:=\mathbb{E}[\phi_{i,0,T}\xi+\int_{0}^T\sum_{i''=1}^d\phi_{i,i'',0,s}f(i'',s,0,0)ds|\mathcal{F}_{a}], \ 0\leq a\leq T,
\dce
$$
is well defined so as its martingale coefficient process $J$:
$
M_{a}=M_{0}
+
\int_{0}^a \sum_{j'=1}^nJ_{j',s}dB_{j',s}, \ 0\leq a\leq T.
$
Consequently, the process $Y$ is well defined. Note that $Y_{i,t}$ coincides with
$$
\dcb
\sum_{i'=1}^d\psi_{i,i',0,t}(M_{i',t}-\int_{0}^t \sum_{i''=1}^d\phi_{i',i'',0,s}f(i'',s,0,0)ds).
\dce
$$
Hence,
$$
\dcb
dY_{i,t}
&=&
-\sum_{i'=1}^d\sum_{i'''=1}^dg_{i,i''',t}\psi_{i''',i',0,t} (M_{i',t}-\int_{0}^t \sum_{i''=1}^d\phi_{i',i'',0,s}f(i'',s,0,0)ds) dt \\
&&
+\sum_{i'=1}^d\sum_{j'=1}^n\sum_{i'''=1}^d\sum_{i^\circ=1}^d h_{i,i''',j',t}h_{i''',i^\circ,j',t}\psi_{i^\circ,i',0,t} (M_{i',t}-\int_{0}^t \sum_{i''=1}^d\phi_{i',i'',0,s}f(i'',s,0,0)ds) dt\\
&& 
- \sum_{i'=1}^d\sum_{j'=1}^n\sum_{i'''=1}^dh_{i,i''',j',t}\psi_{i''',i',0,t} (M_{i',t}-\int_{0}^t \sum_{i''=1}^d\phi_{i',i'',0,s}f(i'',s,0,0)ds) dB_{j',t}\\
&&
+
\sum_{i'=1}^d\sum_{j'=1}^n \psi_{i,i',0,t}J_{i',j',t}dB_{j',t}- \sum_{i'=1}^d\sum_{i''=1}^d\psi_{i,i',0,t}\phi_{i',i'',0,t}f(i'',t,0,0)dt\\
&&
-
\sum_{i'=1}^d \sum_{j'=1}^n \sum_{i'''=1}^dh_{i,i''',j',t}\psi_{i''',i',0,t}J_{i',j',t}dt\\

&=&
-\sum_{i'''=1}^dg_{i,i''',t}Y_{i''',t} dt

+\sum_{j'=1}^n\sum_{i'''=1}^d\sum_{i^\circ=1}^d h_{i,i''',j',t}h_{i''',i^\circ,j',t}Y_{i^\circ,t} dt\\
&&
- \sum_{j'=1}^n\sum_{i'''=1}^dh_{i,i''',j',t}Y_{i''',t} dB_{j',t}
+
\sum_{i'=1}^d\sum_{j'=1}^n \psi_{i,i',0,t}J_{i',j',t}dB_{j',t}\\
&&- f(i,t,0,0)dt
-
\sum_{i'=1}^d \sum_{j'=1}^n \sum_{i'''=1}^dh_{i,i''',j',t}\psi_{i''',i',0,t}J_{i',j',t}dt\\

&=&
- f(i,t,0,0)dt -\sum_{i'''=1}^dg_{i,i''',t}Y_{i''',t} dt 

-\sum_{j'=1}^n\sum_{i'''=1}^dh_{i,i''',j',t}Z_{i''',j'} dt

+ \sum_{j'=1}^n Z_{i,j'} dB_{j',t},
\dce
$$
where$$
\dcb
Z_{i,j'}=-\sum_{i^\circ=1}^dh_{i,i^\circ,j',t}Y_{i^\circ,t}+\sum_{i'=1}^d \psi_{i,i',0,t}J_{i',j',t}.
\dce
$$
The first part of the lemma is proved. 

Consider the second part of the lemma (the uniqueness). Actually, if $(Y,Z)$ is a solution of the linear BSDE(\ref{linear-BSDE}) such that $\sup_{\sigma}\|Y_{\sigma}\|_{\mathbf{L}^2}<\infty$, for any $0\leq a\leq T$, the processes $$
\dcb
\phi_{i,a,t}Y_{t}+\int_{a}^t\sum_{i''=1}^d\phi_{i,i'',a,s}f(i'',s,0,0)ds,\ 1\leq i\leq d,
\dce
$$
are uniformly integrable martingales on $[a,T]$, uniquely determined by their terminal values. The value process $Y$ is hence unique. Consequently, by the computation of Lemma \ref{phi-Y}, the processes $$
\sum_{i''=1}^d\phi_{i,i'',a,t}\sum_{j'=1}^nZ_{i'',j',t}dB_{j',t}, 1\leq i\leq d,
$$ 
are unique on $[a,T]$. This means that $\phi_a Z_{j'}, 1\leq j'\leq n$, are unique on $[a,T]$, and therefore, $Z$ is unique.
\ok

\subsection{Upper bounds on the solution of a linear BSDE}\label{upper-bound-by-phi}

\bcor\label{varrho-estimation}
Under the same conditions as in Lemma \ref{linear-existence}, the solution $(Y,Z)$ of the linear BSDE(\ref{linear-BSDE}) given in Lemma \ref{linear-existence} satisfies the estimations:
$$
\dcb
\|\sum_{i=1}^d|Y_{i,a}|\|_{\infty}
\leq
d e^{d\times\beta(T-a)}\varrho_{2}(h)^{1/4}
\frac{4}{3}\varrho_{28}(h)^{\frac{1}{8}}\ (\|\xi\|_{\infty}+\sqrt{2}d\|f\|_{\centerdot}), \ \forall 0\leq a\leq T,\\
\\
\|Z\|_{\star}^2
\leq
2\sum_{i=1}^d\|\xi_{i}\|_{\infty}^2
+2d^2\ \|f\|_{\centerdot}^2
+2\sup_{0\leq a\leq T}\|\sum_{i'=1}|Y_{i',a}|\|_{\infty}^2(1+2\beta\times d\times T
+
2d\sum_{i=1}^d \|h_{i}\|_{\star}^2).

\dce
$$
\ecor

\proof
We write$$
\dcb
&&
\sum_{i=1}^d|Y_{i,a}|= \sum_{i=1}^d|\ \mathbb{E}[\phi_{i,a,T}\xi+\int_{a}^T\sum_{i''=1}^d\phi_{i,i'',a,s}f(i'',s,0,0)ds|\mathcal{F}_{a}]\ |\\

&\leq&
\sum_{i=1}^d\mathbb{E}[|\phi_{i,a,T}\xi|\ |\mathcal{F}_{a}] +\sum_{i=1}^d\mathbb{E}[|\int_{a}^T\sum_{i''=1}^d|\sum_{i''=1}^d\phi_{i,i'',a,s}f(i'',s,0,0)|ds|\ |\mathcal{F}_{a}].
\dce
$$
Repeat then the computations in Lemma \ref{linear-existence}. We prove the estimation on $Y$.

Consider the estimation on $Z$. Let $0\leq a\leq b\leq T$. By the formula (\ref{ine-z}), we can write, for $a\leq t\leq b$,$$
\dcb
&&\mathbb{E}[\int_{t}^T  \|Z_{s}\|^2 ds\ |\mathcal{F}_{t}] - \sum_{i=1}^d\|\xi_{i}\|_{\infty}^2
=\mathbb{E}[\int_{t}^T \sum_{i=1}^d |Z_{i,s}|^2 ds\ |\mathcal{F}_{t}] - \sum_{i=1}^d\|\xi_{i}\|_{\infty}^2\\

&\leq&
2\sup_{0\leq a\leq T}\|\sum_{i'=1}|Y_{i',a}|\|_{\infty}\mathbb{E}[\int_{t}^T \sum_{i=1}^d |f(i,s,0,0)|+ \sum_{i=1}^d|g_{i,s}Y_{s}|+\sum_{i=1}^d \|h_{i,s}\| \|Z_{s}\| ds\ |\mathcal{F}_{t}]\\

&\leq&
2\sup_{0\leq a\leq T}\|\sum_{i'=1}|Y_{i',a}|\|_{\infty}\times d\times \|f\|_{\centerdot}\

+2\beta\times d\sup_{0\leq a\leq T}\|\sum_{i'=1}|Y_{i',a}|\|_{\infty}^2(T-t)\\
&&+
4\sup_{0\leq a\leq T}\|\sum_{i'=1}|Y_{i',a}|\|_{\infty}^2\mathbb{E}[\int_{t}^T \frac{1}{2}(\sum_{i=1}^d \|h_{i,s}\|)^2 ds\ |\mathcal{F}_{t}] + \mathbb{E}[\int_{t}^T \frac{1}{2} \|Z_{s}\|^2 ds\ |\mathcal{F}_{t}]\\

&\leq&
d^2\ \|f\|_{\centerdot}^2
+
\sup_{0\leq a\leq T}\|\sum_{i'=1}|Y_{i',a}|\|_{\infty}^2(1+2\beta\times d\times T)\\
&&+
2d\sup_{0\leq a\leq T}\|\sum_{i'=1}|Y_{i',a}|\|_{\infty}^2\sum_{i=1}^d \mathbb{E}[\int_{t}^T \|h_{i,s}\|^2 ds\ |\mathcal{F}_{t}] + \mathbb{E}[\int_{t}^T \frac{1}{2} \|Z_{s}\|^2 ds\ |\mathcal{F}_{t}].\ \ok

\dce
$$

\brem\label{valid-for-all-b}
For any $0\leq b\leq T$, the above boundedness can also be established for linear BSDEs defined the interval $[b,T]$. In fact, it is enough to extend the definition of the coefficients $g,h$ and the root process $f(t,0,0)$ in putting $g_t=0, h_t=0, f(t,0,0)=0$ for $0\leq t<b$ and apply the last corollary.
\erem

\section{BSDE(\ref{theBSDE}) under the classical Lipschitzian conditions}\label{epsilon-neighbor}

After Corollary \ref{varrho-estimation}, we consider now Lemma \ref{lambda-hat-lambda}. This lemma gives the infinitesimal variation of solutions of BSDEs. It is the consequence of a series of elementary computations on Lipschitzian BSDEs.
We introduce the following assumptions. 
\begin{equation}\label{bmo-boundedness-y-z}
\dcb
\mbox{\textbf{-}} & \mbox{$\xi$ is a bounded $\mathcal{F}_{T}$-measurable random variable,}\\

\mbox{\textbf{-}} & \mbox{the driver $f(s,y,z)$ is uniformly Lipschitzian in $y$ and in $z$,}\\
&\mbox{ with respectively indices $\beta>0$ and $\eta>0$},\\

\mbox{\textbf{-}} & \|f\|_{\centerdot}=\sup_{1\leq i\leq d}\sup_{0\leq t\leq T}\|\mathbb{E}[\int_{t}^T  |f(i,s,0,0)| ds\ |\mathcal{F}_{t}]\|_{\infty}<\infty \mbox{ (the potential bound)}. 
\dce
\end{equation}

\bpro\label{Lipschitzian-slice-collage}
Suppose the conditions (\ref{bmo-boundedness-y-z}) on the BSDE parameters $[T, \xi, f]$. Then, the BSDE$[T, \xi, f]$(\ref{theBSDE})
has a solution $(Y,Z)$ which is unique in the space $\mathcal{S}_{\infty}[0,T]\times \mathcal{Z}_{\mbox{\tiny BMO}}[0,T]$.
\epro

This proposition is classic. But, what we need here is of explicit estimations on the rate at which the process $(Y_{t}, Z_{t})$ varies in the space $\mathcal{S}_{\infty}[0,T]\times \mathcal{Z}_{\mbox{\tiny BMO}}[0,T]$.

\subsection{Gronwall's estimations}

We begin our computations with the simplest situation where the driver $f$ is free of the $Z$-component. To distinguish this particular situation from the general BSDE(\ref{theBSDE}), we rename it with a new number:
\begin{equation}\label{theBSDE-y-only}
\left\{
\dcb
dY_{t} = -f(t, Y_{t}) dt + Z_{t} dB_{t},\\
\\
Y_{T}=\xi.
\dce
\right.
\end{equation}
By Gronwall's inequality in Lemma \ref{gronwall-ineq}, we have

\bl\label{first-estimate}
Suppose that $(Y,Z)$ is a solution of the BSDE(\ref{theBSDE-y-only}). Suppose that
\ebe 
\item[\textbf{-}]
for every $s$, the driver $f(s,y)$ is Lipschitzian in $y$ with a deterministic Lipschitzian index $g(s)$,
\item[\textbf{-}]
$Z\!\centerdot\! B$ is a true martingale and $\mathbb{E}[\int_{t}^T g(s)  \sum_{i=1}^d|Y_{i,s}| ds]<\infty$,
\item[\textbf{-}]
the terminal value $\xi$ and the driver root $\int_{t}^{T} \sum_{i=1}^d |f(i, s, 0)|ds$ are integrable.
\dbe
Then, 
$
\sum_{i=1}^d|Y_{i,t}|
\leq
e^{d\int_{t}^Tg(s)ds}\ \mathbb{E}[\sum_{i=1}^d|\xi_{i}| + \int_{t}^{T} \sum_{i=1}^d |f(i, s, 0)|ds|\mathcal{F}_{t}]$,
for $
0\leq t\leq T.
$
If $\xi$ is a bounded random variable, if $$
\|f\|_{\centerdot}:=\sup_{1\leq i\leq d} \sup_{0\leq t\leq T}\|\mathbb{E}[\int_{t}^T |f(i,s,0)| ds\ |\mathcal{F}_{t}]\|_{\infty}<\infty,
$$
we have
$
\dcb
\sum_{i=1}^d\|Y_{i,t}\|_{\infty}
\leq
e^{d \int_{t}^{T} g(s)ds }(\sum_{i=1}^d\|\xi_{i}\|_{\infty}+d\times \|f\|_{\centerdot}),
\dce
$
for $0\leq t\leq T$.
\el

\subsection{Local estimations}

In this section we consider the infinitesimal variation of $(Y,Z)$. Let $0<\epsilon \leq T$, that will be called a neighborhood parameter below. We introduce the next BSDE with parameters $[T, \xi, f, \epsilon]$:
\begin{equation}\label{theBSDE-quadratic-epsilon}
\left\{
\dcb
dY_{t} = -\ind_{\{T-\epsilon\leq t\leq T\}}f(t, Y_{t}, Z_{t}) dt + Z_{t} dB_{t},\\
\\
Y_{T}=\xi.
\dce
\right.
\end{equation}
Clearly, the BSDE$[T, \xi, f]$(\ref{theBSDE}) and the BSDE$[T, \xi, f, \epsilon]$(\ref{theBSDE-quadratic-epsilon}) have same solutions on $[T-\epsilon, T]$.

\subsubsection{A transformation on $\mathcal{Z}_{\mbox{\tiny BMO}}[0,T]$}\label{BMO-transformation}

For any BMO martingale coefficient $Z^*$, let us consider the BSDE:
\begin{equation}\label{theBSDE-y-and-z-epsilon}
\left\{
\dcb
dY_{t} = -\ind_{\{T-\epsilon\leq t\leq T\}}f(t, Y_{t}, Z^*_{t}) dt + Z_{t} dB_{t},\\
\\
Y_{T}=\xi.
\dce
\right.
\end{equation}
Notice that the driver of the BSDE(\ref{theBSDE-y-and-z-epsilon}) is free of the martingale coefficient process $Z$. Suppose that the driver $f(t,y,z)$ satisfies the conditions in (\ref{bmo-boundedness-y-z}). Because $Z^*$ is a BMO martingale coefficient, the driver $f(t, y, Z^*_{t})$ also satisfies the conditions in (\ref{bmo-boundedness-y-z}) so that Proposition \ref{Lipschitzian-slice-collage} is applicable to the BSDE(\ref{theBSDE-y-and-z-epsilon}). Hence, the BSDE(\ref{theBSDE-y-and-z-epsilon}) has a unique solution in $(Y,Z)\in \mathcal{S}_{\infty}[0,T]\times \mathcal{Z}_{\mbox{\tiny BMO}}[0,T]$. The map $$
\Phi=\Phi_{\epsilon}:\ Z^* \ \longrightarrow \ Z
$$
is well-defined on the space $\mathcal{Z}_{\mbox{\tiny BMO}}[0,T]$.

\subsubsection{The contraction property of the map $\Phi$}\label{section-bounded-driver}

\bl\label{Phi-bounded-driver}
Suppose the conditions (\ref{bmo-boundedness-y-z}) on the BSDE parameters $[T, \xi,f]$. Suppose $0<\epsilon\leq T$ and  $$
c_{2}
=
c_{2}(d, \beta, \eta,\epsilon)
:=
\sqrt{2(d\times\beta \times \epsilon  +e^{-d\times \beta\times \epsilon})\times d}\ \times e^{d\times \beta\times \epsilon}\ \eta \ \sqrt{\epsilon} <1.
$$
Then, the map $\Phi$ is a contraction on the space $\mathcal{Z}_{\mbox{\tiny BMO}}[0,T]$. Precisely,$$
\|\Phi(Z'^*)-\Phi(Z^*)\|_{\star}
\leq c_{2}\  \|\ind_{[T-\epsilon,T]}(Z'^*-Z^*)\|_{\star}, \ \forall Z, Z^*\in \mathcal{Z}_{\mbox{\tiny BMO}}[0,T].
$$
\el

\proof
Consider two elements $Z^*, Z'^*\in\mathcal{Z}_{\mbox{\tiny BMO}}[0,T]$ 
and the corresponding solution $(Y,Z)$ and $(Y', Z')$ of the BSDE(\ref{theBSDE-y-and-z-epsilon}). 
We apply Lemma \ref{first-estimate} to the BSDE satisfied by $Y'-Y$ with null terminal value and driver $\ind_{\{T-\epsilon\leq t\leq T\}}(f(i,t,Y_{t}+y,Z'^*_{t})-f(i,t,Y_{t},Z^*_{t}))$. We notice that
$$
\dcb
&&\sup_{1\leq i\leq d}\sup_{0\leq t \leq T}\|\mathbb{E}[\int_{t}^T\ind_{\{T-\epsilon\leq s\leq T\}} |f(i,s,Y_{s},Z'^*_{s})-f(i,s,Y_{s},Z^*_{s})| ds\ |\mathcal{F}_{t}]\|_{\infty}\\
&&
\leq
\sup_{1\leq i\leq d}\sup_{0\leq t \leq T}\|\mathbb{E}[\int_{t}^T \ind_{\{T-\epsilon\leq s\leq T\}}\eta \|Z'^*_{s}-Z^*_{s}\| ds\ |\mathcal{F}_{t}]\|_{\infty}\

\leq
\eta\ \sqrt{\epsilon}\
\|\ind_{[T-\epsilon,T]}(Z'^*-Z^*)\|_{\star}.

\dce
$$
Lemma \ref{first-estimate} implies then, for $0\leq t\leq T$,
$$
\dcb
&&
\sum_{i=1}^d\|Y'_{i,t}-Y_{i,t}\|_{\infty}
\leq
e^{d\times \beta\times \epsilon}\ d\times \eta\ \sqrt{\epsilon}\
\|\ind_{[T-\epsilon,T]}(Z'^*-Z^*)\|_{\star}.
\dce
$$
Apply the identity (\ref{ine-z}) to estimate $Z'-Z$.
$$
\dcb
&&\mathbb{E}[\int_{t}^T \sum_{i=1}^d |Z'_{i,s}-Z_{i,s}|^2 ds\ |\mathcal{F}_{t}]\\

&\leq&
2\mathbb{E}[\int_{t}^T \ind_{\{T-\epsilon\leq s\leq T\}}  \sum_{i=1}^d|Y'_{i,s}-Y_{i,s}|  \beta \sum_{i'=1}^d|Y'_{i',s}-Y_{i',s}|\ ds\ |\mathcal{F}_{t}]\\
&&+
2 \mathbb{E}[\int_{t}^T \ind_{\{T-\epsilon\leq s\leq T\}}\sum_{i=1}^d|Y'_{i,s}-Y_{i,s}| \eta \|Z'^*_{s}-Z^*_{s}\|ds\ |\mathcal{F}_{t}]\\

&\leq&
2\beta(e^{d\times \beta\times \epsilon}\ d\times \eta\ \sqrt{\epsilon}\
\|\ind_{[T-\epsilon,T]}(Z'^*-Z^*)\|_{\star})^2\mathbb{E}[\int_{t}^T \ind_{\{T-\epsilon\leq s\leq T\}}  ds\ |\mathcal{F}_{t}]\\
&&+
2\eta\ e^{d\times \beta\times \epsilon}\ d\times \eta\ \sqrt{\epsilon}\
\|\ind_{[T-\epsilon,T]}(Z'^*-Z^*)\|_{\star} \ \mathbb{E}[\int_{t}^T \ind_{\{T-\epsilon\leq s\leq T\}} \|Z'^*_{s}-Z^*_{s}\|ds\ |\mathcal{F}_{t}]\\

&\leq&
2(d\times\beta \times \epsilon  +e^{-d\times \beta\times \epsilon})\times d\times e^{2d\times \beta\times \epsilon}\ \eta^2 \ \epsilon\ \|\ind_{[T-\epsilon,T]}(Z'^*-Z^*)\|_{\star}^2\\

\dce
$$
This proves the contraction property of the map $\Phi$.\ \ok

\subsubsection{The solution of the BSDE$[T, \xi,f, \epsilon]$(\ref{theBSDE-quadratic-epsilon})}

After Lemma \ref{Phi-bounded-driver}, the BSDE$[T, \xi,f, \epsilon]$(\ref{theBSDE-quadratic-epsilon}) can be solved by the fixed-point theorem. We now compute an upper bound for the solution.

\bpro\label{BSDE-epsilon-existence}
Suppose the conditions (\ref{bmo-boundedness-y-z}) on the BSDE parameters $[T, \xi, f]$. Suppose that neighborhood parameter $\epsilon$ satisfies the conditions in Lemma \ref{Phi-bounded-driver}. Then, the BSDE$[T, \xi, f, \epsilon]$(\ref{theBSDE-quadratic-epsilon})
has a unique solution $(Y,Z)$ in the space $\mathcal{S}_{\infty}[0,T]\times\mathcal{Z}_{\mbox{\tiny BMO}}[0,T]$. Moreover, if we define the processes $M=M_{\xi},\mu=\mu_{\xi}$ by the identity: $M_{t}=\mathbb{E}[\xi|\mathcal{F}_{t}]=\mathbb{E}[\xi|\mathcal{F}_{0}]+\int_{0}^t \mu_{s}dB_{s}, 0\leq t\leq T$, if we define
$$
F:=\sup_{1\leq i\leq d}\sup_{0\leq t\leq T}\|\mathbb{E}[\int_{t}^T \ind_{\{T-\epsilon\leq s\leq T\}}|f(i,s,M,0)| ds\ |\mathcal{F}_{t}]\|_{\infty}<\infty,
$$
the solution $(Y,Z)$ of the BSDE$[T, \xi, f, \epsilon]$(\ref{theBSDE-quadratic-epsilon}) satisfy:
$$
\dcb
\sum_{i=1}^d\|Y_{i,t}-M_{i,t}\|_{\infty}

&\leq&
e^{d\times \beta \times \epsilon}\ 
\big(1+d\times \eta \ \sqrt{\epsilon}\ \frac{1}{1-c_{2}}\
e^{d\times\beta\times \epsilon}\sqrt{2 \big(e^{-d\times\beta\times \epsilon}
+
\beta\times\epsilon\big)}\big) \times d\times F\\
&&
+e^{d\times \beta \times \epsilon}\ d\times \eta \ \sqrt{\epsilon}\ \frac{1}{1-c_{2}}\|\ind_{[T-\epsilon, T]}\mu\|_{\star},\
\forall 0\leq t\leq T,
\dce
$$
and
$$
\dcb
&&\|Z-\mu\|_{\star}
\leq
\frac{1}{1-c_{2}}\ e^{d\times\beta\times \epsilon}\sqrt{2 \big(e^{-d\times\beta\times \epsilon}
+
\beta\times\epsilon\big)}\times d\times F
+
\frac{c_{2}}{1-c_{2}}\|\ind_{[T-\epsilon, T]}\mu\|_{\star}.
\dce
$$ 
\epro

\proof
By the contraction property in Lemma \ref{Phi-bounded-driver}, 
the map $\Phi$ has a unique fixed point in the space $\mathcal{Z}_{\mbox{\tiny BMO}}[0,T]$. This fixed point $Z$ together with its associated value process $Y$ in the BSDE(\ref{theBSDE-y-and-z-epsilon}) forms a solution of the BSDE$[T, \xi, f, \epsilon]$(\ref{theBSDE-quadratic-epsilon}). 

To prove the inequalities on $(Y,Z)$, let $(Z^{(k)})_{k\geq 0}$ be the sequence defined by $Z^{(0)}\equiv 0$ and $Z^{(k)} = \Phi(Z^{(k-1)})$ for $k\geq 1$. We know that $Z^{(k)}$ converges to $Z$ in the space $\mathcal{Z}_{\mbox{\tiny BMO}}[0,T]$. According to Lemma \ref{Phi-bounded-driver} we have$$
\dcb
&&\|Z-\mu\|_{\star}\leq
\|Z^{(1)}-\mu\|_{\star}+ \|Z-Z^{(1)}\|_{\star}

\leq
\|Z^{(1)}-\mu\|_{\star}+
\sum_{k=1}^\infty\|Z^{(k+1)}-Z^{(k)}\|_{\star}\\
&&
\leq
\|Z^{(1)}-\mu\|_{\star}+
\sum_{k=1}^\infty c_{2}^{k}\|\ind_{[T-\epsilon, T]}(Z^{(1)}-Z^{(0)})\|_{\star}\\
&&
=
\|Z^{(1)}-\mu\|_{\star}+\frac{c_{2}}{1-c_{2}}\|\ind_{[T-\epsilon, T]}Z^{(1)}\|_{\star}

\leq
\frac{1}{1-c_{2}}\|Z^{(1)}-\mu\|_{\star}+\frac{c_{2}}{1-c_{2}}\|\ind_{[T-\epsilon, T]}\mu\|_{\star}.
\dce
$$
For every $k\geq 1$, let $Y^{(k)}$ be the value process of the BSDE(\ref{theBSDE-y-and-z-epsilon}) when  $Z^*=Z^{(k-1)}$. We apply Lemma \ref{first-estimate} to the BSDE of $Y^{(1)}-M$. With$$
F=\sup_{1\leq i\leq d}\sup_{0\leq t \leq T}\|\mathbb{E}[\int_{t}^T \ind_{\{T-\epsilon\leq s\leq T\}}|f(i,s,M,0)| ds\ |\mathcal{F}_{t}]\|_{\infty}<\infty,
$$
we can write
$
\sum_{i=1}^d\|Y^{(1)}_{i,t}-M_{i,t}\|_{\infty}
\leq
e^{d\times\beta\times \epsilon}\ d\times F,
$
for $0\leq t\leq T$.
By the formula (\ref{ine-z}),
$$
\dcb
&&\mathbb{E}[\int_{t}^T \sum_{i=1}^d |Z^{(1)}_{i,s}-\mu_{i,s}|^2 ds\ |\mathcal{F}_{t}]\\

&\leq&
2\sum_{i=1}^d\sup_{0\leq s\leq T}\|Y^{(1)}_{i,s}-M_{i,s}\|_{\infty}\mathbb{E}[\int_{t}^T  \ind_{\{T-\epsilon\leq s\leq T\}} \ |f(i,s,M,0)|ds\ |\mathcal{F}_{t}]\\
&&
+
2\mathbb{E}[\int_{t}^T \ind_{\{T-\epsilon\leq s\leq T\}} \sum_{i=1}^d|Y^{(1)}_{i,s}-M_{i,s}|\  \beta \sum_{i'=1}^d|Y^{(1)}_{i',s}-M_{i',s}|ds\ |\mathcal{F}_{t}]\\

&\leq&
2 e^{2\times d\times\beta\times \epsilon}\big(e^{-d\times\beta\times \epsilon}
+
\beta\times  \epsilon\big)\times d^2\times F^2.
\dce
$$
We conclude that
$
\|Z^{(1)}-\mu\|_{\star}
\leq
e^{d\times\beta\times \epsilon}\sqrt{2 \big(e^{-d\times\beta\times \epsilon}
+
\beta\times\epsilon\big)}\times d\times F,
$
and consequently
$$
\dcb
&&
\|Z-\mu\|_{\star}
\leq
\frac{1}{1-c_{2}}\|Z^{(1)}-\mu\|_{\star}+\frac{c_{2}}{1-c_{2}}\|\ind_{[T-\epsilon, T]}\mu\|_{\star}\\
&&
\leq
\frac{1}{1-c_{2}}\ e^{d\times\beta\times \epsilon}\sqrt{2 \big(e^{-d\times\beta\times \epsilon}
+
\beta\times\epsilon\big)}\times d\times F
+\frac{c_{2}}{1-c_{2}}\|\ind_{[T-\epsilon, T]}\mu\|_{\star}.
\dce
$$
As for $\sum_{i=1}^d\|Y_{i,t}-M_{i,t}\|_{\infty}$, recall the inequality in the proof of Lemma \ref{Phi-bounded-driver} applied to $Y^{(k)}, Y^{(k+1)}$:
$$
\dcb
&&
\sum_{i=1}^d\|Y^{(k+1)}_{i,t}-Y^{(k)}_{i,t}\|_{\infty}
\leq
e^{d\times \beta\times \epsilon}\ d\times \eta\ \sqrt{\epsilon}\
\|\ind_{[T-\epsilon, T]}(Z^{(k)}-Z^{(k-1)})\|_{\star}.
\dce
$$
Consequently,
$$
\dcb
&&\sum_{i=1}^d\|Y_{i,t}-M_{i,t}\|_{\infty}
\leq
\sum_{i=1}^d\|Y^{(1)}_{i,t}-M_{i,t}\|_{\infty}+\sum_{k=1}^\infty\sum_{i=1}^d\|Y^{(k+1)}_{i,t}-Y^{(k)}_{i,t}\|_{\infty}\\

&\leq&
\sum_{i=1}^d\|Y^{(1)}_{i,t}-M_{i,t}\|_{\infty}+\sum_{k=1}^\infty
e^{d\times \beta \times \epsilon}\ d\times \eta \ \sqrt{\epsilon}\ \|\ind_{[T-\epsilon, T]}(Z^{(k)}-Z^{(k-1)})\|_{\star}\\

&\leq&
\sum_{i=1}^d\|Y^{(1)}_{i,t}-M_{i,t}\|_{\infty}+
e^{d\times \beta \times \epsilon}\ d\times \eta \ \sqrt{\epsilon}\ \frac{1}{1-c_{2}}\|\ind_{[T-\epsilon, T]}Z^{(1)}\|_{\star}\\

&\leq&
\sum_{i=1}^d\|Y^{(1)}_{i,t}-M_{i,t}\|_{\infty}+
e^{d\times \beta \times \epsilon}\ d\times \eta \ \sqrt{\epsilon}\ \frac{1}{1-c_{2}}\|\ind_{[T-\epsilon, T]}(Z^{(1)}-\mu)\|_{\star}\\
&&
+e^{d\times \beta \times \epsilon}\ d\times \eta \ \sqrt{\epsilon}\ \frac{1}{1-c_{2}}\|\ind_{[T-\epsilon, T]}\mu\|_{\star}\\

&\leq&
e^{d\times \beta \times \epsilon}\ 
\big(1+d\times \eta \ \sqrt{\epsilon}\ \frac{1}{1-c_{2}}\
e^{d\times\beta\times \epsilon}\sqrt{2 \big(e^{-d\times\beta\times \epsilon}
+
\beta\times\epsilon\big)}\big) \times d\times F\\
&&
+e^{d\times \beta \times \epsilon}\ d\times \eta \ \sqrt{\epsilon}\ \frac{1}{1-c_{2}}\|\ind_{[T-\epsilon, T]}\mu\|_{\star}.
\dce
$$
The proposition is proved.
\ok

\subsection{Global estimation}

We now consider global estimations on the solution $(Y,Z)$. In fact, we only need to linearize the BSDE(\ref{theBSDE}) and to apply Corollary \ref{varrho-estimation}.

\bpro\label{value-bound}
The value process $Y$ of Proposition \ref{Lipschitzian-slice-collage} satisfies the inequality:$$
\dcb
\|\sum_{i=1}^d|Y_{i,a}|\|_{\infty}
\leq
d e^{d\times\beta(T-a)}e^{4\times d^2\times n\times \eta^2\times (T-a)}
\frac{4}{3}\ (\|\xi\|_{\infty}+\sqrt{2}d\|f\|_{\centerdot}), \ \forall 0\leq a\leq T,\\
\\
\|Z\|_{\star}^2
\leq
2\sum_{i=1}^d\|\xi_{i}\|_{\infty}^2
+2d^2\ \|f\|_{\centerdot}^2
+2\sup_{0\leq a\leq T}\|\sum_{i'=1}|Y_{i',a}|\|_{\infty}^2(1+2\beta\times d\times T
+
2d^3\times n\times \eta^2\times T).

\dce
$$
\epro

\proof
Under the condition of Proposition \ref{Lipschitzian-slice-collage} we rewrite the BSDE$[T,\xi, f]$(\ref{theBSDE}) in a linear form:$$
\dcb
&&dY_{t} = - f(t,Y_{t}, Z_{t}) dt + Z_{t}dB_{t}\\
&=&
- f(t,0,0) dt - (f(t,Y_{t}, 0)-f(t,0, 0)) dt - (f(t,Y_{t}, Z_{t})-f(t,Y_{t}, 0) dt + Z_{t}dB_{t}\\

&=&
- f(t,0,0) dt - \frac{f(t,Y_{t}, 0)-f(t,0, 0)}{\sum_{i'=1}^d|Y_{i',t}|}\sum_{i'=1}^d\mbox{sign}(Y_{i',t})Y_{i',t} dt\\
&& 
- \sum_{j'=1}^n\sum_{i'=1}^d\frac{f(t,Y_{t}, \mathbf{z}'_{i',j'})-f(t,Y_{t}, \mathbf{z}_{i',j'}))}{|Z_{i',j',t}|}\mbox{sign}(Z_{i',j',t})Z_{i',j',t} dt + Z_{t}dB_{t},
\dce
$$
where $\mathbf{z}_{i',j'}$ (resp. $\mathbf{z}'_{i',j'}$) is the matrix whose elements of alphabetical order inferior to (resp. inferior or equal to) $(i',j')$ are the elements of $Z_{t}$ in the same positions, and whose elements of alphabetical order superior or equal to (resp. superior to) $(i',j')$ are null. We see that $(Y,Z)$ is a solution of the linear BSDE(\ref{linear-BSDE}) with$$
\dcb
g_{i,i',t}=\frac{f(i,t,Y_{t}, 0)-f(i,t,0, 0)}{\sum_{i=1}^d|Y_{i,t}|}\mbox{sign}(Y_{i',t}),\\
\\
h_{i,i',j,t}
=
\frac{f(i,t,Y_{t}, \mathbf{z}'_{i',j'})-f(i,t,Y_{t}, \mathbf{z}_{i',j'})}{|Z_{i',j',t}|}\mbox{sign}(Z_{i',j',t}).
\dce
$$
Under the condition of Proposition \ref{Lipschitzian-slice-collage} it is verified that $
|g_{i,i',t}|\leq \beta,\ 
|h_{i,i',j,t}|\leq \eta.
$
Hence, for any constant $p>0$  $$
\varrho_{p}(h)=\sup_{0\leq t\leq T}\|\ \mathbb{E}[e^{p\int_{t}^T \sum_{j=1}^n \|h_{j,s}\|^2ds}|\mathcal{F}_{t}]\ \|_{\infty}
\leq
e^{p\times d^2\times n\times \eta^2\times (T-t)}.
$$
With this remark, the present lemma is a direct consequence of Corollary \ref{varrho-estimation}. \ok

\section{Malliavin derivatives of the BSDE$[T, \xi, f]$(\ref{theBSDE})}\label{malliavin-derivative}

In this section we present the results of \cite[Proposition 5.3]{EKPQ} on the Malliavin derivative of the BSDE$[T, \xi, f]$(\ref{theBSDE}). 
Suppose the conditions (\ref{bmo-boundedness-y-z}) on the parameters $[T, \xi, f]$ and the strong version of conditions (\ref{CR-C1}).
We also need the following ones.
\begin{equation}\label{maliavin-conditions}
\dcb
\mbox{\textbf{-}} & \mbox{$\xi\in\mathbb{D}^{1,2}$ and $\sup_{1\leq j\leq n}\sup_{1\leq i\leq d}\sup_{0\leq \theta\leq t}\|D_{j,\theta}\xi_{i}\|_{\infty}<\infty$.}
\dce
\end{equation}

\bl\label{elkaroui-al}(\cite[Proposition 5.3]{EKPQ})
Suppose the conditions (\ref{bmo-boundedness-y-z}), (\ref{CR-C1}) (strong version) and (\ref{maliavin-conditions}). Let $Y\in\mathcal{S}_{\infty}[0,T]$ and $Z\in \mathcal{Z}_{\mbox{\tiny BMO}}[0,T]$ be the solution in Proposition \ref{Lipschitzian-slice-collage} of the BSDE$[T, \xi, f]$(\ref{theBSDE}). Then, the random variables $Y_{t}, Z_{t}$ are Malliavin differentiable in $\mathbb{D}^{1,2}$ and the following BSDEs are satisfied: for $1\leq j\leq n$, $0\leq \theta \leq T$, 
\begin{equation}\label{theBSDE-malliavin}
\left\{
\dcb
dD_{j,\theta}Y_{t} 
&=& -(\partial_{y}f(t, Y_{t}, Z_{t})D_{j,\theta}Y_{t}+\partial_{z}f(t, Y_{t}, Z_{t})D_{j,\theta}Z_{t} + D_{j,\theta}f(t, Y_{t}, Z_{t}) )dt\\
&& + D_{j,\theta}Z_{t} dB_{t},\ \theta\leq t\leq T,\\
\\
D_{j,\theta}Y_{T}&=&D_{j,\theta}\xi,
\dce
\right.
\end{equation}
where the BSDE(\ref{theBSDE-malliavin}) has a unique solution in $\mathbf{L}^2$-sense. Moreover, the process $D_{\theta}Y_{\theta}, 0\leq \theta\leq T,$ is a version of $Z_{\theta}, 0\leq \theta\leq T$.
\el

In fact, we can have much better integrability bounds on $D_{j,\theta}Y$. We rewrite the BSDE(\ref{theBSDE-malliavin}) in its general form with parameters $[T, \xi, f, j,\theta]$:
\begin{equation}\label{theBSDE-malliavin-*}
\left\{
\dcb
d\hat{Y}_{t} 
&=& -(\partial_{y}f(t, Y_{t}, Z_{t})\hat{Y}_{t}+\partial_{z}f(t, Y_{t}, Z_{t})\hat{Z}_{t} + D_{j,\theta}f(t, Y_{t}, Z_{t}) )dt\\
&& + \hat{Z}_{t} dB_{t},\ 0\leq t\leq T,\\
\\
\hat{Y}_{T}&=&D_{j,\theta}\xi.
\dce
\right.
\end{equation}
(The process $(\hat{Y}_{i,t},\hat{Z}_{i,t})_{t\geq 0}$ depends on the parameter $j,\theta$. If necessary, we write them in the form $\hat{Y}_{j,\theta,i,t},\hat{Z}_{j,\theta,i,t}$.) By the Lipschitzian conditions of (\ref{bmo-boundedness-y-z}), for $1\leq i\leq d, 1\leq i'\leq d,\ 1\leq j'\leq n$, $$
\dcb
|\partial_{y_{i'}}f(i,t, Y_{t}, Z_{t})|\leq \beta,\ \
|\partial_{z_{i',j'}}f(i,t, Y_{t}, Z_{t})| \leq \eta.
\dce
$$
Suppose moreover 
\begin{equation}\label{Df-slice}
\dcb
\mbox{\textbf{-}} 
&\sup_{1\leq j\leq n}\sup_{0\leq \theta<\infty}\sup_{1\leq i\leq d}\sup_{0\leq t\leq T}\|\mathbb{E}[\int_{t}^T  |D_{j,\theta}f(i,s,Y_{s},Z_{s})| ds\ |\mathcal{F}_{t}]\|_{\infty}<\infty.

\dce
\end{equation}
Under the conditions (\ref{bmo-boundedness-y-z}), (\ref{CR-C1}) (strong version), (\ref{maliavin-conditions}) and (\ref{Df-slice}), by Proposition \ref{Lipschitzian-slice-collage}, the BSDE$[T, \xi, f, \theta, j]$(\ref{theBSDE-malliavin-*})
has a unique solution $(\hat{Y},\hat{Z})$ in the space $\mathcal{S}_{\infty}[0,T]\times \mathcal{Z}_{\mbox{\tiny BMO}}[0,T]$.
We have moreover a uniform control on the solution via Proposition \ref{value-bound}. Together with the identity $
Z_{i,j,\theta} = D_{j,\theta}Y_{i,\theta}=\hat{Y}_{\theta}
$
from Lemma \ref{elkaroui-al}, we obtain 

\bcor\label{MD-Z-bound}
Suppose the conditions (\ref{bmo-boundedness-y-z}), (\ref{CR-C1}) (strong version), (\ref{maliavin-conditions}) and (\ref{Df-slice}) on the BSDE parameters $[T, \xi, f]$. Let $Y\in\mathcal{S}_{\infty}[0,T]$ and $Z\in \mathcal{Z}_{\mbox{\tiny BMO}}[0,T]$ be the solution in Proposition \ref{Lipschitzian-slice-collage} of the BSDE$[T, \xi, f]$(\ref{theBSDE}). Then,  the martingale coefficient process satisfies $\sup_{0\leq s\leq T}\|Z_{s}\|_{\infty}<\infty$ and $Y_{t}\in\mathbb{D}^{1,2}$ for $0\leq t\leq T$ with $$
\sup_{1\leq j\leq n}\sup_{0\leq \theta<\infty}\sup_{1\leq i\leq d}\sup_{0\leq t\leq T}\|D_{j,\theta}Y_{i,t}\|_{\infty}<\infty.
$$
\ecor

\section{Infinitesimal estimations under $(\varkappa, \hat{\varkappa})$-rate conditions}

In the previous sections we have discussed the linear BSDEs, the Lipschitzian BSDEs, the Malliavin derivatives of the BSDEs. In this section we join the conditions (\ref{package1}) and (\ref{explicit-Df}) to our discussion. Precisely,
we consider the family of BSDEs which satisfy not only the conditions (\ref{bmo-boundedness-y-z}), (\ref{CR-C1}), (\ref{maliavin-conditions}), (\ref{Df-slice}), but also the conditions (\ref{package1}) and (\ref{explicit-Df}). (Clearly, many conditions are redundant. But we keep them entire to preserve the logical structure of our presentation, except (\ref{maliavin-conditions})). We suppose fixed all the parameters $T$ and $\xi$, $\upsilon$ and $\hat{\upsilon}$, $\beta$ and $\hat{\beta}$, $\varkappa$ and $\hat{\varkappa}$, except the parameter $\eta$ which is undetermined. We search for a uniform estimation on this family of BSDEs.

\subsection{Estimation of $Y$ on a neighborhood of $T$}

\bl\label{Y-phi}
Suppose the conditions (\ref{bmo-boundedness-y-z}). Let $Y\in\mathcal{S}_{\infty}$ and $Z\in \mathcal{Z}_{\mbox{\tiny BMO}}[0,T]$ be the solution in Proposition \ref{Lipschitzian-slice-collage} of the BSDE$[T, \xi, f]$(\ref{theBSDE}). Let $0<\epsilon \leq T$ be a neighborhood parameter. Suppose also that the driver $f$ satisfies the $\varkappa$-rate condition (\ref{kappa-increment}) in $z$. We have the inequality:$$
\dcb
\sum_{i=1}^d\|Y_{i,t}\|_{\infty}
\leq
e^{d \times \beta\times \epsilon}(\sum_{i=1}^d\|\xi_{i}\|_{\infty}+d\times (\upsilon\times\epsilon+ \sup_{T-\epsilon \leq t\leq T}\|\mathbb{E}[\int_{t}^T\varkappa(\|Z_{s}\|)\|Z_{s}\| ds\ |\mathcal{F}_{t}]\|_{\infty}),
\dce
$$
for $T-\epsilon \leq t\leq T$.
\el

\proof
It is Lemma \ref{first-estimate} applied to $Y$ considered as the value process of the BSDE$[T,\xi,f,\epsilon]$(\ref{theBSDE-quadratic-epsilon}) on $[T-\epsilon,T]$. We only need to explicit the estimation of the value $\|\ind_{[T-\epsilon,T]}f(\cdot, 0, Z_{\cdot})\|_{\centerdot}$:
$$
\dcb
\mathbb{E}[\int_{t}^T\ind_{\{T-\epsilon\leq s\leq T\}} |f(i,s,0,Z_{s})| ds\ |\mathcal{F}_{t}]

\leq
\upsilon\times\epsilon+ \mathbb{E}[\int_{t}^T\ind_{\{T-\epsilon\leq s\leq T\}} \varkappa(\|Z_{s}\|)\|Z_{s}\| ds\ |\mathcal{F}_{t}],\
0\leq t\leq T.\  \ok
\dce
$$

\subsection{Estimation of $DY$  on a neighborhood of $T$}\label{DY-at-neighborhood}

We consider now the Malliavin derivatives $D_{j,\theta}Y$, $1\leq j\leq n, 0\leq \theta\leq T$.

\subsubsection{Estimation of the difference $D_{j,\theta}Y-\hat{M}_{j,\theta}$}\label{DY-lipschitzian}

We suppose the conditions (\ref{bmo-boundedness-y-z}), (\ref{CR-C1}), (\ref{Df-slice}) and (\ref{package1}), (\ref{explicit-Df}). We recall that, for given indices $j,\theta$, $(D_{j,\theta}Y_{t}, D_{j,\theta}Z_{t})$ coincides with the solution $(\hat{Y}_{t}, \hat{Z}_{t})$ of the BSDE$[T, \xi, f, j,\theta]$(\ref{theBSDE-malliavin-*}) on the time interval $[\theta, T]$. 
We will now apply Proposition \ref{BSDE-epsilon-existence} to the BSDE$[T, \xi, f, j,\theta]$(\ref{theBSDE-malliavin-*}). To do so, we introduce the processes $(\hat{M},\hat{\mu})=(\hat{M}_{j,\theta},\hat{\mu}_{j,\theta})$ defined by the identity: $$
\hat{M}_{t}=\mathbb{E}[D_{j,\theta}\xi|\mathcal{F}_{t}]=\mathbb{E}[D_{j,\theta}\xi|\mathcal{F}_{0}]+\int_{0}^t \hat{\mu}_{s}dB_{s}, 0\leq t\leq T,
$$
and the constant $\hat{F}=
\hat{F}_{j,\theta,\epsilon}$ defined by
$$
\dcb
&&\hat{F}
=
\sup_{1\leq i\leq d}\sup_{T-\epsilon\leq t\leq T}\|\mathbb{E}[\int_{t}^T \ind_{\{T-\epsilon\leq s\leq T\}} |\partial_{y}f(i,s, Y_{s}, Z_{s})\hat{M}_{j,\theta,s} + D_{j,\theta}f(i,s, Y_{s}, Z_{s})|\ ds\ |\mathcal{F}_{t}]\|_{\infty}\\

&\leq&
\sup_{1\leq i\leq d}\sup_{T-\epsilon\leq t\leq T}\|\mathbb{E}[\int_{t}^T \ind_{\{T-\epsilon\leq s\leq T\}} \beta\sum_{i'=1}^d|\hat{M}_{j,\theta,i,s}| \ ds\ |\mathcal{F}_{t}]\|_{\infty}\\
&&
+\sup_{1\leq i\leq d}\sup_{T-\epsilon\leq t\leq T}\|\mathbb{E}[\int_{t}^T \ind_{\{T-\epsilon\leq s\leq T\}} 
(\hat{\upsilon} + \hat{\beta}\sum_{i'=1}^d|Y_{i,s}| + \hat{\varkappa}(\|Z_{s}\|)\|Z_{s}\|)\ ds\ |\mathcal{F}_{t}]\|_{\infty}.
\dce
$$
We also compute (where $\sigma\leq T$ denotes stopping times),  $$
\dcb
&&
\|\hat{\mu}\|_{\star}=\|\hat{\mu}_{j,\theta}\|_{\star}
=
\sqrt{\sup_{\sigma\leq T}\|\mathbb{E}[\int_{\sigma}^T \sum_{j'=1}^n\sum_{i'=1}^d\mu_{j,\theta,i',j's}^2 ds\ |\mathcal{F}_{\sigma}]\|_{\infty}}
\\
&=&
\sqrt{\sup_{\sigma\leq T}\|\sum_{i'=1}^d \mathbb{E}[ (D_{j,\theta}\xi_{i'})^2\ |\mathcal{F}_{\sigma}] -\sum_{i'=1}^d \hat{M}_{j,\theta,i',\sigma}^2\|_{\infty}}.
\dce
$$
As$$
0\leq \sum_{i'=1}^d \mathbb{E}[ \hat{M}_{j,\theta,i',T}^2\ |\mathcal{F}_{\sigma}] -\sum_{i'=1}^d \hat{M}_{j,\theta,i',\sigma}^2
\leq
\sum_{i'=1}^d \mathbb{E}[ \hat{M}_{j,\theta,i',T}^2\ |\mathcal{F}_{\sigma}],
$$
we obtain
$$
\dcb
\|\hat{\mu}\|_{\star}
\leq
\sqrt{\sup_{\sigma\leq T}\| \mathbb{E}[ \sum_{i'=1}^d(D_{j,\theta}\xi_{i'})^2\ |\mathcal{F}_{\sigma}]\|_{\infty}}

\leq
\|D_{j,\theta}\xi\|_{\infty}
<\infty.
\dce
$$
At last, we recall that the driver of the BSDE$[T, \xi, f, j,\theta]$(\ref{theBSDE-malliavin-*}) is Lipschitzian in $z$ with Lipschitzian index $\sqrt{d\times n}\times \eta$. 

\bl\label{DY-L-K}
Under the conditions (\ref{bmo-boundedness-y-z}), (\ref{CR-C1}), (\ref{Df-slice}) and (\ref{package1}), (\ref{explicit-Df}), for any neighborhood parameter $\epsilon>0$ which satisfies the condition in Lemma \ref{Phi-bounded-driver} with respect to the BSDE$[T, \xi, f, j,\theta]$(\ref{theBSDE-malliavin-*}) (i.e., $c_{2}(d,\beta,\sqrt{d\times n}\times \eta,\epsilon)<1$), for $T-\epsilon\leq t\leq T$,
$$
\dcb
&&
\sup_{1\leq j\leq n}\sup_{0\leq \theta\leq t}\sum_{i=1}^d\|D_{j,\theta}Y_{i,t}-\hat{M}_{j,\theta,i,t}\|_{\infty}\\
&\leq&
e^{d\times \beta \times \epsilon}\ 
\big(1+d\times \sqrt{d\times n}\times \eta \times \sqrt{\epsilon}\ \frac{1}{1-c_{2}(d,\beta,\sqrt{d\times n}\times \eta,\epsilon)}\
e^{d\times\beta\times \epsilon}\sqrt{2 \big(e^{-d\times\beta\times \epsilon}
+
\beta\times\epsilon\big)}\big) \times d\times\\
&& 
\big(
\beta\times \sup_{1\leq j\leq n}\sup_{0\leq \theta\leq T}\sum_{i=1}^d\|D_{j,\theta}\xi_{i}\|_{\infty} \\
&&
+ \hat{\upsilon} +\hat{\beta}\times\sup_{0\leq s\leq T}\sum_{i=1}^d\|Y_{i,s}\|_{\infty} 
+\hat{\varkappa}(\sup_{0\leq s\leq T}\|Z_{s}\|)\sup_{0\leq s\leq T}\|Z_{s}\|
\big)\times \epsilon\\
&&
+e^{d\times \beta \times \epsilon}\ d\times\sqrt{d\times n}\times  \eta \times \sqrt{\epsilon}\ \frac{1}{1-c_{2}(d,\beta,\sqrt{d\times n}\times \eta,\epsilon)}
\sup_{1\leq j\leq n}\sup_{0\leq \theta\leq T}\|D_{j,\theta}\xi\|_{\infty}.
\dce
$$
\el

\subsubsection{Estimation of $(D_{j,\theta}Y_{i})^2$}

Notice that, under the conditions (\ref{package1})  
$$
|\partial_{z_{i',j'}}f(i,t, Y_{t}, Z_{t})| \leq \varkappa(2\|Z_{t}\|),\ 1\leq i\leq d, 1\leq i'\leq d,\ 1\leq j'\leq n.
$$

\bl\label{sumY2}
Suppose the conditions (\ref{bmo-boundedness-y-z}), (\ref{CR-C1}) and (\ref{Df-slice}). Suppose the conditions (\ref{package1}) and (\ref{explicit-Df}). 
We have, for $1\leq j\leq n, 0\leq \theta<\infty, 0\leq t\leq T$,
$$
\dcb

\sum_{i=1}^d(D_{j,\theta}Y_{i,t})^2
&\leq&
\mathbb{E}[\sum_{i=1}^d(D_{j,\theta}\xi_{i})^2|\mathcal{F}_{t}]\\
&&
+\mathbb{E}[\int_{t}^T  \sum_{i=1}^d(D_{j,\theta}Y_{i,s})^2\ \times 
\big(
d\times 2\beta 
+
d^2\times n\times \varkappa(2\|Z_{s}\|)^2\big)
ds\ |\mathcal{F}_{t}]\\

&&
+\mathbb{E}[\int_{t}^T (\sum_{i=1}^d(D_{j,\theta}Y_{i,s})^2)^{1/2} 2\sqrt{d}\big(\hat{\upsilon}+\hat{\beta}\sum_{i'=1}^n|Y_{i',s}| 
+ \hat{\varkappa}(\|Z_{s}\|)\|Z_{s}\|\big)ds\ |\mathcal{F}_{t}].
\dce
$$
\el

\proof
If $t<\theta$, the statement is trivially true, because $D_{\theta}Y_{t}\equiv 0$. Suppose $\theta\leq t$ so that $D_{j,\theta}Y_{t}=\hat{Y}_{t}$. With the formula (\ref{ine-z}) applied to $(\hat{Y}, \hat{Z})$ for $1\leq j\leq n$, $0\leq \theta\leq T$, $1\leq i\leq d$, $0\leq t\leq T$, we write$$
\dcb
&&
\mathbb{E}[\hat{Y}_{i,t}^2 + \int_{t}^T  |\hat{Z}_{i,s}|^2 ds\ |\mathcal{F}_{t}]\\
&=&
\mathbb{E}[(D_{j,\theta}\xi_{i})^2+\int_{t}^T 2\hat{Y}_{i,s}\big(\partial_{y}f(i,s,Y_{s}, Z_{s})\hat{Y}_{s}
+\partial_{z}f(i,s,Y_{s}, Z_{s})\hat{Z}_{s}

+D_{j,\theta}f(i,s,Y_{s}, Z_{s})\big)\ ds\ |\mathcal{F}_{t}]\\

&\leq&
\mathbb{E}[(D_{j,\theta}\xi_{i})^2|\mathcal{F}_{t}]\\
&&
+\mathbb{E}[\int_{t}^T 2 |\hat{Y}_{i,s}| \big(\beta\sum_{i'=1}^d|\hat{Y}_{i',s}|
+\sqrt{d\times n}\varkappa(2\|Z_{s}\|)\|\hat{Z}_{s}\|\big)ds\ |\mathcal{F}_{t}]\\

&&
+\mathbb{E}[\int_{t}^T 2 |\hat{Y}_{i,s}| \big(|D_{j,\theta}f(i,s,0,0)|+\hat{\beta}\sum_{i'=1}^n|Y_{i',s}| 
+ \hat{\varkappa}(\|Z_{s}\|)\|Z_{s}\|\big)ds\ |\mathcal{F}_{t}]\\

&\leq&
\mathbb{E}[(D_{j,\theta}\xi_{i})^2|\mathcal{F}_{t}]\\
&&
+\mathbb{E}[\int_{t}^T |\hat{Y}_{i,s}|\sum_{i'=1}^d|\hat{Y}_{i',s}|\times 2\beta ds\ |\mathcal{F}_{t}]

+\mathbb{E}[\int_{t}^T  |\hat{Y}_{i,s}|^2\times d^2\times n\times \varkappa(2\|Z_{s}\|)^2 + \frac{1}{d}\|\hat{Z}_{s}\|^2ds\ |\mathcal{F}_{t}]\\

&&
+\mathbb{E}[\int_{t}^T 2 |\hat{Y}_{i,s}| \big(\hat{\upsilon}+\hat{\beta}\sum_{i'=1}^n|Y_{i',s}| 
+ \hat{\varkappa}(\|Z_{s}\|)\|Z_{s}\|\big)ds\ |\mathcal{F}_{t}].
\dce
$$
From that inequality, we deduce
$$
\dcb
\sum_{i=1}^d\hat{Y}_{i,t}^2
&\leq&
\mathbb{E}[\sum_{i=1}^d(D_{j,\theta}\xi_{i})^2|\mathcal{F}_{t}]\\
&&
+\mathbb{E}[\int_{t}^T (\sum_{i=1}^d|\hat{Y}_{i,s}|)^2\times 2\beta ds\ |\mathcal{F}_{t}]

+\mathbb{E}[\int_{t}^T  \sum_{i=1}^d|\hat{Y}_{i,s}|^2\times d^2\times n\times \varkappa(2\|Z_{s}\|)^2 ds\ |\mathcal{F}_{t}]\\

&&
+\mathbb{E}[\int_{t}^T (\sum_{i=1}^d |\hat{Y}_{i,s}|) 2\big(\hat{\upsilon}+\hat{\beta}\sum_{i'=1}^n|Y_{i',s}| 
+ \hat{\varkappa}(\|Z_{s}\|)\|Z_{s}\|\big)ds\ |\mathcal{F}_{t}]\\

&\leq&
\mathbb{E}[\sum_{i=1}^d(D_{j,\theta}\xi_{i})^2|\mathcal{F}_{t}]\\
&&
+\mathbb{E}[\int_{t}^T \sum_{i=1}^d|\hat{Y}_{i,s}|^2\times d\times 2\beta ds\ |\mathcal{F}_{t}]

+\mathbb{E}[\int_{t}^T  \sum_{i=1}^d|\hat{Y}_{i,s}|^2\times d^2\times n\times \varkappa(2\|Z_{s}\|)^2 ds\ |\mathcal{F}_{t}]\\

&&
+\mathbb{E}[\int_{t}^T (\sum_{i=1}^d|\hat{Y}_{i,s}|^2)^{1/2} 2\sqrt{d}\big(\hat{\upsilon}+\hat{\beta}\sum_{i'=1}^n|Y_{i',s}| 
+ \hat{\varkappa}(\|Z_{s}\|)\|Z_{s}\|\big)ds\ |\mathcal{F}_{t}].\ \ok
\dce
$$

\subsection{Differential inequalities}\label{differential-inequality}

Recall the functions $\Lambda_{t}$ and $\hat{\Lambda}_{t}$ defined in (\ref{the-lambdas}).

\bl\label{lambda-hat-lambda}
Suppose the conditions (\ref{bmo-boundedness-y-z}), (\ref{CR-C1}), (\ref{Df-slice}), and the conditions (\ref{package1}) and (\ref{explicit-Df}). Then, as function of the time $t$, $\Lambda_{t}$ and $\hat{\Lambda}_{t}$ are continuous on $[0,T]$ and they satisfy the differential inequalities:
$$
\dcb
\limsup_{\epsilon\downarrow 0}\sup_{0\leq t'< t\leq T, t-t'\leq \epsilon}\frac{1}{t-t'}\big(\Lambda_{t'}-\Lambda_{t}\big)

\leq
d \times \beta\times \Lambda_{t}
+d\times \upsilon+ d\times\varkappa(\hat{\Lambda}_{t})
\hat{\Lambda}_{t},
\dce
$$
and
$$
\dcb
\limsup_{\epsilon\downarrow 0}\sup_{0\leq t'< t\leq T, t-t'\leq \epsilon}\frac{1}{t-t'}\big(\hat{\Lambda}_{t'}^2-\hat{\Lambda}_{t}^2\big)
\\

\leq
\hat{\Lambda}_{t}^2\big(d\times 2\beta 
+
d^2\times n\times\varkappa(2\hat{\Lambda}_{t})^2\big)

+\hat{\Lambda}_{t}2\sqrt{d\times n}\big( \hat{\upsilon}+\hat{\beta}\Lambda_{t}+   
\hat{\varkappa}(\hat{\Lambda}_{t})\hat{\Lambda}_{t}\big).\\
\dce
$$
\el

\proof
Notice first of all that under the conditions of the lemma, by Proposition \ref{value-bound} and Corollary \ref{MD-Z-bound}, the functions $\Lambda_{t}$ and $\hat{\Lambda}_{t}$ are finite on $[0,T]$. Suppose that $0\leq t'< t\leq T$ and $t-t'\leq 1$. Suppose $\epsilon=t-t'$ satisfies  the condition in Lemma \ref{Phi-bounded-driver} with $c_{2}(d,\beta,\sqrt{d\times n}\times \eta,\epsilon)\leq \frac{1}{2}$. We begin the proof with the computation of $\Lambda_{t'}-\Lambda_{t}$ for $0\leq t'<t\leq T$. $$
\dcb
\Lambda_{t'}-\Lambda_{t}
&=&
\sup_{t'\leq a\leq t}\sum_{i=1}^d\|Y_{i,a}\|_{\infty}\vee \sup_{t\leq s\leq T}\sum_{i=1}^d\|Y_{i,s}\|_{\infty}
-
\sup_{t\leq s\leq T}\sum_{i=1}^d\|Y_{i,s}\|_{\infty}\\

&\leq&
\sup_{t'\leq a\leq t}(\sum_{i=1}^d\|Y_{i,a}\|_{\infty}
-
\sum_{i=1}^d\|Y_{i,t}\|_{\infty}).
\dce
$$
We consider the BSDE$[t,Y_{t},f]$(\ref{theBSDE}) (the terminal time being $t$ instead of $T$, the terminal value being $Y_{t}$ instead of $\xi$) and its Malliavin derivative BSDE$[t, Y_{t}, f, j,\theta]$(\ref{theBSDE-malliavin-*}). Apply Lemma \ref{Y-phi} with respect to the BSDE parameters $[t, Y_{t}, f]$. $$
\dcb
&&
\sup_{t'\leq a\leq t}\sum_{i=1}^d\|Y_{i,a}\|_{\infty} - \sum_{i=1}^d\|Y_{i,t}\|_{\infty}\\
&\leq&
e^{d \times \beta\times \epsilon}(\sum_{i=1}^d\|Y_{i,t}\|_{\infty}+d\times (\upsilon\times\epsilon+ \sup_{t' \leq a\leq t}\|\mathbb{E}[\int_{a}^t\varkappa(\|Z_{s}\|)\|Z_{s}\| ds\ |\mathcal{F}_{a}]\|_{\infty})
- 
\sum_{i=1}^d\|Y_{i,t}\|_{\infty}\\

&\leq&
(e^{d \times \beta\times \epsilon}-1)\Lambda_{t}

+e^{d \times \beta\times \epsilon}\times d\times \upsilon\times\epsilon+ e^{d \times \beta\times \epsilon}\times d\times\sup_{t' \leq a\leq t}\|\mathbb{E}[\int_{a}^t\varkappa(\|Z_{s}\|)\|Z_{s}\| ds\ |\mathcal{F}_{a}]\|_{\infty}.
\dce
$$
We apply the results of Section \ref{DY-lipschitzian}. Because $f$ is a Lipschitzian driver, the value process $Y$ (cf. Proposition \ref{value-bound}), the martingale coefficient process $Z$, as well as $DY$ (cf. Corollary \ref{MD-Z-bound}), are all bounded (uniformly with respect to the indices $j,\theta$). We consider the martingale $\hat{M}_{j,\theta,i,s}$ defined in Lemma \ref{DY-L-K} with respect to $D_{j,\theta}Y_{t}$, and its martingale coefficient process $\hat{\mu}_{j,\theta}$. We have
$$
\sup_{1\leq j\leq n}\sup_{t\leq \theta\leq T}\|\hat{\mu}_{j,\theta}\|_{\star}
\leq
\sup_{1\leq j\leq n}\sup_{t\leq \theta\leq T}\|D_{j,\theta}Y_{t}\|_{\infty}<\infty.
$$
By Lemma \ref{elkaroui-al}, $Z_{i,j,s}=D_{j,s}Y_{i,s}= \hat{Y}_{j,s,i,s}$ for $1\leq i\leq d,1\leq j\leq n, 0\leq s\leq t$. For all $0<t-s\leq \epsilon$, we deduce from Lemma \ref{DY-L-K} that
$$
\dcb
&&
\sup_{1\leq j\leq n}\sup_{1\leq i\leq d}\|Z_{i,j,s}- \hat{M}_{j,s,i,s}\|_{\infty}
=
\sup_{1\leq j\leq n}\sup_{1\leq i\leq d}\|D_{j,s}Y_{i,s}- \hat{M}_{j,s,i,s}\|_{\infty}\\
&\leq&
\sup_{1\leq j\leq n}\sup_{1\leq i\leq d}\sup_{0\leq \theta\leq s}\sum_{i=1}^d\|D_{j,\theta}Y_{i,s}-\hat{M}_{j,\theta,i,s}\|_{\infty}\\

&\leq&
e^{d\times \beta \times \epsilon}\ 
\big(1+d\times \sqrt{d\times n}\times \eta \times \sqrt{\epsilon}\ \frac{1}{1-c_{2}(d,\beta,\sqrt{d\times n}\times \eta,\epsilon)}\
e^{d\times\beta\times \epsilon}\sqrt{2 \big(e^{-d\times\beta\times \epsilon}
+
\beta\times\epsilon\big)}\big) \times d\times\\
&& 
\big(
\beta\times \sup_{1\leq j\leq n}\sup_{0\leq \theta\leq T}\sum_{i'=1}^d\|D_{j,\theta}Y_{i',t}\|_{\infty} \\
&&
+ \hat{\upsilon} +\hat{\beta}\times\sup_{0\leq s\leq T}\sum_{i'=1}^d\|Y_{i',s}\|_{\infty} 
+\hat{\varkappa}(\sup_{0\leq s\leq T}\|Z_{s}\|)\sup_{0\leq s\leq T}\|Z_{s}\|
\big)\times \epsilon\\
&&
+e^{d\times \beta \times \epsilon}\ d\times\sqrt{d\times n}\times  \eta \times \sqrt{\epsilon}\ \frac{1}{1-c_{2}(d,\beta,\sqrt{d\times n}\times \eta,\epsilon)}
\sup_{1\leq j\leq n}\sup_{0\leq \theta\leq T}\|D_{j,\theta}Y_{t}\|_{\infty}
\dce
$$ 
Notice that $$
\dcb
&&\sum_{j=1}^n\sum_{i=1}^d|\hat{M}_{j,s,i,s}|^2
=
\sum_{j=1}^n\sum_{i=1}^d|\mathbb{E}[D_{j,s}Y_{i,t}|\mathcal{F}_{s}]|^2

\leq
\sum_{j=1}^n\sum_{i=1}^d\mathbb{E}[(D_{j,s}Y_{i,t})^2|\mathcal{F}_{s}]\\
&&
=
\mathbb{E}[\sum_{j=1}^n\sum_{i=1}^d(D_{j,s}Y_{i,t})^2|\mathcal{F}_{s}]
\leq
\|\sum_{j=1}^n\sum_{i=1}^d(D_{j,s}Y_{i,t})^2\|_{\infty}
\leq
\sup_{0\leq \theta<\infty}\|D_{\theta}Y_{t}\|_{\infty}^2.

\dce
$$
Recall that $\varkappa$ is an increasing locally Lipschitzian function. There exists a constant $L$ which depend only on the BSDE parameters $\beta, \hat{\beta}, \upsilon, \eta$, the dimension parameters $d,n$, the Lipschitzian indices of $\varkappa, \hat{\varkappa}$, and$$
\dcb
\sup_{1\leq j\leq n}\sup_{0\leq \theta\leq t\leq T}\sum_{i=1}^d\|D_{j,\theta}Y_{i,t}\|_{\infty} +\sup_{0\leq s\leq T}\sum_{i=1}^d\|Y_{i,s}\|_{\infty} 
+\sup_{0\leq s\leq T}\|Z_{s}\|,
\dce
$$
such that (recalling that $t-t'=\epsilon$ and $c_{2}(d,\beta,\sqrt{d\times n}\times \eta,\epsilon)\leq \frac{1}{2}$)
$$
\dcb
&&
\frac{1}{t-t'}\sup_{t' \leq a\leq t}\|\mathbb{E}[\int_{a}^t\varkappa(\|Z_{s}\|)\|Z_{s}\| \ ds\ |\mathcal{F}_{a}]\|_{\infty}\\

&\leq&
\frac{1}{t-t'}\sup_{t' \leq a\leq t}\|\mathbb{E}[\int_{a}^t|\varkappa(\|Z_{s}\|)\|Z_{s}\|-\varkappa(\|\hat{M}_{s,s}\|)\|\hat{M}_{s,s}\|| \ ds\ |\mathcal{F}_{a}]\|_{\infty}\\
&&
+\frac{1}{t-t'}\sup_{t' \leq a\leq t}\|\mathbb{E}[\int_{a}^t\varkappa(\sqrt{\sum_{j=1}^n\sum_{i=1}^d|\hat{M}_{j,s,i,s}|^2})\sqrt{\sum_{j=1}^n\sum_{i=1}^d|\hat{M}_{j,s,i,s}|^2} \ ds\ |\mathcal{F}_{a}]\|_{\infty}\\

&\leq&
L\times\sqrt{t-t'}
+
\frac{1}{t-t'}\sup_{t' \leq a\leq t}\|\mathbb{E}[\int_{a}^t\varkappa(\sup_{0\leq \theta<\infty}\|D_{\theta}Y_{t}\|_{\infty})
\sup_{0\leq \theta<\infty}\|D_{\theta}Y_{t}\|_{\infty} \ ds\ |\mathcal{F}_{a}]\|_{\infty}\\

&=&
L\times\sqrt{t-t'}
+\varkappa(\sup_{0\leq \theta<\infty}\|D_{\theta}Y_{t}\|_{\infty})
\sup_{0\leq \theta<\infty}\|D_{\theta}Y_{t}\|_{\infty}\ 

\leq
L\times\sqrt{t-t'}
+\varkappa(\hat{\Lambda}_{t})
\hat{\Lambda}_{t}.
\dce
$$
We conclude that
$$
\dcb
&&
\frac{1}{t-t'}\big(\Lambda_{t'}-\Lambda_{t}\big)
\leq
\frac{1}{t-t'}\big(\sup_{t'\leq a\leq t}\sum_{i=1}^d\|Y_{i,a}\|_{\infty} - \sum_{i=1}^d\|Y_{i,t}\|_{\infty}\big)\\

&\leq&
\frac{e^{d \times \beta\times (t-t')}-1}{t-t'} \Lambda_{t}
+e^{d\times\beta\times(t-t')}\times d\times \upsilon

+ e^{d\times\beta\times(t-t')}\times d\times 
\big(L\times\sqrt{t-t'}
+\varkappa(\hat{\Lambda}_{t})
\hat{\Lambda}_{t}\big).
\dce
$$
This proves the continuity of $\Lambda$ and the first differential inequality of the lemma. Consider $\hat{\Lambda}_{t}^2$.$$
\dcb
\hat{\Lambda}_{t'}^2-\hat{\Lambda}_{t}^2

=
\sup_{0\leq \theta<\infty}\big( (\sup_{t'\leq a\leq t}\|D_{\theta}Y_{a}\|_{\infty}^2- \sup_{t\leq s\leq T}\|D_{\theta}Y_{s}\|_{\infty}^2)^+
+
\sup_{t\leq s<\infty}\|D_{\theta}Y_{s}\|_{\infty}^2\big)
-\hat{\Lambda}_{t}\\

\leq
\sup_{0\leq \theta<\infty} (\sup_{t'\leq a\leq t}\|D_{\theta}Y_{a}\|_{\infty}^2- \sup_{t\leq s\leq T}\|D_{\theta}Y_{s}\|_{\infty}^2)^+\\

\leq
\sup_{0\leq \theta<\infty} (\sup_{\theta\vee t'\leq a\leq t}\|D_{\theta}Y_{a}\|_{\infty}^2- \|D_{\theta}Y_{t}\|_{\infty}^2)^+

=
\sup_{0\leq \theta<\infty} \sup_{\theta\vee t'\leq a\leq t}(\|D_{\theta}Y_{a}\|_{\infty}^2- \|D_{\theta}Y_{t}\|_{\infty}^2)^+.

\dce
$$
By Lemma \ref{sumY2} applied to the Malliavin derivative of the BSDE$[t,Y_{t},f]$(\ref{theBSDE}), for $0\leq \theta\leq a\leq t$,
$$
\dcb
&&
\sum_{j=1}^n\sum_{i=1}^d(D_{j,\theta}Y_{i,a})^2\\

&\leq&
\mathbb{E}[\sum_{j=1}^n\sum_{i=1}^d(D_{j,\theta}Y_{i,t})^2|\mathcal{F}_{a}]\\
&&
+\mathbb{E}[\int_{a}^t  \sum_{j=1}^n\sum_{i=1}^d(D_{j,\theta}Y_{i,s})^2\ \times 
\big(
d\times 2\beta 
+
d^2\times n\times \varkappa(2\|Z_{s}\|)^2\big)
ds\ |\mathcal{F}_{a}]\\

&&
+\mathbb{E}[\int_{a}^t \sum_{j=1}^n(\sum_{i=1}^d(D_{j,\theta}Y_{i,s})^2)^{1/2} 2\sqrt{d}\big(\hat{\upsilon}+\hat{\beta}\sum_{i'=1}^n|Y_{i',s}| 
+ \hat{\varkappa}(\|Z_{s}\|)\|Z_{s}\|\big)ds\ |\mathcal{F}_{a}].
\dce
$$
Note that $\sum_{j=1}^n(\sum_{i=1}^d(D_{j,\theta}Y_{i,s})^2)^{1/2}\leq (\sum_{j=1}^n\sum_{i=1}^d(D_{j,\theta}Y_{i,s})^2)^{1/2}\sqrt{n}$. Hence,$$
\dcb
&&
(\|D_{\theta}Y_{a}\|_{\infty}^2- \|D_{\theta}Y_{t}\|_{\infty}^2)^+\\

&\leq&
\|\ \mathbb{E}[\int_{a}^t  \sum_{j=1}^n\sum_{i=1}^d(D_{j,\theta}Y_{i,s})^2\ \times 
\big(
d\times 2\beta 
+
d^2\times n\times \varkappa(2\|Z_{s}\|)^2\big)
ds\ |\mathcal{F}_{a}]\ \|_{\infty}\\

&&
+\|\ \mathbb{E}[\int_{a}^t (\sum_{j=1}^n\sum_{i=1}^d(D_{j,\theta}Y_{i,s})^2)^{1/2} 2\sqrt{d\times n}\big(\hat{\upsilon}+\hat{\beta}\sum_{i'=1}^n|Y_{i',s}| 
+ \hat{\varkappa}(\|Z_{s}\|)\|Z_{s}\|\big)ds\ |\mathcal{F}_{a}]\ \|_{\infty}\\

&\leq&
\hat{\Lambda}_{a}^2\big(d\times 2\beta\times (t-a) 
+
d^2\times n\times\|\ \mathbb{E}[\int_{a}^t 
 \varkappa(2\|Z_{s}\|)^2
ds\ |\mathcal{F}_{a}]\ \|_{\infty}\big)\\

&&
+\hat{\Lambda}_{a}2\sqrt{d\times n}\big( (\hat{\upsilon}+\hat{\beta}\Lambda_{a})(t-a)+\|\ \mathbb{E}[\int_{a}^t   
\hat{\varkappa}(\|Z_{s}\|)\|Z_{s}\|ds\ |\mathcal{F}_{a}]\ \|_{\infty}\big).
\dce
$$
By the same argument as that we do for $\Lambda_{t}$, there is a constant $\hat{L}$, independent of $t,t'$, such that
$$
\dcb
\frac{1}{t-t'}\big(\hat{\Lambda}_{t'}^2-\hat{\Lambda}_{t}^2\big)

&\leq&
\hat{L}\hat{\Lambda}_{t'}^2\times\sqrt{t-t'}

+\hat{\Lambda}_{t'}^2\big(d\times 2\beta 
+
d^2\times n\times\varkappa(2\hat{\Lambda}_{t})^2\big)\\

&&
+\hat{\Lambda}_{t'}2\sqrt{d\times n}\big( \hat{\upsilon}+\hat{\beta}\Lambda_{t'}+   
\hat{\varkappa}(\hat{\Lambda}_{t})\hat{\Lambda}_{t}\big).
\dce
$$
This proves the continuity of $\hat{\Lambda}^2$ and then the second inequality of the lemma. \ok

\

\section{Notes on the Lyapunov function and on the sliceability}\label{L-S}

We collect some immediate consequences of the Lyapunov function and of the sliceability.

\subsection{Lyapunov function and its consequences}\label{about-LF}

We begin with the following lemma.

\bl
Let $\mathtt{h}$ be a smooth positive uniformly strictly convex function. Let $f$ be a driver such that $$
|\partial_{y_i}\mathtt{h}(y)f(i, t, k^{-1}y, z)| \leq C(1+\|z\|^2)
$$
for all $t, y, z$ and for $k>0$ big enough. Then, for some $k>0$, $\mathtt{h}(ky)$ defines a Lyapunov function for $f$.
\el

As an example of such a driver, we can take$$
f(i, t, y, z) = e^{-|\partial_{y_i}\mathtt{h}(y)|}(1+\cos(|y|)+\|z\|^{1+\alpha}),
$$
for $0\leq \alpha\leq 1$.

\proof
For a given $k>0$, we compute
$$
\dcb
&&
\frac{1}{2}\sum_{i=1}^d\sum_{i'=1}^d k^2\partial_{y_{i}}\partial_{y_{i'}}\mathtt{h}(ky)\sum_{j=1}^nz_{i,j}z_{i',j}
-
\sum_{i=1}^dk\partial_{y_{i}}\mathtt{h}(ky)f(i,t, y, z)\\

&=&
k^2\frac{1}{2}\sum_{j=1}^n\sum_{i=1}^d\sum_{i'=1}^d \partial_{y_{i}}\partial_{y_{i'}}\mathtt{h}(ky)z_{i,j}z_{i',j}
-
k\sum_{i=1}^d\partial_{y_{i}}\mathtt{h}(ky)f(i,t, y, z)\\

&\geq& 
k^2\frac{1}{2}\sum_{j=1}^n c \sum_{i=1}^d|z_{i,j}|^2
-
kC(1+\|z\|^2)\\
&&\mbox{where $c>0$ is a constant of the uniform convexity of $\mathtt{h}$,}\\

&=&
k^2c \ \frac{1}{2}\|z\|^2
-
kC(1+\|z\|^2)
\geq 
\|z\|^2+kC,\\

\dce
$$
for large enough $k>0$. \ok

\

Let us consider the consequences of a Lyapunov function. The basic reason to introduce Lyapunov function is to control the martingale coefficient process $Z$ of the BSDE$[T,\xi, f]$(\ref{theBSDE}) by the value process $Y$.

\bl\label{Z2-bound}
Consider the BSDE$[T,\xi, f]$(\ref{theBSDE}). Suppose that the driver $f$ has a global Lyapunov function $\mathtt{h}$ (together with $\mathtt{k}$ the lower bound function). Then, for any $0\leq b\leq T$, for any solution $(Y,Z)$ on the time interval $[b,T]$ of the BSDE$[T, \xi, f]$(\ref{theBSDE}) such that $Y\in\mathcal{S}_{\infty}[b,T]$ (and therefore $\sum_{i=1}^d\|\xi_{i}\|_{\infty}<\infty$), for any pair of stopping times $b\leq \sigma\leq \sigma'\leq T$, we have$$
\dcb
\mathbb{E}[\int_{\sigma}^{\sigma'}\|Z_{s}\|^2\ ds|\mathcal{F}_{\sigma}]
\leq
\mathbb{E}[\mathtt{h}(Y_{\sigma'}) -\mathtt{h}(Y_{\sigma}) + \int_{\sigma}^{\sigma'} \mathtt{k}(s, Y_{s})ds |\mathcal{F}_{\sigma}].
\dce
$$
\el

\proof
By Ito's formula,
$$
\dcb
\mathtt{h}(Y_{\sigma'}) - \mathtt{h}(Y_{\sigma})

&=&
\int_{\sigma}^{\sigma'}\sum_{i=1}^d\partial_{y_{i}}\mathtt{h}(Y_{s})dY_{i,s} + \frac{1}{2}\int_{\sigma}^{\sigma'}\sum_{i=1}^d\sum_{i'=1}^d\partial_{y_{i}}\partial_{y_{i'}}\mathtt{h}(Y_{s})d\cro{Y_{i},Y_{i'}}_{s}\\

&=&
-\int_{\sigma}^{\sigma'}\sum_{i=1}^d\partial_{y_{i}}\mathtt{h}(Y_{s})f(i,s,Y_{s},Z_{s})ds

+ \int_{\sigma}^{\sigma'}\sum_{i=1}^d\partial_{y_{i}}\mathtt{h}(Y_{s})\sum_{j=1}^n Z_{i,j}dB_{j,s}\\
&&
+\frac{1}{2}\int_{\sigma}^{\sigma'}\sum_{i=1}^d\sum_{i'=1}^d\partial_{y_{i}}\partial_{y_{i'}}\mathtt{h}(Y_{s})\sum_{j=1}^n Z_{i,j}Z_{i',j}ds\\

&\geq&
\int_{\sigma}^{\sigma'}\|Z_{s}\|^2 - \mathtt{k}(s, Y_{s})ds

+ \int_{\sigma}^{\sigma'}\sum_{i=1}^d\partial_{y_{i}}\mathtt{h}(Y_{s})\sum_{j=1}^n Z_{i,j}dB_{j,s}.
\dce
$$
Suppose firstly that the stopping times $\sigma, \sigma'$ are chosen in such a way that every above stochastic integrals are uniformly integrable. We have 
$$
\dcb
\mathbb{E}[\int_{\sigma}^{\sigma'}\|Z_{s}\|^2\ ds|\mathcal{F}_{\sigma}]
\leq
\mathbb{E}[\mathtt{h}(Y_{\sigma'}) -\mathtt{h}(Y_{\sigma})+
\int_{\sigma}^{\sigma'} \mathtt{k}(s, Y_{s})ds\ |\mathcal{F}_{\sigma}].
\dce
$$
As $\|Z_{s}\|^2$ and $\mathtt{k}(s, Y_{s})$ are positive, as $\mathtt{h}(Y_{s})$ is uniformly bounded, we can pass to limit to establish the inequality for a general couple of stopping times $\sigma,\sigma'$. \ok 

\

The upper bound on the value process $Y$ in Proposition \ref{value-bound} is no more valid, if the driver of the BSDE is not Lipschitzian in $z$ (so that $\eta=\infty$). Lyapunov function allows us to write a new upper bound on the value process $Y$.

\bcor\label{bound-Lyapunov}
Suppose the same assumptions on BSDE$[T,\xi, f]$(\ref{theBSDE}) as in Lemma \ref{Z2-bound}. Suppose moreover that
\ebe
\item
the driver $f$ is uniformly Lipschitzian in $y$ with index $\beta>0$,
\item
$f$ has an increment rate in $z$ uniformly less than $\varkappa$ and the rate function $\varkappa$ satisfies the inequality: $\varkappa(r)r\leq \gamma(1+r^2), r\geq 0,$ for some constant $\gamma>0$,
\item
$\|f\|_{\centerdot}<\infty$ (the potential bound). 
\dbe
Then, for any $0\leq b\leq T$, for any solution $(Y,Z)$ on $[b,T]$ of the BSDE$[T, \xi, f]$(\ref{theBSDE}) such that $Y\in\mathcal{S}_{\infty}[b,T]$, we have
$$
\dcb
\sup_{b\leq  t\leq T}\|\sum_{i=1}^d|Y_{i,t}|\ \|_{\infty}
\leq
e^{d\times (\beta+\gamma \bar{\beta})\times T}(\sum_{i=1}^d\|\xi_{i}\|_{\infty}+d\times (\|f\|_{\centerdot}+\gamma\times T+
\gamma\|\mathtt{h}(\xi)\|_{\infty}+\gamma\|\mathtt{k}\|_{\centerdot})),
\dce
$$
and$$
\dcb
\|\ind_{[b,T]}Z\|_{\star}^2
\leq
\|\mathtt{h}(\xi)\|_{\infty}+\|\mathtt{k}\|_{\centerdot}+  \bar{\beta}\times T \sup_{b\leq a\leq T}\|\sum_{i'=1}|Y_{i',a}|\|_{\infty}.
\dce
$$
\ecor

\

One should read the corollary in the following way: either a solution $(Y,Z)$ is unbounded, or it is bounded, indifferently with respect to $b$ and to $(Y,Z)$, by a same constant. No intermediate situation is possible. 

\

\proof
Starting from the formula in Lemma \ref{Z2-bound}, under the conditions of the corollary, for $b\leq t\leq T$, we have (recalling $\mathtt{h}\geq 0$ and the conditions (\ref{k-beta}))
$$
\dcb
\mathbb{E}[\int_{t}^{T}\|Z_{s}\|^2\ ds|\mathcal{F}_{t}]
&\leq&
\mathbb{E}[\mathtt{h}(Y_{T})+
\int_{t}^{T} \mathtt{k}(s, Y_{s})ds\ |\mathcal{F}_{t}]

\leq
\|\mathtt{h}(\xi)\|_{\infty}+ \mathbb{E}[\int_{t}^{T} \mathtt{k}(s, 0) + \bar{\beta}\sum_{i'=1}^d|Y_{i',s}| ds |\mathcal{F}_{t}]\\

&\leq&
\|\mathtt{h}(\xi)\|_{\infty}+\|\mathtt{k}\|_{\centerdot}+\bar{\beta}\times T \sup_{b\leq a\leq T}\|\sum_{i'=1}|Y_{i',a}|\|_{\infty}.
\dce
$$
This proves the second inequality of the corollary. Consider now the first inequality.
By formula (\ref{ine-y}) (which is applicable on the time interval $[b,T]$ under the assumptions of the corollary), by Lemma \ref{Z2-bound}, for $1\leq i\leq d$, $b\leq t\leq T,$
$$
\dcb
&&|Y_{i,t}|
\leq 
\mathbb{E}[|\xi_{i}|+\int_{t}^T |f(i,s, Y_{s},Z_{t})| ds|\mathcal{F}_{t}]\\

&\leq&
\mathbb{E}[|\xi_{i}|+\int_{t}^T |f(i,s, 0,0)| ds + \int_{t}^T \beta \sum_{i'=1}^d|Y_{i',s}| ds + \int_{t}^T \varkappa(\|Z_{s}\|)\|Z_{s}\| ds|\mathcal{F}_{t}]\\

&\leq&
\|\xi_{i}\|_{\infty}+\|f\|_{\centerdot}+\mathbb{E}[\int_{t}^T \beta \sum_{i'=1}^d|Y_{i',s}| ds + \int_{t}^T \gamma(1+\|Z_{s}\|^2) ds|\mathcal{F}_{t}]\\

&\leq&
\|\xi_{i}\|_{\infty}+\|f\|_{\centerdot}+\mathbb{E}[\int_{t}^T \beta \sum_{i'=1}^d|Y_{i',s}| ds|\mathcal{F}_{t}]+\gamma(T-t)
+
\gamma\mathbb{E}[\mathtt{h}(Y_{T}) + \int_{t}^{T} \mathtt{k}(s, Y_{s}) ds |\mathcal{F}_{t}]\\

&\leq&
\|\xi_{i}\|_{\infty}+\|f\|_{\centerdot}+\gamma\times T+\mathbb{E}[\int_{t}^T \beta \sum_{i'=1}^d|Y_{i',s}| ds|\mathcal{F}_{t}]
+
\gamma\|\mathtt{h}(\xi)\|_{\infty}+\gamma \mathbb{E}[\int_{t}^{T} \mathtt{k}(s, 0) + \bar{\beta}\sum_{i'=1}^d|Y_{i',s}| ds |\mathcal{F}_{t}]\\

&\leq&
\|\xi_{i}\|_{\infty}+\|f\|_{\centerdot}+\gamma\times T+
\gamma\|\mathtt{h}(\xi)\|_{\infty}+\gamma\|\mathtt{k}\|_{\centerdot}+\mathbb{E}[\int_{t}^T \beta \sum_{i'=1}^d|Y_{i',s}| ds|\mathcal{F}_{t}]
+\mathbb{E}[\int_{t}^{T} \gamma \bar{\beta}\sum_{i'=1}^d|Y_{i',s}| ds |\mathcal{F}_{t}],
\dce
$$
i.e.,
$$
\dcb
\sum_{i=1}^d|Y_{i,t}|
&\leq&
\sum_{i=1}^d\|\xi_{i}\|_{\infty}+d\times (\|f\|_{\centerdot}+\gamma\times T+
\gamma\|\mathtt{h}(\xi)\|_{\infty}+\gamma\|\mathtt{k}\|_{\centerdot})\\
&&
+\mathbb{E}[\int_{t}^T d\times (\beta+\gamma \bar{\beta}) \sum_{i'=1}^d|Y_{i',s}| ds|\mathcal{F}_{t}].
\dce
$$
By Gronwall's inequality Lemma \ref{gronwall-ineq}, 
$$
\dcb
\sum_{i=1}^d|Y_{i,t}|
&\leq&
e^{d\times (\beta+\gamma \bar{\beta}) (T-t)}(\sum_{i=1}^d\|\xi_{i}\|_{\infty}+d\times (\|f\|_{\centerdot}+\gamma\times T+
\gamma\|\mathtt{h}(\xi)\|_{\infty}+\gamma\|\mathtt{k}\|_{\centerdot})).\ \ok
\dce
$$

\subsection{Exponential integrability via uniform sliceability}

\bl\label{slice-bound}
Consider $0\leq a\leq T$ and a predictable matrix valued process $X$, which is uniformly sliceable on $(a,T]$ (an interval open at $a$).
Let $\delta>0$ be a constant which makes $\|\ind_{(b,b')}.X\|_\star^2\leq \frac{1}{2}$ for any $a< b<b'\leq T, b'-b\leq \delta$. Let $N$ be an integer number bigger than $\frac{T-a}{\delta}$. Then, $$
\dcb
\mathbb{E}[e^{\int_{a}^T\|X_{s}\|^{2}ds}|\mathcal{F}_{a}] \leq 2^{N}.
\dce
$$
\el

\proof
Let $a< b\leq b'\leq T$ with $b'-b\leq \delta$. By John-Nirenberg inequality \cite[Chapitre VI, $n^\circ$105, Theorem]{DM}) or \cite[Theorem 2.2]{K},$$
\dcb
\mathbb{E}[e^{\int_{b}^{b'}\|X_{s}\|^{2}ds} |\mathcal{F}_{b}]
\leq
\frac{1}{1-\|\ind_{(b,b')}.X\|^{2}_{\star}}
\leq 2.
\dce
$$
Hence, if $T-a\leq \delta$, 
$$
\dcb
&&\mathbb{E}[e^{\int_{a}^T\|X_{s}\|^{2}ds} |\mathcal{F}_{a}]\leq 2.
\dce
$$
Otherwise, 
$$
\dcb
&&\mathbb{E}[e^{\int_{a}^T\|X_{s}\|^{2}ds} |\mathcal{F}_{a}]
=
\mathbb{E}[e^{\int_{a}^{T-\delta}\|X_{s}\|^{2}ds} \mathbb{E}[e^{\int_{T-\delta}^T\|X_{s}\|^{2}ds} |\mathcal{F}_{T-\delta}]\ |\mathcal{F}_{a}]\\
&\leq&
2\mathbb{E}[e^{\int_{a}^{T-\delta}\|X_{s}\|^{2}ds} |\mathcal{F}_{a}]\leq \ldots \leq 2^{N}. \ \ok
\dce
$$

\


\begin{thebibliography}{99}
\bibitem{bahlali}
\textsc{Bahlali K.} (2002) \texttt{"}Existence and uniqueness of solutions for BSDEs with locally Lipschitz coefficient\texttt{"} \textit{Elect. Comm. in Probability} \textbf{7} 169-179 

\bibitem{bahlali-al}
\textsc{Bahlali, K., Essaky, E.H., Hassani, M.} (2010) \texttt{"}p-integrable solutions to multidimensional BSDEs and degenerate systems of PDEs with logarithmic nonlinearities\texttt{"} \textit{arXiv}1007.2388


\bibitem{BEK}
\textsc{Barrieu P. and El Karoui N.} (2013) \texttt{"}Monotone stability of quadratic semimartingales with applications to unbounded general quadratic BSDEs\texttt{"} \textit{The Annals of Probability} \textbf{41}(3B) 1831-1863 

\bibitem{BH}
\textsc{Briand P. and Hu Y.} (2006) \texttt{"}BSDE with quadratic growth and unbounded terminal value\texttt{"} \textit{Probability Theory and Related Fields} \textbf{136} 604-618 


\bibitem{CN}
\textsc{Cheridito P. and Nam K.} (2015) \texttt{"}Multidimensional quadratic and subquadratic BSDEs with special structure\texttt{"} \textit{arXiv}1309.6716v4

\bibitem{DM}
\textsc{Dellacherie C. and Meyer A.P.} (1975-1992) \textit{Probabilit\' es et Potentiel} Hermann 


\bibitem{DHB}
\textsc{Delbaen F. and Hu Y. and Bao X} (2011) \texttt{"}Backward SDEs with superquadratic growth\texttt{"} \textit{Probability Theory and Related Fields} \textbf{150}(1-2) 145-192


\bibitem{EKPQ}
\textsc{El Karoui N. and Peng S. and Quenez M.C.} (1997) \texttt{"}Backward stochastic differential equations in finance\texttt{"}\textit{Math. Finance} \textbf{7}(1) 1?71
 

\bibitem{FR}
\textsc{Frei C. and Reis D.} (2011) \texttt{"}A financial market with interacting investors: does an equilibrium exist?\texttt{"} \textit{Math. Financ. Econ.} \textbf{4} 161-182


\bibitem{HR}
\textsc{Harter J. and Richou A.} (2017) \texttt{"}A stability approach for solving multidimensional quadratic BSDEs\texttt{"} \textit{arXiv}1606.08627v2


\bibitem{HT}
\textsc{Hu Y. and Tang S.} (2016) \texttt{"}Multi-dimensional backward stochastic differential equations of diagonally quadratic generators\texttt{"} \textit{Stochastic Processes and their Applications} \textbf{126}(4) 1066-1086


\bibitem{HWY}
\textsc{He, S.W. and Wang, J.G. and Yan, J.A.} (1992) \textit{Semimartingale Theory and Stochastic Calculus} Science Press and CRC Press, Beijing 



\bibitem{IMPR}
\textsc{Imkeller P. and Mastroli T. and Possamai D. and Reveillac A.} (2016)  \texttt{"}A note on the Malliavin-Sobolev spaces\texttt{"} \textit{Statistics and Probability Letters} \textbf{109} 45-53



\bibitem{JKL}
\textsc{Jamneshan A. and Kupper M. and Luo P.} (2016) \texttt{"}Solvability of multidimensional quadratic BSDEs\texttt{"} \textit{arXiv}1612.02698v1


\bibitem{K}
\textsc{Kazamaki N.} (1994) \textit{Continuous exponential martingales and BMO} Lecture Notes in Mathematics \textbf{1579} Springer-Verlag

\bibitem{MW}
\textsc{Mocha M. and Westray N.} (2011) \texttt{"}Quadratic semimartingale BSDEs under an exponential moments condition\texttt{"} \textit{arXiv}1101.2582v1


\bibitem{N}
\textsc{Nualart D.} (2006) \textit{The Malliavin calculus and related topics} Springer-Verlag, Berlin


\bibitem{PR} 
\textsc{Pardoux E. and Rascanu A.} (2014) \textit{Stochastic Differential Equations, Backward SDEs, Partial Differential Equations} Springer


\bibitem{RY}
\textsc{Revuz D. and Yor M.} (1999) \textit{Continuous martingales and Brownian motion} Springer 

\bibitem{schachermayer}
\textsc{Schachermayer W.} (1996) \texttt{"}A characterisation of the closure of $H^\infty$ in BMO\texttt{"} \textit{S\'eminaire de Probabilit\'es de Strasbourg} \textbf{XXX} 344-356 

\bibitem{T}
\textsc{Tevzadze R.}  (2008) \texttt{"}Solvability of backward stochastic differential equations with quadratic growth\texttt{"} \textit{Stochastic Processes and their Applications} \textbf{118}(3) 503-515

\bibitem{xing}
\textsc{Xing H.} (2016) A talk at
\textsf{Stochastic Analysis and Mathematical Finance-A Fruitful Partnership, Oaxaca, Mexico, May 25}

\bibitem{XZ}
\textsc{Xing H. and Zitkovic G.} (2017) \texttt{"}A class of globally solvable Markovian quadratic BSDE systems and applications\texttt{"} \textit{arXiv}1603.00217v2


\end{thebibliography}
\end{document}